% Ceci est le fichier PlainTeX de l'article
% "Groupes de Garside", precedemment initule
% "Petits groupes gaussiens"; cette version remplace
% celle transmise en mars 2001

\magnification1000
\tolerance 3000
\hsize = 125mm
\vsize= 200mm
\hoffset-0truecm
\voffset-1truecm
\overfullrule=0pt
\parindent=10pt

\font\sc=cmr10
\font\go=eufm10

\catcode`\ç=\active\defç{\c c}        % option c
\catcode`\ä=\active\defä{\"a}        % option u, then  a
\catcode`\ë=\active\defë{\"e}        % option u, then  e
\catcode`\ï=\active\defï{\"{\i}}     % option u, then  i
\catcode`\ö=\active\defö{\"o}        % option u, then  o
\catcode`\ü=\active\defü{\"u}        % option u, then  u
\catcode`\ÿ=\active\defÿ{\"y}        % option u, then  y
\catcode`\é=\active\defé{\'e}        % option e, then  e
\catcode`\É=\active\defÉ{\'E}        % option e, then  E
\catcode`\à=\active\defà{\`a}        % option `, then  a
\catcode`\è=\active\defè{\`e}        % option `, then  e
\catcode`\ù=\active\defù{\`u}        % option `, then  u
\catcode`\À=\active\defÀ{\`A}        % option `, then  A
\catcode`\â=\active\defâ{\^a}        % option i, then  a
\catcode`\ê=\active\defê{\^e}        % option i, then  e
\catcode`\î=\active\defî{\^{\i}}     % option i, then  i
\catcode`\ô=\active\defô{\^o}        % option i, then  o
\catcode`\û=\active\defû{\^u}        % option i, then  u
\catcode`\¦=\active\def¦{\oe }        % option i, then  u

\catcode`\:=\active\catcode`\!=\active\catcode`\?=\active
\catcode`\:=\active\def:{\penalty500\hskip2pt\string:}
\catcode`\!=\active\def!{\penalty500\hskip2pt\string!}
\catcode`\?=\active\def?{\penalty500\hskip2pt\string?}
\frenchspacing

\font\cyrIX=wncyr9
\def\<<{\penalty-500{\cyrIX\char'074}\penalty500}
\def\>>{\penalty500{\cyrIX\char'076\ }\penalty-500}

\def\resp{{\it resp. }}
\def\Block{\hbox{\vrule width 7pt height 7pt depth 0pt}}

\newcount\SecNo
\SecNo=0
\newcount\PropNo
\newcount\EqNo

% COMMENT: \Displaylines and \Eqalign are similar to TeX's
% standard \displaylines and \eqalign, excepted that each line
% can receive an equation number, by using a second & sign.
% The boolean \Displaylines is true inside a \Displaylines or
% \Eqalign group.

\newif\ifDisplaylines\Displaylinesfalse
{\catcode`\@=11
\gdef\Displaylines#1{\Displaylinestrue\displ@y\halign{
\hbox to\displaywidth{$\@lign\hfil\displaystyle##\hfil$}%
&\llap{$##$}\crcr#1\crcr}}
\gdef\Eqalign#1{\Displaylinestrue\displ@y \tabskip=\centering
\halign to \displaywidth{\hfil$\@lign\displaystyle{##}$\tabskip=0pt
&$\@lign\displaystyle{{}##}$\hfil\tabskip=\centering
&\llap{$\@lign##$}\tabskip=0pt\crcr#1\crcr}}}

% COMMENT: cross-reference system. To give give a name or a
% number to a formula or a proposition. The command \«54»
% causes the content of \box0 to be stored in \box54, and then
% «54» lets this content be printed. Similarly \«\name» does 
% the same with a box named \name, and «\name» lets it be
% printed. These are the only possible parameters, namely a
% number (between 1 and 255, preferably between 100 and 255),
% or a command name, i.e., \ followed by letters.

\catcode`\«=\active
{\catcode`\@=11
\gdef\«#1#2»{\ifcat\noexpand#1\noexpand\%
\global\alloc@4\box\chardef\insc@unt{#1#2}\fi%
\global\setbox#1#2\box0\ignorespaces}}
\def«#1»#2{\leavevmode\copy#1\ifcat#2a\ \fi#2}

\def\Sec#1\par{\vskip0pt plus.2\vsize
\penalty-250\vskip0pt plus-.2\vsize 
\advance\SecNo by 1\global\PropNo=0\global\EqNo=0
\vskip30pt\centerline{\sc\the\SecNo. #1}}

\def\Parag#1\par{\vskip0pt plus.1\vsize
\penalty-250\vskip0pt plus -.1\vsize
\bigbreak{\bf #1}}

\def\Def{\global\advance\PropNo by 1%
\setbox0=\hbox{\the\SecNo.\the\PropNo}%
\bigbreak{\bf Définition \the\SecNo.\the\PropNo.} \ignorespaces}

\def\Rem{\bigbreak{\bf Remarque. }}

\def\Conj{\global\advance\PropNo by 1%
\setbox0=\hbox{\the\SecNo.\the\PropNo}%
\bigbreak{\bf Conjecture \the\SecNo.\the\PropNo.} \ignorespaces}

\def\Question{\global\advance\PropNo by 1%
\setbox0=\hbox{\the\SecNo.\the\PropNo}%
\bigbreak{\bf Question \the\SecNo.\the\PropNo.} \ignorespaces}

\def\Cor{\global\advance\PropNo by 1%
\setbox0=\hbox{\the\SecNo.\the\PropNo}%
\bigbreak{\bf Corollaire \the\SecNo.\the\PropNo.} \ignorespaces}

\def\Ex{\global\advance\PropNo by 1%
\setbox0=\hbox{\the\SecNo.\the\PropNo}%
\bigbreak{\bf Exemple \the\SecNo.\the\PropNo.} \ignorespaces}

\def\Lem{\global\advance\PropNo by 1%
\setbox0=\hbox{\the\SecNo.\the\PropNo}%
\bigbreak{\bf Lemme \the\SecNo.\the\PropNo.} \ignorespaces}

\def\Prop{\global\advance\PropNo by 1%
\setbox0=\hbox{\the\SecNo.\the\PropNo}%
\bigbreak{\bf Proposition \the\SecNo.\the\PropNo.} \ignorespaces}

\def\Dem{\bigbreak{\it Démonstration. }}

\def\Eq{\global\advance\EqNo by 1
\global\setbox0=\hbox{\rm(\the\SecNo.\the\EqNo)}
\ifDisplaylines&\copy0\else\eqno{\copy0}\fi}

\def\FinDem\par{{\unskip\nobreak\hfil\penalty50%
\hskip1em\hbox{}\nobreak\Block
\parfillskip=0pt\finalhyphendemerits=0\par}}

%%%%%%%%%%%%%%%%%%%%%%%%%%%%%%%%%%%%%%%%
% List processing
%%%%%%%%%%%%%%%%%%%%%%%%%%%%%%%%%%%%%%%%

								% COMMENT: \rank(#1#2) writes the rank of the first
								% occurrence of the token #1 in the token list #2 (0 if #1 
								% does not occur.

\newif\iffound
\def\rank#1#2{{\foundfalse\count0=0 \getrank#1#2\end
\iffound\number\count0\else0\fi}} 
\def\getrank#1#2{\ifx#2\end\def\next#1{\relax}%
\else\iffound\else\advance\count0 by1\fi\let\next=\getrank
\ifx#2#1\foundtrue\fi\fi\next#1}

								% COMMENT: \List{ \a \b \c .... } causes \\\x to print the
								% rank of \x in the above list.  This strange macro is very
								% useful for references: one declares a list of references
								% at the beginning of the paper, and then typing a
								% reference with \\ in front gives the number of that
								% reference. In particular this automatically changes the
								% numbers when references are added or deleted.       					

\def\List#1{\def\\##1{\def##1{##1}\rank{##1}{#1}}}

\def\Ref#1; #2; #3; #4\par{\smallskip\noindent
\item{ [#1]}{\sc #2}, {\sl #3}, #4\par}
\def\Reff#1; #2; #3; #4; #5; #6; #7\par{\smallskip\noindent
\item{ [#1]}{\sc #2}, {\sl #3}, #4 {\bf #5} (#6) #7\par}

%%%%%%%%%%%%%%%%%%%%%%%%%%%%%%%%%%%%%%%%
% 9. Basic Postscript Commands
%%%%%%%%%%%%%%%%%%%%%%%%%%%%%%%%%%%%%%%%

									% COMMENT: These TeX commands code for postcript
         					% commands and are used to describe simple drawings.
                    % All measures are in millimeters. The recommended
                    % syntax is to use a graphic box for each picture. Note
                    % that there is no PS previewer on Macintosh, and therefore
                    % the result of PS commands will not appear on the screen.
                    % However \Write is NOT a PS command.
                    % A typical example, where an external picture called "figure1"
                    % is inserted (from the file ....pictures)
                    %         $$\OpenGraphicBox width 65mm height 30mm depth 0mm;
                    %         \leavevmode\special{picture figure1 scaled 500} 
                    %         \Translate(0, 0);
                    %         \Write(35, 31, $i$);
                    %         \Write(42, 31, $i+1$);
                    %         \Write(-7, 25, $\s_i$);
                    %         \Write(-8, 10, $\s_i^{-1}$);
                    %         \CloseGraphicBox$$

\newdimen\Dwidth
\newdimen\Dheight
\newdimen\Ddepth

\def\OpenGraphicBox width#1 height#2 depth#3;
{\dimen100=\parindent\parindent=0pt
\Dwidth=#1\Dheight=#2\Ddepth=#3
\dimen104=\Dheight\advance\dimen104 by \Ddepth
\catcode`\;=14\leavevmode
\raise-\Ddepth\vbox to\dimen104\bgroup\hsize=\Dwidth
\vglue\Dheight}

\def\Translate(#1, #2){\dimen105=-#2mm\vglue\dimen105\hskip#1mm}

\def\CloseGraphicBox{\vfill\egroup\parindent=\dimen100\catcode`\;=12}

\def\Write(#1, #2, #3){\rlap{\smash{\raise#2mm\hbox{\kern#1mm #3}}}}

%%%%%%%%%%%%%%%%%%%%%%%%%%%%%%%%%%%%%%%

\let\( = \langle 
\let\) = \rangle
\def\;{\, ; \,}

\def\card{{\rm card}}
\def\CCI{*}
\def\CCII{**}
\def\cl#1{\overline{\vrule height 5.5pt width 0pt #1}}
\def\cL{\mathord{/}}
\def\cLf{\cL_{\scriptscriptstyle \!\!\tilde f}\,}
\def\comp{\mathbin{\scriptscriptstyle\circ}}
\def\cR{\mathord{\backslash}}
\def\cRf{\cR_{\scriptscriptstyle f}}

\let\D=\Delta
\def\df{d_{\scriptscriptstyle f}\!}
\def\DivL{{\rm Div}}
\def\Dom{{\rm Dom}}
\def\dL{{}^*\hskip-0.4pt}
\def\DR{\Til D}

\let\e=\varepsilon
\def\eqf{\mathrel{\mathstrut\equiv}_{\!\scriptscriptstyle f}^{\!\scriptscriptstyle}}
\def\eqfp{\mathrel{\mathstrut\equiv}_{\!\scriptscriptstyle f}^{\!\scriptscriptstyle +}}
\def\eqfpp{\mathrel{\mathstrut\equiv}_{\!\scriptscriptstyle f}^{\!\scriptscriptstyle +\!+}}

\let\f = \phi

\def\Fr(#1){G_{\!_{#1}}}

\let\G=\Gamma
\def\GG(#1; #2){\langle #1; {#2}\rangle}
\def\gL{\mathbin{{\scriptstyle\wedge}}}
\def\gR{\mathbin{\Til{\scriptstyle\wedge}}}
\def\Gr#1{\(#1\)}

\def\id{{\rm id}}
\def\ii{^{-\!1}}
\let\ince=\subseteq

\def\jL{\mathbin{\Til{\scriptstyle\vee}}}
\def\jLf{\mathbin{\jL_{\!\scriptscriptstyle \tilde f}}}
\def\jR{\mathbin{\scriptstyle\vee}}
\def\JR{\mathbin{\bigvee}}
\def\jRf{\mathbin{\jR_{\!\scriptscriptstyle f}}}
\def\Mon#1{\(#1\)^{\scriptscriptstyle +}}

\def\Mo(#1){{\rm Mo}(#1)}

\def\norm(#1){\vert\!\vert#1\vert\!\vert}
\def\NR{\Til N}

\let\O = \Omega

\def\om{\mathord{\cdot}}
\def\op{}

\def\pp{\hbox{$\ldots$}}
\let\ppp=\cdots

\let\ra=\rightarrow
\def\Rf{R_{\! f}}

\let\s=\sigma
\def\S{\Sigma}
\def\spe(#1){\mathbin{\mathstrut\equiv}_{\!\scriptscriptstyle
f, Y}^{(#1)}}
\def\SPE{\mathbin{\mathstrut\equiv}_{\!\scriptscriptstyle f, Y}^{\scriptscriptstyle +++}}
\def\Sym{\hbox{\go S}}

\let\Til=\widetilde

\def\vLf{\mathrel{\Til{\vR}_{\scriptscriptstyle \!\tilde f}}}
\def\vR{\mathrel{\hbox{\sb\char'171}}}  % reversing relation
\def\vR{\rightarrow}  % reversing relation
\def\vRf{\vR_{\scriptscriptstyle \!f}}
\def\vv{\hskip-1pt\raise5pt\hbox{$\scriptscriptstyle\vee$}}

\def\WX{\Mo(\Sigma \cup \Sigma\ii)}
 
\hyphenation{par-ti-cu-lier se-con-de}

%%%%%%%%%%%%%%%%%%%%%%%%%%%%%%%%%%%%%%%%%

\List{\Adj \Adk \BaM \BDM \BrS \BMR \BuM \Cha \Chb \ClP
\Dfa \Dfb \Dff \Dfz \Dgc \Dgd \Dfx \Dlg \ElM \Eps \Gar \Mic
\Pic \Pid \Rmm \Tat \Thu}

%%%%%%%%%%%%%%%%%%%%%%%%%%%%%%%%%%%%%%%%%

\vskip 1 true cm 

\centerline{\bf GROUPES DE GARSIDE}

\bigskip\centerline{Patrick DEHORNOY}

\vskip 1 true cm 

{\narrower{\bf Abstract}. Define a Garside monoid
to be a cancellative monoid where right and left lcm's exist and
that satisfy additional finiteness assumptions, and a Garside
group to be the group of fractions of a Garside monoid.
The family of Garside groups contains the braid groups,
all spherical Artin groups, and various generalizations
previously considered\footnote{$^1$}{In particular, the
Garside groups considered in~[\\\Dfx] are special cases of
those considered here; the latter had been called "small
Gaussian" in earlier works; the name has been changed in
order to uniformize the terminology with works by Bessis,
Charney, Michel, and other authors}. Here we prove that
Garside groups are biautomatic, and that being a Garside
group is a recursively enumerable property, i.e., there exists an
algorithm constructing the (infinite) list of all small Gaussian
groups. The latter result relies on an effective, tractable
method for recognizing those presentations that define a 
Garside monoid.
\par}

{\narrower\bigskip{\bf Résumé}. Les monoïdes de Garside
sont introduits comme monoïdes simplifiables où existent
ppcm et pgcd et où sont satisfaites des conditions convenables
de finitude, et les groupes de Garside comme groupes de
fractions de monoïdes de Garside. La famille des groupes
de Garside  contient les groupes de tresses, les groupes d'Artin
sphériques, et diverses généralisations considérées
antérieurement\footnote{$^2$}{En particulier, les groupes de
Garside introduits dans~[\\\Dfx] sont un cas particulier de
ceux considérés ici, qui, eux, avaient été appelés petits
groupes gaussiens; ce changement de nom a été 
décidé afin d'uniformiser la terminologie avec d'autres
travaux récents de Bessis, Charney, Michel, entre autres.}.
Dans cet article, nous montrons que les groupes de Garside
sont bi-automatiques, et qu'être un groupe de Garside est une
propriété récursivement énumérable, c'est-à-dire qu'il existe
un algorithme énumérant la liste infinie de tous ces groupes.
Ce résultat repose sur une méthode effective pour reconnaître
les présentations des monoïdes de Garside.\par}

{\narrower\bigskip Key words: Presentations of monoids
and groups; word problem; rewriting systems; Cayley graph;
braid groups; Artin groups.

\par}

\vskip 1 true cm

 L'objet de cet article est l'étude par des méthodes algébriques et
combinatoires d'une classe de groupes, appelés groupes de Garside, qui contient en parti\-culier  les groupes de tresses
classiques, tous les groupes d'Artin sphériques, et
diverses extensions de ces groupes précédemment introduites.
Les groupes de Garside sont définis comme groupes de
fractions des monoïdes de Garside, ces derniers étant
définis par l'existence de notions convenables de plus petit
commun multiple (ppcm) et de plus grand commun diviseur
(pgcd), et la satisfaction de conditions de noethérianité et de
génération finie. On démontre ici:

\bigskip {\bf Théorème A.} {\sl Tout groupe de Garside est
bi-automatique.}

\bigskip {\bf Théorème B.} {\sl Etre un groupe de Garside et un 
monoïde de Garside sont des propriétés~$\S_1^0$,
c'est-à-dire récursivement énumérables.}

\bigskip Le théorème~A, qui résout positivement une
conjecture de~[\\\Dfx] (énoncée en termes de petit groupe
gaussien), implique en particulier que tout groupe de Garside
a un problème de mot et un problème de conjugaison
décidables, et qu'il satisfait une inégalité isopérimétrique
quadratique.

Le théorème~B signifie qu'il existe un algorithme (théorique) qui
énumère systématiquement toutes les présentations de groupes
 de Garside et de monoïdes de Garside. Il
repose sur l'existence d'un critère effectif caractérisant
certaines présentations des monoïdes de Garside.
Plus précisément, il existe des conditions explicites~$(\CCII)$
de complexité~$\S_1^0$ telle qu'on ait les résultats suivants,
dont les termes seront définis plus bas:

\bigskip {\bf Théorème B'.} {\sl (i) Si $M$ est un monoïde de Garside, et si $f$ est un sélecteur de ppcm
sur une partie génératrice~$\Sigma$ de~$M$, alors $M$ admet la
présentation complémentée~$\Mon{\Sigma; \Rf}$, et celle-ci
satisfait  aux conditions~$(\CCII)$.

(ii) Inversement, si $\Mon{\Sigma; \Rf}$ est une présentation
complémentée satisfaisant aux conditions~$(\CCII)$, alors le
monoïde qu'elle définit est un monoïde de Garside, et
$f$ est un sélecteur de ppcm sur~$\Sigma$.}

\bigskip Si la démonstration du théorème~B', qui occupe
la plus grande partie de l'article, est assez délicate, les
conditions~$(\CCII)$ qu'il met en jeu sont simples, et, en
particulier, leur implantation sur ordinateur est
aisée. La conjonction des théorèmes~A et~B' permet
donc  de construire de nombreux exemples de groupes
automatiques.

\bigskip Les propriétés algébriques des groupes de tresses et de leurs
extensions ont fait l'objet de nombreux travaux. Dans le cas des
groupes de tresses, la plupart des résultats classiques sur les
problèmes de mots et de conjugaison, les formes normales, et
l'automaticité sont présents ou implicites dans les travaux de Garside
[\\\Gar] et Thurston [\\\Thu] --- voir aussi [\\\Adk] [\\\ElM].
On sait que les méthodes et les résultats s'étendent à des
classes de groupes plus vastes: groupes d'Artin sphériques
[\\\BrS] [\\\Dlg] [\\\Cha] [\\\Chb], groupes des tresses des
groupes de réflexions complexes [\\\BMR], groupes de
Garside au sens de~[\\\Dfx] et~[\\\BDM] --- qu'on appellera
ici \<<groupes de Garside au sens restreint\>>. La classe des
groupes de Garside considérés ici, et introduits dans~[\\\Dfx]
sous l'appellation \<<petits groupes gaussiens\>>, inclut
strictement toutes les classes précédentes. Par exemple, les
groupes $\Gr{a  , b  \; a  b^pa   = b^q}$ sont pour $q \ge p
\ge 1$ des groupes de Garside n'appartenant à aucune des
familles précédentes.

Un groupe de Garside est le groupe de fractions d'un
monoïde dans lequel existe une bonne théorie de la
divisibilité et des ppcm (la définition précise sera
donnée au début de la section~1). Déjà relevée chez Garside et
Deligne, l'importance des ppcm dans les monoïdes de tresses a
été dégagée explicitement par Thurston dans la construction de
formes normales et, de là, de structures automatiques [\\\Thu]
[\\\Eps], et elle est au centre de l'étude développée
dans~[\\\Dfx] pour les groupes dits de Garside. En un sens, le
théorème~A n'est qu'une extension naturelle du résultat
analogue établi dans~[\\\Dfx] pour les groupes de Garside au
sens restreint, lequel ne fait que reprendre l'argument
de~[\\\Chb] pour les groupes d'Artin sphériques, qui, à son
tour, n'est qu'une adaptation de l'argument original
de~[\\\Thu]. Le point nouveau consiste ici à s'affranchir d'une
hypothèse technique superflue incluse dans la définition des
groupes de Garside au sens restreint, à savoir la propriété
(vérifiée dans le cas des tresses et des groupes d'Artin) que
l'élément fondamental~$\D$ est le ppcm des générateurs
minimaux. L'intérêt de la présente approche ne tient pas tant à
l'affaiblissement des hypothèses obtenu qu'à l'argument
utilisé pour établir le résultat, à savoir l'utilisation
systématique des compléments: si $M$ est un monoïde
admettant la simplification à gauche, et où deux éléments
quelconques admettent un unique ppcm à droite, nous
appellerons {\it complément} (à droite) de~$x$ dans~$y$, et
noterons~$x  \cR y$, l'unique élément~$z $ tel que
$x z $ est le ppcm à droite de~$x$ et~$y$. L'étude
pour elles-mêmes de l'opération~$\cR$ et de sa contrepartie à
gauche~$\cL$ semble nouvelle. Nous montrons ici comment
construire à l'aide de ces opérations des formes normales
explicites et des structures automatiques. De telles
constructions sont déjà présentes dans~[\\\Thu] et~[\\\Eps]
dans le cas des tresses, et le passage aux groupes de Garside serait une extension facile si  toutes les hypothèses
avaient été conservées. La tâche en fait n'est pas si aisée, car
les constructions initiales utilisaient au moins trois conditions
que nous éliminons ici: préservation de la longueur dans les
relations définissant le groupe, existence d'une structure de
groupe de Coxeter sous-jacente (le groupe symétrique et ses
opérations de treillis dans le cas des tresses), centralité de
l'élément~$\D^2$. S'affranchir de telles hypothèses exige de
reprendre l'étude du début, avec des arguments et un
enchaînement différents, tout en aboutissant à des
démonstrations plus simples et naturelles que celles
de~[\\\Cha] et~[\\\Dfx].  Il semble donc que les groupes de Garside constituent un cadre bien adapté, ce que confirment
d'autres travaux récents: [\\\Dfz], où il est montré que les
groupes de Garside, et, plus généralement, les groupes
gaussiens sont sans torsion, [\\\Pic], où est adaptée la solution
du problème de conjugaison décrite par Morton et ElRifai
[\\\ElM] dans le cas des tresses, [\\\Pid], où il est montré que
tout groupe de Garside se décompose en produit croisé de
groupes de Garside avec centre monogène.

Sauf dans des cas triviaux, il est très difficile de reconnaître les
propriétés combinatoires d'un groupe à partir d'une présentation, et de
nombreux résultats d'indécidabilité sont connus~[\\\BaM]. Il existe
en particulier très peu de critères permettant de reconnaître qu'un
monoïde se plonge dans le groupe de même présentation [\\\Adj]
[\\\Rmm], et on sait que la preuve d'un tel résultat par Garside dans
le  cas des groupes de tresses est à l'origine de la plupart des
développements algébriques ultérieurs sur ces groupes. Le
théorème~B repose sur une nouvelle méthode dans cette
direction. L'outil essentiel utilisé ici est une opération
combinatoire sur les mots appelée {\it redressement}.
Introduite dans~[\\\Dfa], et considérée sous une forme
similaire dans~[\\\Tat], puis étudiée sytématiquement
dans~[\\\Dff], cette opération est une contrepartie syntaxique
à l'opération de complément: à partir d'une certaine
fonction~$f$ déterminée par la présentation lorsque celle-ci
est d'un certain type dit complémenté, le redressement permet
de définir une opération binaire~$\cRf$ sur les mots de sorte
que, si $M$ est un monoïde engendré par un
ensemble~$\Sigma$, et si $u, v$ sont des mots sur~$\Sigma$
représentant respectivement les éléments~$x$ et~$y$ de~$M$,
alors le mot~$u \cRf v$ représente l'élément~$x  \cR
y$, sous la condition que la fonction~$f$, c'est-à-dire la
présentation considérée, vérifie une certaine condition dite
propriété du cube, déjà implicite dans le théorème~H
de~[\\\Gar]. Il a été montré dans~[\\\Dfx] que tous les
groupes de Garside admettent une présentation
complémentée, et, réciproquement, on y a donné des
conditions suffisantes pour qu'une présentation
complémentée définisse un groupe de Garside. Ces conditions
au demeurant sont d'intérêt pratique limité, car elles
nécessitent d'une part d'établir la noethérianité d'une relation
et, d'autre part, de vérifier la propriété du cube pour une
infinité de mots: dans les deux cas, il s'agit de conditions de
type infinitaire (conditions respectivement~$\Pi_1^1$ et
$\Pi_1^0$ dans le langage de la théorie de la récursivité, 
c'est-à-dire mettant en jeu une quantification universelle
portant respectivement sur les suites de mots et les mots).
Dans certains cas particuliers, comme celui des tresses et, plus
généralement, lorsque les relations de la présentation
considérée préservent la longueur, la condition de
noethérianité est triviale, et la propriété du cube se réduit à
une vérification finie, qui constitue le travail de Garside
dans~[\\\Gar]. Notre travail ici consiste à montrer que, dans le
cas général, la condition de noethérianité peut être déduite de
conditions effectives plus faibles, et qu'il est suffisant de
vérifier la propriété du cube pour des triplets de mots pris
dans un certain ensemble fini. De la sorte, et comme énoncé
dans le théorème~B', nous obtenons un critère algorithmique
--- et simple: si, par exemple, nous entrons la
présentation standard d'un monoïde de tresses dans
l'algorithme, celui-ci répondra en un temps fini qu'il
s'agit d'un monoïde de Garside, et il fournira une
description explicite d'une structure automatique. Le
prix à payer pour ces résultats est une analyse fine du
processus de redressement des mots --- lequel jouait
déjà un rôle technique essentiel dans la construction
d'un ordre total sur les tresses~[\\\Dgd].

\bigskip Le plan de l'article est le suivant. Dans la section~1, nous
définissons les notions de monoïde de Garside et
de groupe de Garside; nous établissons des propriétés de base
des opérations de complément, et nous décrivons  un treillis fini
attaché à chaque monoïde de Garside et qui le caractérise.
Dans la section~2,  on établit des formules de dualité reliant
les opérations
 de ppcm, pgcd, et complément à gauche et à droite dans tout monoïde de Garside. En fait, on travaille ici dans un cadre plus
général: ainsi, on obtient à la fois des formules de dualité utiles pour la
preuve du théorème~A, et une caractérisation des monoïdes de Garside par des hypothèses plus faibles (mais moins naturelles) que
celles de la définition initiale, lesquelles seront prises comme
point de départ pour la preuve du théorème~B'. Dans la
section~3, on démontre le théorème~A. Grâce au
pgcd, on définit une forme normale unique pour les éléments
de tout groupe de Garside: cette forme, la \<<mixed
canonical form\>> de~[\\\Eps], est bien
connue dans le cas des tresses, le point spécifique ici étant la
simplicité de la preuve d'automaticité, ainsi que la possibilité
de construire effectivement des automates (ou des
transducteurs) calculant les formes normales. En particulier,
on montre comment calculer le pgcd à gauche {\it et} le pgcd à
droite avec l'élément~$\D$ au moyen d'un automate dont les
états sont des fonctions sur la clôture des atomes par
complément (et non des diviseurs de~$\D$ ainsi qu'il est
classique pour le pgcd à gauche). Dans la section~4, on montre
que tout groupe de Garside admet une présentation de
forme syntactique parti\-culière, dite complémentée, et que
toutes les opérations du monoïde, ppcm, pgcd, compléments,
peuvent se calculer à partir de l'opération de redressement
des mots déduite de la présentation. Ceci donne en particulier
une solution effective au problème de mot (du monoïde et du
groupe). Dans la section~5, et grâce à la caractérisation de la
section~2, on montre que les présentations complémentées
de groupes de Garside satisfont à certaines conditions
nécessaires et suffisantes, dont la propriété du cube.
A ce stade, le théorème~B' n'est pas établi, car le critère
obtenu requiert de vérifier la propriété du cube
pour tous les triplets de mots: l'objet de la
section~6 est alors de montrer qu'il suffit d'effectuer cette
vérification pour un ensemble fini de mots, à savoir la
clôture des générateurs par l'opération~$\cRf$. Ceci
établit le résultat cherché, et fournit un critère
pratique, ainsi que l'illustre un exemple.

\bigskip Si $\Sigma$ est un alphabet, et $R$ une famille de relations
sur~$\Sigma$, c'est-à-dire une famille de paires de mots sur~$\Sigma$, nous
noterons $\Mon{\Sigma \; R}$ le monoïde engendré par~$\Sigma$ et
présenté par les relations~$R$, et $\Gr{\Sigma \; R}$ le groupe
correspondant.

\bigskip L'auteur remercie Matthieu Picantin
pour ses remarques judicieuses.

\Sec Calcul des compléments

\bigskip Si $M$ est un monoïde, et $x$, $y$ des éléments
de~$M$, on dira que $x$ est diviseur à gauche de~$y$, ou, de
façon équivalente, que $y$ est multiple à droite de~$x$, s'il
existe~$z $ vérifiant~$x  z  = y$; on parlera de diviseur ou de
multiple propre si, en outre, on a~$z  \not= 1$. De là une
notion naturelle de plus petit multiple commun (ppcm) à
droite: $z $ est un ppcm à droite de~$x$ et~$y$ si
$z $ est un multiple à droite de~$x$ et de~$y$, et un diviseur
à gauche de tout multiple à droite commun à~$x$ et~$y$. De
la même façon, $z $ est un plus grand commun diviseur à
gauche, ou pgcd à gauche, de~$x$ et~$y$ si $z $ est un
diviseur à gauche de~$x$ et~$y$, et un multiple à droite de
tout diviseur à gauche commun à~$x$ et~$y$.

Il n'y a en général aucune raison pour qu'un ppcm~$z $
de~$x$ et~$y$, s'il existe, soit unique, non plus
que l'élément~$y '$ vérifiant $z  = x  y '$. Cependant,
des conditions assez faibles garantissent cette unicité.

\Lem\«\Uniqueness»
{\sl Soit $M$ un monoïde simplifiable à
gauche où $1$ est le seul élément inversible. Alors la
relation \<<être un diviseur à gauche propre\>> est un
ordre partiel sur~$M$, et les ppcm à droite et pgcd à 
gauche de deux éléments sont uniques quand ils existent.}

\Dem Par hypothèse, la conjonction de $x  \not= 1$ et $y 
\not= 1$ implique $x y  \not= 1$ dans~$M$, donc la
relation \<<être un diviseur à gauche propre\>> est
transitive; que $M$ soit simplifiable à gauche entraîne que
cette relation est irréflexive, et c'est donc un ordre (strict).
Quand ils existent, le ppcm à droite et le pgcd à gauche de
deux éléments~$x , y$ sont la borne supérieure et la
borne inférieure de~$x$ et~$y$ dans l'ordre
ci-dessus.\FinDem

\Def Sous les hypothèses du lemme précédent, nous
noterons $x  \jR y$ le ppcm à droite de~$x$ et~$y$, quand il
existe; dans ce cas, l'unique élément~$z $ vérifiant $x  \jR
y  = x  z $ sera noté~$x  \cR y$, lu \<<$x$
sous~$y$\>>, et appelé {\it complément à droite} de~$y$
sur~$x$. 

\bigskip On a donc
$$x  \jR y  = x  \op (x  \cR y ) = y  \op (y  \cR x ) \Eq\«100»$$
dès que $x  \jR y$ est défini. De façon symétrique (donc sous
l'hypothèse que le monoïde est simplifiable à droite), nous
notons $x  \jL y$ le ppcm à gauche de~$x$ et~$y$ quand il
existe, et $x  \cL y$ l'unique élément~$z $ vérifiant
$x  \jL y  = z  x$. La contrepartie de~«100» est alors
$$x  \jL y  = (x  \cL y ) \op y  = (y  \cL x ) \op x 
\Eq\«200»$$ 
(Figure~1.1). Noter que $x$ est un diviseur à gauche
de~$y$ si et seulement si
$y  \cR x$ est défini et égal à~$1$: dans ce cas, $x  \cR
y$ est le quotient~\<<$x \ii y$\>> et, de même, $x 
\cL y$ est le quotient~\<<$x  y \ii$\>> dans le cas où
$y$ est un diviseur à droite de~$x$ --- ce qui explique et
justifie nos notations.

\midinsert
$$\OpenGraphicBox width 80mm height 30mm
depth0mm;
\leavevmode\special{illustration prod4.eps}
\CloseGraphicBox$$ 
\centerline{{\bf Figure 1.1.} Complément à droite et à gauche}
\endinsert

L'objet de notre étude est la famille des monoïdes appelés  gaussiens
dans~[\\\Dfx] et~[\\\Dfz]. Ceux-ci sont définis en termes de ppcm et
de pgcd: essentiellement, un monoïde gaussien est un monoïde avec
uniques ppcm et pgcd à droite et à gauche, et un monoïde de Garside est un monoïde gaussien satisfaisant une condition forte de
génération finie. Pour une définition précise, nous aurons besoin
d'une notion supplémentaire:

\Def
Soit $M$ un monoïde. On dit qu'un élément~$x$ de~$M$
est un {\it atome} si  $x$ est  distinct de~$1$ et que $x  =
y z $ entraîne $y  = 1$ ou $z  = 1$. On dit que $M$ est
{\it atomique} si $M$ est engendré par ses atomes et si, de
plus, pour tout~$x$ dans~$M$, la {\it
norme}~$\norm(x )$ de~$x$, définie comme la borne
supérieure des longueurs des décompositions de~$x$ en
produit d'atomes, est finie.

\bigskip On remarquera que $1$ est nécessairement le seul
élément inversible dans un monoïde atomique, puisque,
par définition, on a $\norm(x  y ) \ge \norm(x ) +
\norm(y )$ pour tous~$x , y$, et $\norm(x ) \ge 1$ pour $x 
\not= 1$. 

La définition suivante d'un monoïde gaussien apparaît
dans~[\\\Dfx]:

\Def
(i) On dit qu'un monoïde~$M$ est {\it gaussien} si $M$ est
atomique, simplifiable, et si deux éléments quelconques
de~$M$ admettent un ppcm à droite et à gauche, et un pgcd à
droite et à gauche. 

(ii) On dit qu'un monoïde gaussien~$M$ est un {\it monoïde
de Garside} s'il contient un {\it élément de Garside}, ce
dernier étant défini comme un élément~$\D$ tel que les
diviseurs à gauche de~$\D$ coïncident avec les diviseurs à
droite de~$\D$, ils soient en nombre fini, et ils
engendrent~$M$.

(iii) Tout monoïde gaussien satisfait aux conditions de Ore,
et admet donc un groupe de fractions. On dit qu'un
groupe~$G$ est un {\it groupe gaussien} (\resp un {\it
groupe de Garside}) s'il existe un monoïde gaussien~$M$
(\resp un monoïde de Garside~$M$) tel que $G$ soit le groupe
de fractions de~$M$.

\Ex 
Les groupes de tresses [\\\Gar] et, plus généralement,
tous les groupes d'Artin sphériques, c'est-à-dire associés
à un groupe de Coxeter fini, sont des groupes de Garside [\\\BrS] 
[\\\Dlg], l'élément dit fondamental étant un élément de
Garside. Plus généralement, tous les groupes de Garside
au sens de~[\\\Dfx] sont des groupes de Garside au
sens ci-dessus, mais la réciproque est fausse: un groupe de
Garside au sens restreint de~[\\\Dfx] est un groupe de
Garside dans lequel, de surcroît, le ppcm à droite des atomes
est un élément de Garside. Un contre-exemple typique est
$\Mon{a  , b  \; a   b  a   = b^2 }$: il y a ici deux atomes, à
savoir~$a  $ et~$b$, dont le ppcm est~$b^2$, qui n'est pas un
élément de Garside, puisque $a  b $ est un diviseur à gauche
mais pas à droite de~$b^2$; par contre, l'élément $\D = b^3$
est un élément de Garside.

\bigskip Nous allons dans la suite établir divers résultats sur les
monoïdes de Garside, en particulier l'équivalence de
la définition originale avec plusieurs définitions
alternatives. Tous ces résultats reposent sur les propriétés
algébriques des opérations de ppcm et de complément, que
nous allons établir maintenant. Il sera utile de disposer de
telles propriétés dans un cadre plus général que celui des monoïdes de Garside. En particulier, nous considérerons le cas où
l'opération de ppcm (à droite) n'est pas nécessairement partout
définie. Pour ce faire,  nous utiliserons un symbole~$\bot$
signifiant \<<non défini\>>: de la sorte,  $x 
\jR y  = \bot$ signifie que $x$ et $y$ n'ont pas de ppcm,
et une égalité telle que $x  \jR y  = x ' \jR y '$ signifie
que soit $x  \jR y$ et $x ' \jR y '$ sont définis et ils
sont égaux, soit ni $x  \jR y$, ni $x ' \jR y '$ ne sont
définis. On convient qu'on a toujours $x  \bot = \bot x  =
\bot$. Ainsi les égalités~«100» et~«200» sont toujours
satisfaites, que les ppcm et compléments mentionnés soient
définis ou non.

Pour formuler de façon compacte nos hypothèses,
nous utiliserons les abréviations suivantes:

\Def On dira qu'un monoïde~$M$ satisfait la
condition

-  $(C_0)$ si l'élément~$1$ est le seul
inversible de~$M$;

-  $(C_0^+)$ si $M$ est atomique;

-  $(C_1)$ si $M$ admet la simplification à gauche;

-  $(\Til C_1)$ si $M$ admet la simplification à droite;

-  $(C_2)$ si deux éléments de~$M$ admettant un multiple
à droite commun admettent un ppcm à droite;

-  $(C_2^+)$ si deux éléments quelconques de~$M$ admettent
un ppcm à droite;

-  $(C_3)$ si $M$ possède une partie génératrice finie~$P$
close par~$\cR$, c'est-à-dire telle que $x , y  \in P$
entraîne~$x  \cR y  \in P$; 

-  $(C_3^+)$ si $M$ possède une partie génératrice finie~$S$
close par~$\cR$ et~$\jR$, c'est-à-dire telle que $x , y  \in
S$ entraîne~$x  \cR y  \in S$ et $x  \jR y  \in S$. 

\bigskip Noter que $(C_0^+)$ entraîne~$(C_0)$, et, de même,
$(C_2^+)$ entraîne~$(C_2)$, et $(C_3^+)$ entraîne~$(C_3)$.
La conjonction de~$(C_0)$, $(C_1)$ et $(C_2)$ est le cadre
naturel pour que les opérations de ppcm et de complément à
droite soient bien définies --- donc, en particulier, pour
que $(C_3)$ et $(C_3^+)$ aient un sens --- et $(C_2^+)$ est
la condition additionnelle garantissant que ces opérations
soient partout définies. Sous l'hypothèse que $(C_0)$,
$(C_1)$ et $(C_2)$ sont satisfaites, l'égalité~«100» est
toujours valable, de même que $x  \jR y  = y  \jR x$, et
$$(x  y ) \jR (x  z ) = x  \op (y  \jR z ).\Eq\«101»$$ 

\Lem\«\Rules»
{\sl Soit $M$ un monoïde satisfaisant~$(C_0)$, $(C_1)$
et~$(C_2)$. Pour tous $x ,
y , z $ dans~$M$, on a:
$$\Displaylines{
(x  y ) \cR z  = y  \cR (x  \cR z ), \qquad 
z  \cR (x  y ) = (z  \cR x ) \cdot ((x  \cR z ) \cR y ) 
\Eq\«103»\cr 
(x  \jR y ) \cR z  = (x  \cR y ) \cR (x 
\cR z ) = (y  \cR x ) \cR (y  \cR z ), \qquad 
z  \cR (x  \jR y ) = (z  \cR x ) \jR (z  \cR y )
\qquad \Eq\«105»\cr }$$ 

}\Dem Pour~«103», on a, en utilisant~«100»
et~«101»,
$$\displaylines{
x  y  ((x  y ) \cR z ) 
= (x  y ) \jR z  
= (x  y ) \jR (x  \jR z ) 
= (x  y ) \jR (x  (x  \cR z ))
= x (y  \jR (x  \cR z ))
= x  y (y  \cR (x  \cR z )), \cr
z  ((x  y ) \cR z ) 
= (x  y ) \jR z  
= x  y ((z  \cR x ) \cR y ) 
= x (z  \cR x )(y  \cR (z  \cR x )) 
= z (x \cR z )(y  \cR (z  \cR x )), \cr}$$
et on simplifie à gauche par~$x  y$ dans le premier
cas, et par~$z $ dans le second. Pour~«105»,
appliquer~«103»\ à $x  \jR y  = x (x  \cR y )$ donne
directement
$(x  \jR y ) \cR z  = (x  \cR y ) \cR (x  \cR z )$.
Comme $\jR$ est symétrique, cette expression est
aussi~$(y  \cR x )
\cR (y  \cR z )$. Enfin, on trouve, en appliquant~«103»,
l'égalité précédente, et~«100»,
$$z  \cR (x  \jR y )
= z  \cR (x (x  \cR y ))
= (z  \cR x ) ((x  \cR z ) \cR (x  \cR y )) 
= (z  \cR x ) ((z  \cR x ) \cR (z  \cR y )) 
= (z  \cR x ) \jR (z  \cR y ).$$
(Figure~1.2). Remarquer que les formules sont valides
lorsqu'un des ppcm n'est pas défini, chacune des expressions
ayant alors la valeur~$\bot$, c'est-à-dire étant non
définie.\FinDem

\midinsert
$$\OpenGraphicBox width 65mm height 52mm
depth0mm;
\leavevmode\special{illustration comp4.eps}
\CloseGraphicBox$$ 
\centerline{{\bf Figure 1.2.} Complément itéré}
\endinsert

\bigskip Dans le contexte précédent, pour $X, Y \ince M$ et
$k \ge 1$, on pose $X \cR Y = \{ x  \cR y  ~; x  \in X,
y  \in Y\}$, $X^k =
\{x_1 \ppp x_k \,;\, x_1, \pp, x_k \in X\}$,
et $X^0 = \{1\}$.

\Lem\«\Compl» 
{\sl Soit $M$ un monoïde satisfaisant aux conditions~$(C_0)$,
$(C_1)$ et $(C_2)$, et $P$ une partie génératrice de~$M$ 
close par~$\cR$. Alors, pour tout~$k$, on a $M \cR P^k
\ince P^k$, $P^k$ est close par~$\cR$ et par diviseur à droite,
et $M$ satisfait à la condition~$(C_2^+)$.}

\Dem On montre par induction sur~$i + k$ que la conjonction de $x 
\in P^i$ et
$y  \in P^k$ entraîne $x  \cR y  \in P^k$. Pour $k = 0$, c'est-à-dire pour $y  = 1$, on a $x  \cR y  =
1$, et le résultat est vrai. Supposons $k = 1$. Pour $i =
0$, soit $x  = 1$, le résultat est vrai. Pour $i = 1$, le résultat
est l'hypothèse que $P$ est clos par~$\cR$. Pour $i \ge 2$,
écrivons $x  = x_1x '$ avec $x_1 \in P$ et $x ' \in P^{i-1}$.
Par~«103», on a $x  \cR y  = x ' \cR (x_1 \cR y )$. L'hypothèse
d'induction donne $x_1 \cR y  \in P$, d'où $x ' \cR
(x_1 \cR y ) \in P$. Supposons enfin  $k \ge 2$. On écrit $y  =
y_1 y '$ avec $y ' \in P$ et $y ' \in P^{k-1}$. Par~«103», on a
$x  \cR y  = (x  \cR y_1) \op ((y_1 \cR x ) \cR
y ')$. L'hypothèse d'induction donne $x  \cR y_1 \in
P$ et $y_1 \cR x  \in P^i$. On déduit $(y_1 \cR x ) \cR
y ' \in P^{k-1}$, d'où $x  \cR y  \in P^k$. 

On a ainsi montré $P^i \cR P^k \ince P^k$ pour
tout~$i$, d'où $M \cR P^k \ince P^k$ si $P$
engendre~$M$, c'est-à-dire que $M$ est la réunion des
ensembles~$P^i$. On a aussi $P^k \cR P^k \ince P^k$,
donc $P^k$ est clos par l'opération~$\cR$.
De plus, supposons $y  \in P^k$ et $y  =z x$. On a alors
$x  = z  \cR y$, donc $x  \in P^k$.

Les calculs précédents montrent que $x  \cR
y$, donc aussi $x  \jR y$, existent pour tout~$x$
dans~$M$ et tout~$y$ dans~$P^k$, pour tout~$k$.
Comme $P$ engendre~$M$, on déduit que $x  \cR y$
et $x  \jR y$ existent pour tous~$x , y$ dans~$M$,
et $M$ satisfait donc à la condition~$(C_2^+)$. \FinDem

\Lem\«\Closure» 
{\sl Soit $M$ un monoïde satisfaisant~$(C_0)$, $(C_1)$
et~$(C_2)$, et  $P$ un sous-ensemble de~$M$ clos
par~$\cR$. Soit $S$ la clôture de~$P$ par~$\jR$.

(i) L'ensemble~$S$ est clos par~$\cR$ et~$\jR$. Si $P$ est fini,
alors $S$ est fini, on a $\card(S) \le 2^{\card(P)}$; par
conséquent, la condition~$(C_3)$ entraîne la
condition~$(C_3^+)$.

(ii) Pour tous~$i, k$, on a $S^i \cR S^k \ince S^k$ et
$S^i \jR S^k \ince S^{\sup(i, k)}$; en
particulier, $S^k$ est clos par~$\cR$ et~$\jR$ pour tout~$k$.}

\Dem (i) Par construction, tout élément de~$S$ peut être
exprimé comme ppcm à droite d'un sous-ensemble fini
de~$P$, donc, si $P$ a $n$~éléments, $S$ a au plus
$2^n$~éléments. Par construction, $S$ est clos par~$\jR$, et
il s'agit de montrer qu'il est également clos par~$\cR$.
Supposons $x  = x_1 \jR \ppp \jR x_p$, $y  = y_1
\jR \ppp \jR y_q$ avec $x_1$, \pp, $y_q \in P$. Nous
voulons établir $x  \cR y  \in S$. D'abord, «105»
entraîne $x  \cR y  = (x  \cR y_1) \jR \ppp \jR (x 
\cR y_q)$, donc il suffit de montrer  $x  \cR y_j
\in P$ pour tout~$j$. On utilise une induction sur~$p \ge 1$.
Pour $p = 1$, le résultat est l'hypothèse que $P$ est clos
par ppcm à droite. Sinon, posons $x ' = x_1 \jR \ppp \jR
x_{p-1}$. Par~«105», nous obtenons
$$x  \cR y_j 
= (x ' \jR x_p) \cR y_j 
= (x ' \cR x_p) \cR (x ' \cR y_j);$$  
L'hypothèse d'induction entraîne $x ' \cR x_p
\in P$ et $x ' \cR y_j \in P$, d'où $(x ' \cR x_p) \cR
(x ' \cR y_j) \in P$ puisque $P$ est clos par~$\cR$. 

(ii) L'ensemble~$S$ étant clos par~$\cR$, le lemme~«\Compl»
donne $S^i \cR S^k \ince S^k$ directement. On établit la
relation $S^i \jR S^k \ince S^{\sup(i, k)}$ par récurrence sur
$\inf(i, k)$. Le résultat est trivial pour
$\inf(i, k) = 0$. Pour $\inf(i, k) = 1$, soit par exemple $k
= 1$, on utilise une récurrence sur~$i$: supposant $x \in S^i$
et $y \in S$, on pose $x = x_1 x'$ avec $x_1 \in S$ et $x' \in
S^{i-1}$. Par les formules du lemme~«\Rules», on trouve
$$x \jR y = (x_1 \jR y) ((x_1 \cR y) \cR x') \in S \cdot S^{i-1}
= S^i.$$
Supposons enfin $\inf(i, k) \ge 2$. On
écrit $x = x_1 x'$, $y = y_1 y'$ avec $x_1, y_1 \in S$, $x_1
\in S^{u-1}$, et $y' \in S^{k-1}$. Les règles du 
lemme~«\Rules» donnent
$$x \jR y = (x_1 \jR y_1) \, (((x_1 \cR y_1) \cR x') \jR
((y_1 \cR x_1) \cR y')).$$
On a alors $x_1 \jR y_1 \in S$, puis, d'après ce qui
précède, $(x_1 \cR y_1) \cR x' \in S^{i-1}$, et $(y_1 \cR x_1)
\cR y' \in S^{k-1}$, d'où $((x_1
\cR y_1) \cR x') \jR ((y_1 \cR x_1) \cR y') \in S^{\inf(i-1,
k-1)}$ par hypothèse de récurrence, et, finalement,
$x \jR y \in S^{\sup(i, ,)}$.\FinDem

\bigskip Sous les hypothèses précédentes, tout élément
de~$S$ qui s'exprime comme le ppcm à droite de
$k$~éléments de~$P$ peut aussi s'exprimer comme produit
de $k$~éléments de~$P$. En effet, on peut établir par
récurrence une formule générale du type
$$x  \jR y  \jR z  \jR \ppp = x  \cdot (x  \cR y ) \cdot
((x  \cR y ) \cR (x  \cR z )) \cdot \ppp.$$
Nous allons dans la suite donner plusieurs
définitions équivalentes des monoïdes de Garside. La
première provient de la remarque que, si $M$ est un
monoïde de Garside et que $\D$ est un élément de
Garside dans~$M$, alors l'ensemble~$S$ des diviseurs
de~$\D$ est clos par les opérations de ppcm, pgcd et
complément à droite et à gauche. L'existence d'un tel ensemble peut en
fait être prise comme définition:

\Prop\«\Alternative»
{\sl Soit $M$ un monoïde gaussien, $\D$ est
élément de~$M$, et $S$ une partie de~$M$. Alors il y a
équivalence entre

(i) $\D$ est élément de Garside, et $S$ est l'ensemble des
diviseurs (à gauche ou à droite) de~$\D$;

(ii) $S$ est une partie génératrice finie de~$M$ qui est close
par diviseur, ppcm, complément et pgcd à droite et à
gauche, et $\D$ est le ppcm à droite de~$S$.}

\Dem Supposons~(i). Supposons $x  \in S$ et $x  =
y z $: alors il existe $x '$ vérifiant $x x ' = \D$ , d'où
$y  (z x ') = \D$, et $y$ est dans~$S$. Ainsi, $S$ est clos
par diviseur à gauche, et, par un argument symétrique,
par diviseur à droite.

Soient maintenant $x$, $y$ dans~$S$. Par hypothèse, $\D$
est multiple à droite de~$x$ et~$y$, donc de~$x  \jR y$. Par
conséquent, $x  \jR y$ est dans~$S$, et il en est de même
de~$x  \cR y$, qui en est un diviseur à droite.
Le cas du ppcm et du complément à gauche est
symétrique; le cas des pgcd est trivial. Enfin, $\D$, étant
élément de~$S$, en est nécessairement ppcm à droite et à
gauche puisque $S$ est clos par ces opérations.

Réciproquement, supposons~(ii). Le premier problème est
de montrer que l'ensemble~$\Til S$ des diviseurs à gauche
de~$\D$ coïncide avec~$S$. D'abord, par construction, tout
élément de~$S$ est diviseur à gauche de~$\D$, et $S$ est
donc inclus dans~$\Til S$. Soit $x$ un élément quelconque
de~$\Til S$. Il existe donc~$x '$ vérifiant $x x ' = \D$, d'où
$x  = \D \cL x '$. Par le lemme~«\Closure», on a $M \cR S
\ince S$, et donc, symétriquement, $S \cL M \ince S$, d'où en
particulier $x  = \D \cL x '  \in S$. Donc $S$ est
l'ensemble des diviseurs à gauche de~$\D$. 

Le second problème est de montrer de même que $S$
coïncide avec l'ensemble des diviseurs à droite de~$\D$. Ce
point résultera du calcul précédent (renversé par
symétrie) pourvu qu'on montre que $\D$ est le ppcm à
gauche de~$S$. Comme $\D$ est dans~$S$, il s'agit de
montrer que tout élément~$x$ de~$S$ est un diviseur à
droite de~$\D$. Soit $x$ un tel élément. Par hypothèse, 
$x  \jL \D$ existe et est dans~$S$: il existe $y$ vérifiant $x 
\jL \D = y  \D$, et, puisque $y  \D$ est dans~$S$, il
existe~$z $ vérifiant $y  \D z  = \D$. Comme $\D$ est dans~$S$,
et que $S$ est clos par diviseur à gauche et à droite, $y$ et~$z $
sont dans~$S$. Il existe donc $y '$ dans~$S$, à savoir $y ' = y 
\cR \D$, vérifiant
$y  y ' = \D$, et, puisque la simplification à gauche est
possible, on déduit $y ' = \D z $. Puisque $y '$ est
dans~$S$, il existe ensuite $y ''$ dans~$S$, à savoir $y '' =
y ' \cR \D$, vérifiant $y ' y '' = \D$, d'où $\D = \D z 
y ''$. Ceci entraîne $z  y '' = 1$, d'où $z  = y '' = 1$,
donc $y  \D = \D$, et, finalement, $y  = 1$, ce qui montre
que $x$ est diviseur à droite de~$\D$. \FinDem

\bigskip On notera que l'argument précédent utilise
seulement~$(C_0)$, et non~$(C_0^+)$.

\Def Si $M$ est un monoïde de Garside, on appelle
{\it primitifs à droite} les éléments de la clôture des atomes
de~$M$ par l'opération~$\cR$, et {\it simples} les éléments
de la clôture des éléments primitifs par l'opération~$\jR$ ---
qui est aussi la clôture des atomes par~$\cR$ et~$\jR$.

\bigskip Si $M$ est un monoïde de Garside, alors l'ensemble
des atomes de~$M$ est la partie génératrice minimale de~$M$,
l'ensemble des éléments primitifs à droite est la partie
génératrice close par~$\cR$ minimale, l'ensemble~$S$ des
éléments simples est la partie génératrice close par~$\cR$
et~$\jR$ minimale, et son ppcm est l'élément de Garside de
norme minimale\footnote{$^1$}{dans le cas des groupes
d'Artin, les éléments simples sont aussi appelés {\it réduits}
dans la littérature}. Avec les notations précédentes, la
structure~$(S, \jR, \gL)$ est un treillis (fini), et il en est de
même de la structure symétrique~$(S, \jL, \gR)$.

\Prop
{\sl Un monoïde de Garside est déterminé
par l'ensemble de ses éléments primitifs à droite et
la restriction de l'opération~$\cR$ à ces éléments.}

\Dem Soit $M$ un monoïde de Garside, et $P$ l'ensemble
de ses éléments primitifs à droite. Par hypothèse,  $P$
engendre~$M$, donc il suffit de montrer que, pour tous
$x_1, \pp, x_p$, $y_1, \pp, y_q$ dans~$P$, la
condition $x_1 \ppp x_p = y_1 \ppp y_q$ peut
s'exprimer en terme de la restriction de~$\cR$ à~$P$. Or
$x_1 \ppp x_p = y_1 \ppp y_q$ équivaut à $(x_1
\ppp x_p) \cR (y_1 \ppp y_q) = (y_1
\ppp y_q) \cR (x_1 \ppp x_p) = 1$, et les règles de calcul du
lemme~«\Rules» montrent que $(x_1 \ppp x_p) \cR (y_1
\ppp y_q)$ s'exprime comme produit de $q$~éléments
de~$P$ déterminés de proche en proche à partir de~$x_1,
\pp, y_q$ à l'aide de l'opération~$\cR$. Par exemple, pour
$p = q = 2$, on trouve
$$(x_1 x_2) \cR (y_1 y_2) = (x_2 \cR (x_1 \cR y_1)) \om
(((x_1 \cR y_1) \cR x_2) \cR ((y_1 \cR x_1) \cR y_2)).$$
En vertu de la condition~$(C_0)$, le produit de ces
$q$~éléments vaut~$1$ si et seulement chacun d'eux
vaut~$1$, ce qui donne une condition du type requis. Par
exemple, $x_1x_2 = y_1 y_2$ est équivalent
à la conjonction des quatre égalités
$$\Displaylines{
x_2 \cR (x_1 \cR y_1) = 1, \quad
((x_1 \cR y_1) \cR x_2) \cR ((y_1 \cR x_1) \cR y_2) = 1, \cr
y_2 \cR (y_1 \cR x_1) = 1, \quad
((y_1 \cR x_1) \cR y_2) \cR ((x_1 \cR y_1) \cR x_2) = 1. 
&\Block\cr
}$$

\Prop
{\sl Un monoïde de Garside	$M$ est déterminé par son graphe
caractéristique, défini comme la restriction aux éléments simples
du graphe de Cayley atomique, c'est-à-dire la famille des triplets
$(x , a  , y )$ où $x$ et $y$ sont simples et $a  $ est un
atome vérifiant $x  a   = y$.}

\Dem Notons~$S$ l'ensemble des éléments simples de~$M$, et~$\G$
son graphe caractéristique. D'après la proposition précédente, il suffit
de montrer que, pour $x$, $y$ dans~$S$, la valeur de~$x  \cR
y$ est déterminée par~$\G$. Soient $x$, $y$ des éléments
de~$S$, donc des sommets de~$\G$. Le ppcm à droite de~$x$
et~$y$ est la borne supérieure de~$x$ et~$y$ dans~$\G$,
c'est-à-dire l'unique sommet~$z $ accessible depuis~$x$ et~$y$
dans~$\G$ tel que, pour tout sommet~$z '$ accessible
depuis~$x$ et~$y$, $z '$ soit accessible
depuis~$z $. Alors $x  \cR y$ est l'unique
sommet~$y '$ de~$\G$ tel qu'il existe un chemin de~$1$
à~$y '$ qui porte les mêmes étiquettes que le chemin
de~$x$ à~$z $. En effet, si $a  _1, \pp, a  _k$ est la
suite des étiquettes du chemin de~$x$ à~$z $
dans~$\G$, on a
$x \cR y  = a  _1 \ppp a  _k$ dans~$M$. Le seul point à
justifier est l'existence d'un chemin étiqueté~$a  _1, \pp,
a  _k$ depuis~$1$ dans~$\G$: or, pour $i \le k$, $a  _1 \ppp
a  _i$ est un diviseur à gauche de~$x  \cR y$, donc un
élément de~$S$, et donc l'arête $(a  _1 \ppp a  _{i-1}, a  _i,
a  _1 \ppp a  _i)$ appartient bien à~$\G$.\FinDem

\Ex\«\ExMon»
La figure~1.3 montre le graphe caractéristique dans
le cas du monoïde $\Mon{a  , b  \; a   b  a   = b^2}$,
dont les critères de la section~6 montreront qu'il est un
monoïde de Garside: il y a ici 8 éléments simples. Noter
que le ppcm à droite des atomes, à savoir~$b^2$, n'est
pas l'élément~$\D$, qui est~$b^3$. 

\midinsert
$$\OpenGraphicBox width 48mm height 45mm
depth0mm;
\leavevmode\special{illustration 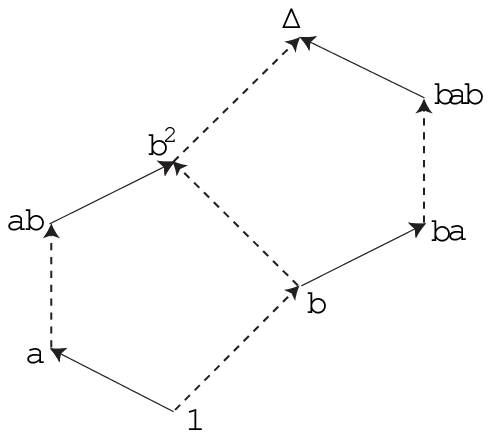}
\CloseGraphicBox$$ 
\centerline{{\bf Figure 1.3.} Graphe caractéristique du monoïde
$\Mon{a  , b  \; a   b  a   = b^2}$}
\centerline{(Les arêtes pleines représentent~$a  $,
les tiretées~$b $)} 
\endinsert

\bigskip En effaçant les étiquettes du graphe caractéristique, on obtient
le diagramme de Hasse du treillis des éléments simples. Ce graphe
ne détermine pas le monoïde à isomorphisme près. Par exemple, les
monoïdes $\Mon{a  ,
b  \; a  b  = b a  }$ et $\Mon{a  , b  \; a  ^2 =
b^2}$ sont tous deux des monoïdes de Garside; les
graphes de Cayley associés sont respectivement
$$\OpenGraphicBox width 65mm height 21mm
depth0mm;
\leavevmode\special{illustration 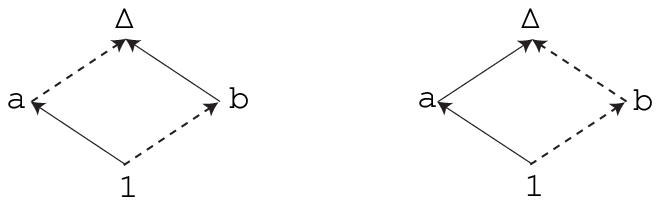}
\CloseGraphicBox$$ 
Les graphes non étiquetés sous-jacents sont identiques, alors que les
monoïdes ne sont pas isomorphes, puisque les graphes étiquetés ne le
sont pas.

Lorsque $M$ est le monoïde d'Artin associé à un groupe de
Coxeter fini~$W$, les éléments simples de~$M$ sont en bijection avec
les éléments de~$W$, et le graphe caractéristique de~$M$
s'obtient à partir du graphe de Cayley de~$W$ en orientant les
flèches, c'est-à-dire en supprimant les relations de torsion~$x^2 =
1$. Cette propriété ne s'étend pas au cas d'un monoïde de Garside
général, puisque, dans l'exemple~«\ExMon», tous les sommets du
graphe n'ont pas le même degré, ce qui exclut que ce dernier
soit la projection du graphe de Cayley d'un groupe. 

\Question
Le graphe caractéristique d'un monoïde de Garside est-il  la
projection du graphe de Cayley d'un monoïde quotient
de~$M$ obtenu en ajoutant des relations de torsion?

\bigskip Le réponse est positive dans le cas de l'exemple~«\ExMon»: si
$W$ est le monoïde obtenu en ajoutant à la présentation les relations
de torsion $a  ^2 = 1$ et $b^4 = b $ (la seconde est en fait
conséquence de la première), alors le graphe caractéristique
de~$M$ s'obtient à partir du graphe de Cayley de~$W$ en supprimant
des arêtes de torsion (Figure~1.4). Le groupe
$\Gr{a  ,
b  \; a  b a   = b^2}$ est le groupe de tresses~$B_3$, présenté à
partir des générateurs $a   = \s_1$ et $b  = \s_2\s_1$. Le
groupe~$B_3$ est également le groupe de fractions du monoïde de Garside~$B_3^+$, dont le groupe de Coxeter associé est le groupe
symétrique~$\Sym_3$. Remarquer que ce dernier admet la
présentation $\Gr{a  , b  \; a   b  a   = b^2,
a  ^2 = 1, b^3 = 1}$, et qu'il s'obtient donc à partir du monoïde
$\Mon{a  , b  \; a   b  a   = b^2, a  ^2 = 1, b^4 = b }$ en
quotientant par la relation~$b^3 = 1$.

\midinsert
$$\OpenGraphicBox width 48mm height 45mm
depth0mm;
\leavevmode\special{illustration 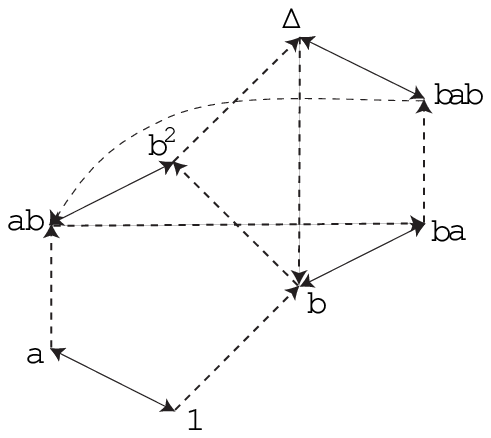}
\CloseGraphicBox$$ 
\centerline{{\bf Figure 1.4.} Graphe de Cayley
du monoïde $\Mon{a  , b  \; a   b  a   = b^2, a  ^2 =
1, b^4 = b }$}
\centerline{(Les arêtes pleines représentent~$a  $,
les tiretées~$b $)} 
\endinsert

\Sec Dualité dans les monoïdes de Garside

\bigskip L'objet de cette section est à la fois d'établir des formules de
calcul valables dans tout monoïde de Garside (lemmes~2.5
et~2.6), et de donner une nouvelle ca\-rac\-térisation de ceux-ci par
des conditions plus faibles que celles de la définition
initiale qui seront utilisées dans la section~5. Le résultat est le suivant:

\Prop\«\Char»
{\sl Un monoïde est un monoïde de Garside si, et
seulement si, il vérifie les conditions~$(C_0)$, $(C_1)$,
$(C_2)$, $(C_3)$, et $(\Til C_1)$.}

\bigskip Tout monoïde gaussien satisfait par définition aux
conditions~$(C_0)$, $(C_1)$, $(C_2)$ et $(\Til C_1)$, et tout
monoïde de Garside satisfait en outre la
condition~$(C_3)$ d'après la proposition~«\Alternative», donc les
conditions de la proposition~«\Char» sont nécessaires. Le problème est
de montrer que ces conditions sont aussi suffisantes, et, en particulier,
de montrer qu'elles entraînent l'existence des ppcm à gauche
qui n'y sont pas mentionnés, non plus que
l'atomicité. La démonstration comporte plusieurs résultats
intermédiaires, et repose sur l'existence d'une dualité
échangeant division à gauche et à droite pour les éléments
simples.

\Lem\«\Demep» 
{\sl Soit $M$ un monoïde satisfaisant aux
conditions~$(C_0)$, $(C_1)$ et~$(C_2^+)$, 
$S$ une partie génératrice de~$M$ close par~$\cR$
et~$\jR$, et admettant un ppcm à droite~$\D$.

(i) Pour tout entier~$k$, l'ensemble~$S^k$ est clos
par~$\cR$ et~$\jR$, tout diviseur à droite de~$\D^k$
appartient à~$S^k$, et tout élément de~$S^k$ est un
diviseur à gauche de~$\D^k$.

(ii) Posons $x^* = x  \cR \D$ pour~$x  \in S$. Alors,
on a $x^* \in S$ et $x  x^* = \D$. Pour $x , y$
dans~$S$, $y^*$ est un diviseur à droite de~$x^*$ si
$x$ est un diviseur à gauche de~$y$, et l'implication
réciproque est vraie si $M$ admet la simplification à droite.

(iii) La fonction $x  \mapsto x^{**}$ s'étend en un
endomorphisme~$\f$ de~$M$ qui envoie~$S^k$ dans
lui-même pour tout~$k$, et, pour tout~$x$ dans~$M$, on
a
$$x  \op \D = \D \op \f(x ). \Eq\«111»$$

}\Dem (i) Le lemme~«\Closure»(iii) affirme que $S^k$
est clos par~$\cR$ et~$\jR$, et le lemme~«\Compl» que tout
diviseur à droite de~$\D^k$ appartient à~$S^k$, puisque
$\D^k$ est élément de~$S^k$ par construction.

Supposons $x  \in S^k$. On démontre par récurrence
sur~$k$ que $x$ est un diviseur à gauche de~$\D^k$. Pour
$k = 0$, on a $x  = 1$, et le résultat est trivial. Pour $k = 1$,
le résultat découle de la définition de~$\D$. Supposons
$x  \in S^k$ avec $k \ge 2$. On écrit $x  =
x_1 x '$ avec $x_1 \in S$ et $x ' \in S^{k-1}$.
Par~«103», on a
$$\D^k \cR x  = (\D^k \cR x_1) \cdot ((x_1 \cR
\D^k) \cR x ').$$ 
Puisque $x_1$ est dans~$S$, il est diviseur à gauche
de~$\D$, donc de~$\D^k$, donc on a $\D^k \cR x_1 = 1$, et,
d'autre part, toujours par~«103»,
$$x_1 \cR \D^k = (x_1 \cR \D) \om ((\D \cR x_1)
\cR \D^{k-1}) = (x_1 \cR \D) \om (1 \cR \D^{k-1}) = (x_1
\cR \D) \om \D^{k-1}.$$
On déduit
$$\D^k \cR x  = ((x_1 \cR \D)  \op \D^{k-1}) \cR x '
= \D^{k-1} \cR ( (x_1 \cR \D) \cR x ' ).$$
Par le lemme~«\Compl»,  $(x_1 \cR \D) \cR x '$ est
dans~$S^{k-1}$, donc, par hypothèse de récurrence, cet
élément est un diviseur à gauche de~$\D^{k-1}$. On déduit
$ \D^{k-1} \cR ( (x_1 \cR \D) \cR x ' ) = 1$, soit
finalement $\D^k \cR x  = 1$, et $x$ est
diviseur à gauche de~$\D^k$.

(ii) Supposons $x  \in S$. Comme $\D$ est dans~$S$, on a
$x^* = x  \cR \D \in S$ puisque, par hypothèse, $S$ est
clos par~$\cR$. Par définition, $x$ est diviseur à gauche
de~$\D$, donc on a $x  x^* = x  \jR \D = \D$. Supposons
que $x$ et $y$ soient dans~$S$ et $x$ est un diviseur
à gauche de~$y$, disons $y  =x z $. On a alors $x 
x^*= \D = y  y^* = x  z  y^*$, d'où
$x^* = z  y^*$, et $y^*$ est diviseur à droite
de~$x^*$. Inversement, $x^* = z  y^*$ entraîne $y 
y^* = x  z  y^*$, donc $x  z  = y$ si la
simplification à droite par~$y^*$ est permise.

(iii) Supposons $x  \in S$. Alors $x^*$ est dans~$S$,
donc $x^{**}$ est défini, et on obtient
$$x  \op \D = x  \op (x^* \op x^{**}) = (x  \op x^*)
\op x^{**} = \D \op x^{**}. \Eq\«98»$$
Soit $x$ un élément quelconque de~$M$. Puisque $S$
engendre~$M$, il existe une décomposition, en général non
unique, de~$x$ comme produit d'éléments de~$S$.
Supposons $x  = x_1 \ppp x_p = y_1 \ppp y_q$, 
avec $x_1, \pp, x_p, y_1, \pp, y_q \in S$. Par~«98»,
on obtient
$$\D \op x_1^{**} \ppp x_p^{**}
= x_1 \ppp x_p \op \D 
= x  \op \D 
= y_1 \ppp y_q \op \D 
= \D \op y_1^{**} \ppp y_q^{**},$$
et, puisque $\D$ est simplifiable à gauche, il n'y a pas
d'ambiguïté à définir~$\f$ par $\f(x ) = x_1^{**} \ppp
x_k^{**}$. Alors, par construction, $\f$ est un endomorphisme
de~$M$, et «111» est vérifiée pour tout~$x$ dans~$M$. Enfin, $x 
\in S^k$ entraîne
$\f(x )
\in S^k$ par construction.\FinDem

\bigskip Nous introduisons désormais la condition de
finitude~$(C_3)$ dans les hypothèses. Notons que, même si
$M$ est un monoïde de type fini, les conditions~$(C_0)$,
$(C_1)$ et $(C_2^+)$ ne garantissent pas~$(C_3)$,
c'est-à-dire l'existence d'une partie génératrice finie close par~$\cR$: le
monoïde de Baumslag--Solitar $\Mon{a  , b  \; a  b^2 = b a  }$
est un contre-exemple, car on a alors $a   \cR
b^n = b^{2n}$ pour tout~$n$, et la clôture
par~$\cR$ de l'ensemble des atomes~$\{a  , b \}$ est infinie.

D'après le lemme~«\Closure», si $M$ est un monoïde
vérifiant les conditions~$(C_0)$, $(C_1)$, $(C_2)$ et $(C_3)$, il
existe une partie finie~$S$ close par~$\cR$ et~$\jR$. En
particulier, $S$ admet un ppcm à droite~$\D$, qui est
élément de~$S$, et $M$ vérifie les hypothèses du
lemme~«\Demep».

\Lem\«\Bijection»
{\sl Soit $M$ un monoïde vérifiant les conditions~$(C_0)$,
$(C_1)$, $(C_2)$ et~$(C_3)$, et $S$ une partie
génératrice finie de~$M$ close par~$\cR$ et~$\jR$. Soient
$\dL$, $\f$ comme dans le lemme~«\Demep». Alors les
conditions suivantes sont équivalentes:

(i) Le monoïde~$M$ satisfait~$(\Til C_1)$, c'est-à-dire admet la
simplification à droite;

(ii) L'application $x  \mapsto x^*$ est une permutation
de~$S$;

(iii) Il existe~$m$ tel que $\f^m(x ) = x$ soit vérifié
pour tout~$x$ dans~$S$;

(iv) L'application~$\f$ est un automorphisme de~$M$, et
une permutation de~$S^k$ pour tout~$k$.}

\Dem Notons $\D$ le ppcm à droite de~$S$. Puisque $S$ est
fini, (ii) et (iii) sont équivalents, et, puisque $\f$ est un
endomorphisme, que $S$ engendre~$M$, et que $\f$
envoie~$S^k$ dans lui-même pour tout~$k$, il en est de
même de~(ii) et~(iv).

Supposons $x , y  \in S$ et $x^* = y^*$. On a alors $x 
x^* = \D = y  y^* = y  x^*$, d'où $x  = y$ si $M$
admet la simplification à droite. Donc (i) entraîne
l'injectivité de~${}^*$ sur~$S$,  donc~(ii) puisque $S$ est
fini. 

Inversement, supposons $x z  = y z $. Pour $k$ assez grand,
$x$, $y$ et $z $ sont dans~$S^k$, donc $z $ est un diviseur à
gauche de~$\D^k$, et $x z  = y z $ entraîne $x  \D^k = y 
\D^k$. Par~«111», on trouve
$$\D^k \op \f^k(x ) = x  \op \D^k = y  \op \D^k = \D^k \op
\f^k(y ),$$
d'où $\f^k(x ) = \f^k(y )$, et $x  = y$ si (iv) est
vérifiée. Donc (iv) entraîne~(i).\FinDem

\bigskip Nous supposons désormais, jusqu'à la fin de cette
section, que $M$ est un monoïde satisfaisant aux
conditions~$(C_0)$, $(C_1)$, $(C_2)$, $(C_3)$, et $(\Til C_1)$,
et que $S$ est une partie génératrice finie de~$M$ close
par~$\cR$ et~$\jR$. Gardant les conventions précédentes,
nous notons $\D$ le ppcm à droite de~$S$, et, pour $x$
dans~$S$, nous posons~$x^* = x  \cR \D$, et nous
notons~$\f$ l'automorphisme de~$M$ induit par~$x 
\mapsto x^{**}$.

\Lem\«\Simple»
{\sl (i) Les implications du lemme~«\Demep»(i) sont des
équivalences: pour tout~$k$, il y a équivalence entre être
diviseur à gauche de~$\D^k$, être diviseur à droite
de~$\D^k$, et appartenir à~$S^k$.

(ii) Tout élément de~$S^k$ admet au plus $\card(S)^k$
diviseurs à gauche.

(iii) Le monoïde~$M$ est atomique, et $x  \in S^k$ entraîne
$\norm(x ) \le \card(S)^k$; en parti\-cu\-lier, on a 
$\norm(\D) \le \card(S)$.

(iv) Deux éléments quelconques de~$M$ ont un pgcd à
gauche. La structure $(M, \jR, \gL)$ est un treillis, et $(S,
\jR, \gL)$ est un treillis fini de minimum~$1$ et de
maximum~$\D$.}

\Dem (i) Par le lemme~«\Bijection», $\f$ est un
automorphisme de~$M$. Supposons que $x$ divise à
gauche~$\D^k$, soit $x z  = \D^k$. On a alors
$$\f^{-k}(z ) \op x  \op z  
= \f^{-k}(z ) \op \D^k = \D^k \op z ,$$
d'où, en simplifiant à droite, $\D^k = \f^{-k}(z ) \op x$, et
$x$ divise à droite~$\D^k$.

(ii) D'après~(i), tout diviseur à gauche d'un élément
de~$S^k$ est lui-même élément de~$S^k$, d'où le
résultat puisqu'on a $\card(S^k) \le \card(S)^k$.

(iii) Supposons $x  \in S^k$, et $x  = x_1 \ppp x_n$
avec $x_1, \pp, x_n \not= 1$. Par construction, $x_1
\ppp x_i$ est un diviseur à gauche de~$x$, donc
de~$\D^k$, pour tout~$i$, ce qui entraîne $x_1 \ppp x_i
\in S^k$. Les conditions~$(C_0)$ et~$(C_1)$ entraînent
$x_1 \ppp x_i \not= x_1 \ppp x_j$ pour $i \not= j$.
On a donc $n \le \card(S^k) \le \card(S)^k$, soit $\norm(x )
\le \card(S)^k$.

(iv) D'après le lemme~«\Simple»(ii), l'ensemble des
diviseurs à gauche de tout élément de~$M$ est fini.
Donc, pour $x$, $y$ quelconques dans~$M$, l'ensemble
des diviseurs à gauche communs de~$x$ et~$y$ est fini,
et il admet un ppcm à droite, lequel est un pgcd à gauche
de~$x$ et~$y$ par construction. Que $(M, \jR, \gL)$ soit
un treillis est alors standard. Pour~$S$, nous remarquons
que, pour $x , y  \in S$, l'ensemble des diviseurs à
gauche de~$x$ et~$y$ est inclus dans~$S$ et donc que le
ppcm à droite de cet ensemble, qui est $x  \gL y$,
appartient à~$S$.\FinDem

\bigskip L'étape suivante consiste à utiliser la dualité $x 
\mapsto x^*$ pour montrer la symétrie de la structure.
Par le lemme~«\Bijection», nous savons que l'application 
$x  \mapsto x^*$ est une permutation de~$S$.
Pour~$x$ dans~$S$, nous noterons $\dL x$ l'unique
élément de~$S$ vérifiant $(\dL x )^* = x$. Alors, pour
tout~$x$ dans~$S$, nous avons
$$\D = x  \cdot x^* = \dL x  \cdot x ,\Eq\«119»$$
et, donc, $\dL(x^*) = x$. D'après le lemme~«\Demep»(ii),
si $x$ et $y$ sont dans~$S$, $x$ est diviseur à gauche
de~$y$ si et seulement si $y^*$ est diviseur à droite
de~$x^*$, et, par conséquent, $x$ est diviseur à droite
de~$y$ si et seulement si $\dL y$ est diviseur à gauche
de~$\dL x$.

\Lem\«\Dual»
{\sl (i) Deux éléments quelconques de~$M$ admettent un
ppcm à gauche et un pgcd à droite.

(ii) L'ensemble~$S$ est clos par les opérations~$\cL$, $\jL$
et~$\gR$, et, pour $x , y  \in S$, on a
$$x  \cL y  = (\dL x  \gL \dL y ) \cR \dL y , \qquad 
x  \jL y  = (\dL x  \gL \dL y )^*, \qquad
x  \gR y  = (\dL x  \jR \dL y )^*.  \Eq\«122»$$

(iii) L'élément~$\D$ est ppcm à gauche de~$S$, et, pour
tout~$x$ dans~$S$, on a\  $\dL x  = \D \cL x$.}

\Dem (i) Soient $x$, $y$ quelconques dans~$M$. Alors il
existe~$k$ tel que $x$ et $y$ appartiennent à~$S^k$, ce
qui, par le lemme~«\Simple», entraîne que $\D^k$ est un
multiple à gauche de~$x$ et de~$y$. D'après~[\\\Dfx,
Proposition~2.4], ceci est suffisant pour assurer l'existence
d'un ppcm à gauche et d'un pgcd à droite pour~$x$
et~$y$

(ii) Par définition, $\dL x  \gL \dL y$ est un diviseur à
gauche de~$\dL x$ et~$\dL y$, donc $(\dL x  \gL \dL
y )^*$ est un multiple à gauche de~$x$ et~$y$, et donc
de~$x  \jL y$. Donc, en particulier, $x  \jL y$ est
dans~$S$, et $\dL(x  \jR y )$ est défini. Ensuite, $x$ et
$y$ sont des diviseurs à droite de~$x  \jL y$, donc
$\dL(x  \jL y )$ est un diviseur à gauche de~$\dL x$ et
de~$\dL y$, donc de~$\dL x  \gL \dL y$. Par
conséquent, $(\dL x  \gL \dL y )^*$ est un diviseur à
droite de~$x  \jL y$, et nous déduisons $x 
\jL y  = (\dL x  \gL \dL y )^*$.

On trouve ensuite
$$\dL y  \op y  = \D = \dL(x  \jL y ) \op (x  \jL y )  =
(\dL x  \gL \dL y ) \op (x  \cL y ) \op y ,$$ 
d'où $\dL y  = (\dL x  \gL \dL y ) \op (x  \cL y )$, qui
donne $x  \cL y  = (\dL x  \gL \dL y ) \cR \dL y$.

L'argument est semblable pour le pgcd à droite. Comme
$\dL x$ et~$\dL y$ sont des diviseurs à gauche de
$\dL x  \jR \dL y$,  $(\dL x  \jR \dL y )^*$ est diviseur à
droite de~$x$ et de~$y$, donc de~$x  \gR y$. D'un
autre côté,  $x  \gR y$ est un diviseur à droite d'un
élément de~$S$, donc il est élément de~$S$. Puisque $x 
\gR y$ est diviseur à droite de~$x$ et~$y$, $\dL x$ et
$\dL y$ sont diviseurs à gauche de $\dL(x  \gR y )$,
donc $\dL x  \jR \dL y$ est diviseur à gauche de $\dL(x 
\gR y )$, et $x  \gR y$ est diviseur à droite de $(\dL x 
\jR \dL y )^*$. 

(iii) On a déjà noté que $\D$ est multiple à gauche de tout
élément de~$S$. Puisque $\D$ appartient à~$S$, ceci
entraîne que $\D$ est ppcm à gauche de~$S$. Soit $x$
quelconque dans~$S$: on a alors $\dL x  x   = \dL x 
(\dL x )^*= \D = \D \jL x  = (\D \cL x ) x$, d'où $\dL x  =
\D \cL x$.\FinDem

\bigskip Nous avons ainsi complété la démonstration de la
proposition~«\Char». En effet, il s'agissait de montrer qu'un
monoïde~$M$ comme ci-dessus est un monoïde de Garside. Or, d'après le lemme~«\Simple»(iii), $M$ est
atomique; par définition, deux éléments quelconques
de~$M$ admettent un ppcm à droite; ils admettent un pgcd
à gauche par le lemme~«\Simple»(iii), et un ppcm à gauche
et un pgcd à droite par le lemme~«\Dual», donc $M$ est un
monoïde gaussien. Enfin, la partie~$S$ considérée
ci-dessus est finie, elle engendre~$M$, et on a vu qu'elle est
close par chacune des opérations~$\cR$, $\cL$, $\jR$
et~$\jL$. Par la proposition~«\Alternative», on déduit que
$M$ est un monoïde de Garside.\FinDem

\bigskip Les formules de dualité~«122» sont valables dans
tout monoïde de Garside, et elles déterminent
les opérations~$\cL$, $\jL$ et~$\gR$ en termes de leurs
contreparties~$\cR$, $\jR$ et~$\gL$. Ces formules ne
s'appliquent qu'aux éléments de l'ensemble~$S$ considéré, ce
qui n'est pas une restriction véritable: pour des
éléments~$x$, $y$ arbitraires, on peut déterminer les
valeurs de~$x  \cL
y$, $x  \jL y$ et $x  \gR y$ soit en utilisant une
décomposition de~$x$ et~$y$ en produit d'éléments
de~$S$, soit --- ce qui est essentiellement équivalent ---
en remplaçant~$S$ par une puissance~$S^k$ telle que $x$
et $y$ appartiennent à~$S^k$; on sait alors que
l'ensemble~$S^k$ a la mêmes propriétés de clôture que~$S$.
Nous terminerons cette section par une formule générale
exprimant le pgcd de deux éléments arbitraires en termes des
opérations de complément~$\cR$ et~$\cL$.

\Lem\«\Gcd»
{\sl Soit $M$ un monoïde gaussien. Alors,
pour tous~$x , y$ dans~$M$, on a
$$x  \jR y  = (x  \gL y ) \om ((x  \cR y ) \jL (y  \cR x )),\Eq\«120»$$
et, donc, $x  \gL y  = (x  \jR y ) \cL ((x  \cR y ) \jL (y  \cR x ))$.}

\Dem  Posons $x ' = y  \cR x$, $y ' = x  \cR y$, $x '' = x ' \cL y '$,
et $y '' = y ' \cL x '$. Par définition, nous avons
$$x  \jR y  = x  y ' = y  x ' 
\hbox{\quad et \quad} x ' \jL y ' = x '' y ' = y '' x '.$$
D'après la première égalité, $x  \jR y$ est un multiple à
gauche de~$x '$ et~$y '$, donc il existe~$z $ vérifiant
$x  \jR y  = z  (x ' \jL y ')$. On déduit $x y ' = z 
x '' y '$, donc $x  = z  x ''$, et, de même, $y  =
z y ''$. Donc $z $ est un diviseur à gauche de~$x$
et~$y$, donc de~$x  \gL y$. 

Réciproquement, supposons $x  = z_1 x_1$ et $y  =
z_1 y_1$. On a $z_1 x_1 y ' = z_1 y_1 x ' = x 
\jR y$, donc $x_1 y ' = y_1 x '$. Il existe donc $z ''$
vérifiant $x_1 = z '' x ''$ et $y_1 = z '' y ''$. On
déduit $z_1 z '' x '' y ' =z_1 z '' y '' x ' = x  \jR
y$, d'où $z_1z '' = z $. Par conséquent, $z_1$ est un
diviseur à gauche de~$z $, et on conclut que $z $ est le
pgcd à gauche de~$x$ et~$y$.\FinDem

\Sec Structure automatique

\bigskip Nous établissons ici que tout groupe de Garside est
bi-automatique (Théorème~A de l'introduction), et, de surcroît, que
cette structure automatique s'explicite très simplement en termes des
opérations de complément.

Dans un premier temps, nous allons construire une forme
normale dans tout groupe gaussien (de Garside ou non).
Soit $M$ un monoïde quelconque, et $S$ une
partie génératrice de~$M$. Alors tout élément de~$M$
s'écrit comme produit fini d'éléments de~$S$. L'idée pour obtenir une
décomposition distinguée, classique pour les groupes de tresses
depuis  [\\\Dlg] [\\\Adk] [\\\Eps] [\\\ElM], consiste à pousser
les éléments de~$S$ par exemple vers la gauche, de sorte que le
premier élément soit maximal. 

\Def
Soit $M$ un monoïde, et $S$ une partie de~$M$. Pour $x ,
y  \in M$, on dira que $x  \perp_S y$ est vérifié si, pour
tout diviseur à gauche~$s $ de~$y$ appartenant à~$S
\setminus \{1\}$, on a $x s  \notin S$.

\bigskip Le point crucial de la construction est le résultat
suivant, qui est une conséquence immédiate des règles
du calcul des compléments:

\Lem\«\Main»
{\sl Soit $M$ un monoïde vérifiant les
conditions~$(C_0)$, $(C_1)$ et $(C_2)$, et $S$ une
partie de~$M$ close par~$\cR$ et~$\jR$. Soit
$(x_1, \pp, x_p)$ une suite d'éléments de~$S$ vérifiant
$x_i \perp_S x_{i+1}$ pour tout~$i$.
Alors on a $x_1 \perp_S x_2 \ppp x_p$.}

\Dem (Figure~3.1) Supposons  $x_1 s_1 \in S$,
avec $s_1 \in S$  et $s_1 \not=1$. Il s'agit de montrer que
$s_1$ n'est pas un diviseur à gauche de~$x_2 \ppp
x_p$. Posons de proche en proche $s_i = x_i \cR s_{i-1}$
pour $2 \le i \le p$, et montrons par récurrence sur~$i$
qu'on a $x_i s_i \in S$, $s_i \in S$ et $s_i \not= 1$. Pour $i
= 1$, ce sont les hypothèses posées plus haut. Pour $i \ge
2$, on a
$s_i = x_i \cR s_{i-1}$ et $x_i s_i = x_i \jR s_{i-1}$, d'où $s_i
\in S$ et $x_i s_i \in S$, puisque $x_i$ est dans~$S$, $s_{i-1}$
aussi par hypothèse de récurrence, et $S$ est clos
par~$\cR$ et~$\jR$; par ailleurs, la conjonction de $x_{i-1}
s_{i-1} \in S$ et $x_{i-1} \perp_S x_i$ entraîne que $s_{i-1}$
n'est pas un diviseur à gauche de~$x_i$, c'est-à-dire que $x_i \cR
s_{i-1}$, qui est~$s_i$, ne vaut pas~$1$. On a donc $s_p \not=
1$. Or, par les règles du lemme~«\Rules», on a $s_p = (x_2 \ppp
x_p) \cR s_1$, et
$s_p \not= 1$ signifie donc que $s_1$ n'est pas un diviseur
à gauche de~$x_2 \ppp x_p$.\FinDem

\midinsert
$$\OpenGraphicBox width 67mm height 18mm depth0mm;
\leavevmode\special{illustration 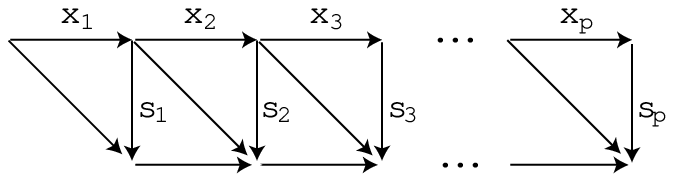}
\CloseGraphicBox$$
\centerline{{\bf Figure 3.1.} Décomposition normale}
\endinsert

\Prop\«\Decomposition»
{\sl Soit $M$ un monoïde gaussien --- ou,
plus généralement, un monoïde vérifiant $(C_0^+)$,
$(C_1)$, et $(C_2^+)$ --- et $S$ une partie
génératrice de~$M$ close par~$\cR$ et~$\jR$. Alors tout
élément de~$M$ possède une unique
décomposition~$x_1 \ppp x_p$ avec $x_1,
\pp, x_p \in S$ et $x_i \perp_S x_{i+1}$
pour tout~$i$.}

\Dem On prouve d'abord l'existence d'une
décomposition comme ci-dessus pour tout élément~$x$
de~$M$ par récurrence sur la norme~$\norm(x )$. Pour
$\norm(x ) = 0$, on a $x  = 1$, et le résultat est vrai.
Supposons $x  \not= 1$. Notons $\DivL(x )$ l'ensemble
des diviseurs à gauche de~$x$. Soit $x_1$ un élément
de norme maximale dans~$\DivL(x )$. Comme $y  \in
\DivL(x )$ entraîne $\norm(y ) \le \norm(x )$, un tel
élément~$x_1$ existe, et, comme $\DivL(x )$ contient au
moins un atome, on a $x_1
\not= 1$. Ecrivons $x  = x_1 x '$. On a alors $\norm(x ') \le
\norm(x ) - \norm(x_1) < \norm(x )$, donc, par hypothèse
de récurrence, $x '$ possède une décomposition $x ' = x_2
\ppp x_p$ avec $x_i \in S$ et $x_i \perp_S x_{i+1}$ pour~$i
\ge 2$. On a alors $x  = x_1 x_2 \ppp x_p$, et il s'agit de
montrer $x_1 \perp_Sx_2$. Soit $x$ un diviseur à gauche
propre de~$x_2$ appartenant à~$S$. Alors $x_1 x$ est un
diviseur à gauche de~$x_1 x_2$, donc de~$x$, et on a
$\norm(x_1 x) \ge \norm(x_1) + \norm(x) > \norm(x_1)$,
ce qui, par définition de~$x_1$, entraîne $x_1 x \notin S$.

Pour l'unicité, il suffit de montrer que, si $(x_1, \pp, x_p)$
est une suite d'éléments de~$S$ vérifiant $x_i \perp_S x_{i+1}$
pour~$1 \le i < p$, alors $x_1$ est déterminé par le produit~$x_1
\ppp x_p$. En effet, puisque $M$ est simplifiable à gauche,
une récurrence montre ensuite que $x_2$, \pp, $x_p$ sont
déterminés de même. Posons $x  = x_1 \ppp x_p$. Par
construction, $x_1$ appartient à $\DivL(x ) \cap S$. Supposons
que $x_1s $ soit un diviseur à gauche de~$x$: alors
$s $ est un diviseur à gauche de~$x_2 \ppp x_p$, par le
lemme~«\Main», on a $x_1 \perp_Sx_2 \pp
x_p$, et donc on déduit $x_1s  \notin S$. Ceci signifie
que $x_1$ est un élément de norme maximale
dans~$\DivL(x ) \cap S$, ce qui détermine $x_1$ de
façon unique, car $\DivL(x ) \cap S$ est clos par ppcm à
droite et, si $x_1$ et $x '_1$ étaient deux éléments
distincts de norme maximale dans $\DivL(x ) \cap S$,
l'élément $x_1 \jR x '_1$ contredirait cette
maximalité.\FinDem

\bigskip On considère maintenant le cas des monoïdes de Garside.
Supposons que $M$ soit un monoïde de Garside. Nous
savons que $M$ admet une unique partie génératrice
minimale, à savoir l'ensemble~$A$ de ses atomes, et que la
clôture~$S$ de~$A$ par~$\cR$ et~$\jR$ est la plus petite
partie génératrice de~$M$ qui soit close par ces opérations.
Notons que $S$ est également la clôture de~$A$ par les
opérations à gauche~$\cL$ et~$\jL$. En effet, d'après le
lemme~«\Dual»(ii), l'ensemble~$S$ est clos par~$\cL$
et~$\jL$, donc, si $\Til S$ désigne la clôture de~$A$
par~$\cL$ et~$\jL$, nous avons $\Til S \ince S$. Comme les
hypothèses sont symétriques, un argument similaire donne $S
\ince
\Til S$, d'où $S = \Til S$. On note~$\D$ le ppcm de~$S$, et les
opérations de dualité réfèrent au complément dans~$\D$: pour~$x$
dans~$S$, on pose $x^* = x  \cR \D$, et $\dL x  = \D \cL x$. 

Nous allons appliquer la proposition~«\Decomposition» en
prenant pour~$S$ l'ensemble des éléments simples. En fait, la
minimalité n'est requise nulle part, et nous pourrions aussi bien
utiliser à la place de~$S$ une partie génératrice quelconque qui soit
close pour~$\cR$ et~$\jR$, en particulier toute partie~$S^k$ avec $k
\ge 1$.

La condition d'orthogonalité mentionnée dans la
construction s'exprime simplement en termes de pgcd.

\Lem\«\Orthogonal»
{\sl Soit $M$ un monoïde de Garside, et $S$
l'ensemble de ses éléments simples. Pour $x , y$ dans~$S$, $x 
\perp_Sy$ équivaut à $x^* \gL y  = 1$.}

\Dem On a $x  \perp_S y$ si et seulement si, pour tout diviseur à
gauche simple~$s $ de~$y$ distinct de~$1$, $x s $ n'est pas
simple, c'est-à-dire que $x s $ n'est pas un diviseur à gauche
de~$\D$ de~$S$, soit encore que $s $ n'est pas un diviseur à
gauche de~$x  \cR \D$, qui est~$x^*$.\FinDem

\Prop\«\Normal»
{\sl Soit $M$ un monoïde de Garside.
Tout élément~$x$ de~$M$ admet une unique
décomposition de la forme~$x_1 \ppp x_p$ avec
$x_1, \pp, x_p$ simples distincts de~$1$ et $x_i^* \gL
x_{i+1} = 1$ pour tout~$i$. Pour $x  \not= 1$, le premier
facteur de la décomposition est ~$x  \gL \D$.}

\Dem D'après la proposition~«\Decomposition» et le
lemme~«\Orthogonal», seule la dernière assertion reste à
démontrer. Supposons $x  \not= 1$, et soit $x_1 \ppp
x_p$ la décomposition de~$x$ donnée par la
proposition. Par construction, $x_1$ est un diviseur à
gauche simple de~$x$ de norme maximale. Puisqu'il est
simple, $x_1$ est un diviseur à gauche de~$\D$, donc
de~$x  \gL \D$, et, de là, nous avons $\norm(x_1) \le
\norm(x  \gL \D)$. D'un autre côté, $x  \gL \D$ est
un diviseur à gauche simple de~$x$, donc, par définition
de~$x_1$, on doit avoir $\norm(x_1) = \norm(x  \gL
\D)$, d'où $x_1 = x  \gL \D$  puisque $x_1$ est
un diviseur à gauche de~$x  \gL \D$. \FinDem

\bigskip La forme normale précédente sera appelée forme
normale à gauche pour~$M$. De façon symétrique, nous
avons une forme normale à droite:

\Prop\«\Normal»
{\sl Soit $M$ un monoïde de Garside. Tout
élément~$x$ de~$M$ admet une unique décomposition
de la forme~$x_q \ppp x_1$ avec $x_1, \pp, x_q$ simples et
$\dL x_i \gR x_{i+1} = 1$ pour tout~$i$. Pour $x  \not=
1$, le dernier facteur de la décomposition est ~$x 
\gR \D$.}

\bigskip Considérons maintenant le cas des groupes de Garside.
 Il se ramène à celui des monoïdes grâce à
l'existence de fractions irréductibles.

\Lem\«\Fraction»
{\sl Soit $M$ un monoïde de Garside, et $G$ son groupe
de fractions. Alors tout élément~$z $ de~$G$ admet une
unique décomposition $z  = x \ii y$ avec
$x , y  \in M$ et $x  \gL y  = 1$, et, symétriquement, une
unique décomposition $z  = x ' {y '}\ii$ avec $x ', y ' \in M$ et $x '
\gR y ' = 1$.}

\Dem Par définition, tout élément~$z $ de~$G$ admet une
décomposition $z  = {x '}\ii y '$ avec $x ' , y ' \in M$. 
Posant $x ' = z ' x$ et $y ' = z ' y$ avec $z ' = x '
\gL y '$, on a $x  \gL y  = 1$ et $z  = x \ii y$ par
construction. Par le lemme~«\Gcd», on obtient
$x  \jR y  = x  (x  \cR y ) = y  (y  \cR x ) = (x  \cR
y ) \jL (y  \cR x )$.  Supposons alors $z  = {x ''}\ii
y ''$ avec $x '', y '' \in M$. De $x \ii y  = (x  \cR y )
(y  \cR x )\ii$ on déduit $x '' (x  \cR y ) = y '' (y 
\cR x )$, donc, puisque $x  (x  \cR y )$  est le ppcm à
gauche de $x  \cR y$ et $y  \cR x$, l'existence
de~$z ''$ vérifiant $x '' = z '' x$ et
$y '' = z '' y$. Alors $x '' \gL y '' = 1$ entraîne $z '' =
1$, soit $x '' = x$ et $y '' = y$.\FinDem

\Def
Dans le contexte du lemme précédent, les éléments~$x$
et~$y$ seront appelés le {\it dénominateur} et
le {\it numérateur gauches} de~$z $, et notés~$D (z )$
et~$N (z )$; symétriquement, $x '$ et~$y '$
seront appelés numérateurs et dénominateurs à droite, et
notés~$\NR(z )$ et~$\DR(z )$. Ainsi, pour tout~$z $ dans~$G$,
on a
$$z  = D (z )\ii N (z ) = \NR(z ) \DR(z )\ii,$$
avec $D (z ) \gL N (z ) = \NR(z ) \gR \DR(z ) = 1$.

\bigskip En rassemblant les éléments, nous obtenons
l'existence des formes normales appelées \<<mixed canonical
forms\>> dans~[\\\Eps, Chap.~9]:

\Prop
{\sl Soit $G$ le groupe de fractions d'un monoïde de Garside~$M$.

(i) Tout élément de~$G$ admet une unique
décomposition de la forme $x_p\ii \ppp x_1\ii y_1 \ppp
y_q$, avec $x_1, \pp, x_p, y_1, \pp, y_q$ éléments
simples de~$M$ vérifiant
$x_1 \gL y_1 = 1$,  et, pour tout~$k$, $x_k^* \gL
x_{k+1} = 1$, et $y_k^* \gL y_{k+1} = 1$.

(ii)  Tout élément de~$G$ admet une unique
décomposition de la forme $x_p \ppp x_1 y_1\ii \ppp
y_q\ii$, avec $x_1, \pp, x_p, y_1, \pp, y_q$ éléments
simples de~$M$ vérifiant $x_1 \gR y_1
= 1$, et, pour tout~$k$, ${}^*x_k \gR x_{k+1} = 1$, et
${}^*y_k \gR y_{k+1} = 1$.}

\bigskip Dans le contexte précédent, nous dirons que la suite
$(x_p\ii, \pp,  x_1\ii , y_1, \pp,  y_q)$ est la forme
normale à gauche de l'élément~$z $, et que $(x_p, \pp,
x_1, y_1\ii, \pp, y_q\ii)$ en est la forme normale à
droite.

Nous allons montrer que les formes normales précédentes sont
associées à une structure (bi)-automatique. La démonstration est
extrêmement simple, qui contraste avec les calculs de~[\\\Chb], repris
dans~[\\\Dfx]. Les seuls résultats dont nous aurons besoin
sont les suivants:

\Lem\«\Prelim»
{\sl Soit $M$ un monoïde gaussien.

(i) Supposons que $x '$ divise à gauche~$z  x$, et $y '$
divise à gauche~$z  y$. Alors $x ' \gL y '$ divise à
gauche $z  (x  \gL y )$; donc, en particulier, $x  \gL y 
= 1$ entraîne $x ' \gL y ' = x ' \gL y ' \gL z $.

(ii) Si $\D$ est un élément de Garside de~$M$, on a $(x y )
\gL \D = (x (y  \gL \D)) \gL \D$ pour tous~$x, y$
dans~$M$.}

\Dem (i) Par définition, $x ' \gL y '$ divise à
gauche $z x$ et $z y$, donc il divise leur pgcd à
gauche, qui est $z (x  \gL y )$. Pour $x  \gL y  = 1$, on
déduit que $x ' \gL y '$ divise à gauche $z $, donc
$x ' \gL y ' \gL z $.

(ii) Posons $y ' = y  \gL \D$. Alors $(x y )ý
\gL \D$ est le ppcm à droite des diviseurs de~$\D$ qui
divisent à gauche~$x y$, soit $(x y ) \gL \D = \JR \{s \; \D
\cR s = (x y ) \cR s = 1\}$. On a $(x y ) \cR s = y  \cR
(x  \cR s)$; or, pour $s$ divisant~$\D$, $x  \cR s$
divise~$\D$, donc la condition $y  \cR (x  \cR s) = 1$ dans
la formule précédente est équivalente à $y ' \cR (x  \cR s) =
1$, d'où
$$x y  \gL \D 
= \JR \{s \; \D \cR s = y  \cR (x  \cR s) = 1\}
= \JR \{s \; \D \cR s = y ' \cR (x  \cR s) = 1\}
= x y ' \gL \D. \eqno{\Block}$$

\bigskip Il est alors très facile de calculer, pour $s $ simple, la
forme normale d'un élément~$z s ^{\pm 1}$ à partir de celle
de~$z $. Dans le cas de la forme normale à gauche, il est
commode de décrire le passage de la forme normale de~$z $ à
celle de~$z  s \ii$, et à celle de~$z  \D$:  le passage
de~$z $ à~$z s $ s'en déduit en écrivant $z  s  = z  \D
(s ^*)\ii$.

\Lem\«\Normale»
{\sl Soit $G$ le groupe de fractions d'un monoïde de Garside~$M$, et $(y_q\ii, \pp, y_1\ii,
x_1, \pp, x_p)$ la forme normale à gauche d'un
élément~$z $ de~$G$. 

(i) Pour $s $ simple, la forme normale à gauche
de~$z s \ii$ est $({y '}_{\!\!q+1}\ii, \pp, y '_1{}\ii,
x '_1, \pp, x '_p)$, avec $s_{p+1} = s $, puis $x '_i = x_i \cL
s_{i+1}$ et $s_i = s_{i+1} \cL x_i$ pour
$p \ge i \ge 1$, ensuite $t_1 = s_1$, $y '_j = t_jy_j \gL
\D$ et $t_{j+1} = y_j \cR (t_jy_j)$ pour $1 \le j
\le q$, et, enfin, $y '_{q+1} = t_{q+1}$.

(ii) La forme normale à gauche de~$z \D$ est
$(y_q\ii, \pp, y_2\ii, y_1^*, x_1^{**}, \pp,
x_p^{**})$.

(iii) La forme normale à gauche de~$z \D\ii$ est
$(y_q\ii, \pp, y_1\ii, ({}^*x_1)\ii, {}^{**}x_2, \pp,
{}^{**}x_p)$.}

\Dem (i) Par construction, on a $x '_i s_{i+1} = s_i
x_i$ pour tout~$i$, et $y '_j t_{j+1} = t_j
y_j$ pour tout~$j$, donc le diagramme
$$\OpenGraphicBox width 102mm height 19mm depth0mm;
\leavevmode\special{illustration 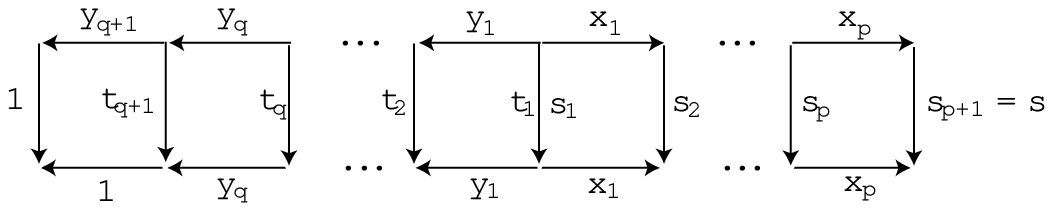}
\CloseGraphicBox$$ 
\vglue-3mm\noindent
commute, et on obtient ${y '}_{\!\!q+1}\ii
\ppp y '_1{}\ii x '_1 \ppp x '_p =  y_{q}\ii \ppp
y_1\ii x_1 \ppp x_p s \ii = z  s \ii$. La seule
question est de montrer que la suite  $({y '}_{\!\!q+1}\ii, \pp,
y '_1{}\ii, x '_1, \pp, x '_p)$ est bien une forme normale à
gauche. Par construction, tous les facteurs sont simples, et il
s'agit de vérifier les conditions de pgcd. 

Considérons d'abord le cas de~$x '_i$ et~$x '_{i+1}$.
Par définition, on a 
$$x '_i x '_i{}^* s_i^{**} 
= \D s_i^{**} 
= s_i \D
= s_i x_i x_i^*
= x '_i s_{i+1} x_i^*,$$
donc $x '_i{}^* s_i^{**} = s_{i+1} x_i^*$,
ce qui montre que $x '_i{}^*$ est un diviseur à
gauche de~$s_{i+1} x_i^*$. Par ailleurs, on a $x '_{i+1}
s_{i+2} = s_{i+1} x_i$, donc $x '_{i+1}$ est un diviseur à
gauche de~$s_{i+1} x_{i+1}$. Comme on a $s_{i+1} \gL
x '_{i+1} = 1$ par construction, et $x_i^* \gL x_{i+1} =
1$ par hypothèse, le lemme~«\Prelim»(i) donne $x '_i{}^* \gL
x '_{i+1} = 1$, comme souhaité.

Considérons maintenant le cas de~$y '_1$ et~$x '_1$. A
nouveau, $x '_1$ est diviseur à gauche de~$s_1 x_1$, et,
de même, $y '_1$ est diviseur à gauche de~$s_1 y_1$.
On a $s_1 \gL x '_1 = 1$ par construction, et $x_1 \gL
y_1 = 1$ par hypothèse, d'où $x '_1 \gL y '_1 = 1$ par le
lemme~«\Prelim»(i).

Considérons finalement le cas de~$y '_j$ et~$y '_{j+1}$. Par
définition, on a $y '_j = t_j y_j \gL \D$. Or, par
hypothèse, on a aussi $y_j = (y_j \ppp y_q) \gL \D$. Par
le lemme~«\Prelim»(ii), on  déduit $y '_j = (t_j y_j \ppp
y_q) \gL \D = (y '_j y '_{j+1} \ppp y '_{q+1}) \gL \D$,
donc, en particulier, $y '_j{}^* \gL (y '_{j+1} \ppp
y '_{q+1}) = 1$, et, a fortiori, $y '_j{}^* \gL y '_{j+1} =1$.

(ii) On a $x_i \D = \D x_i^{**}$ pour tout~$i$, et,
par ailleurs, $y_1 y_1^* = \D$, donc le diagramme
$$\OpenGraphicBox width 92mm height 19mm depth0mm;
\leavevmode\special{illustration 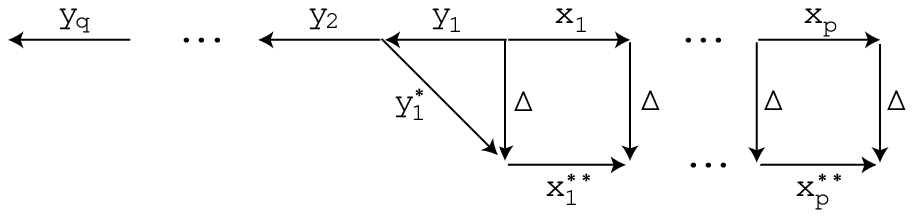}
\CloseGraphicBox$$ 
\vglue-3mm\noindent
est commutatif, et on obtient $y_q\ii \ppp y_2\ii
y_1^* x_1^{**} \ppp x_p^{**} = z  \D$. La seule
question est à nouveau de montrer que la suite
$(y_q\ii, \pp, y_2\ii, y_1^*, x_1^{**}, \pp,
x_p^{**})$ est une forme normale à gauche, soit, les facteurs
étant clairement simples, de montrer que les conditions de pgcd
entre deux termes consécutifs sont vérifiées. Or, pour tout~$i$,
on a $(x_i^{**})^* \gL x_{i+1}^{**} = (x_i^* \gL
x_{i+1})^{**} = 1^{**} = 1$, puisque l'application $x 
\mapsto x^{**}$ induit un automorphisme de~$M$. Ensuite,
on a $y_1^* \gL y_2 = 1$ par hypothèse, puis, pour $j \ge
2$, $y_j^* \gL y_{j+1} = 1$, toujours par hypothèse.

(iii) L'argument est similaire. On a cette fois ${}^{**}x_i \D =
\D x_i$ pour tout~$i$, et ${}^*x_1 x_1 = \D$, d'où le
diagramme commutatif
$$\OpenGraphicBox width 92mm height 19mm depth0mm;
\leavevmode\special{illustration 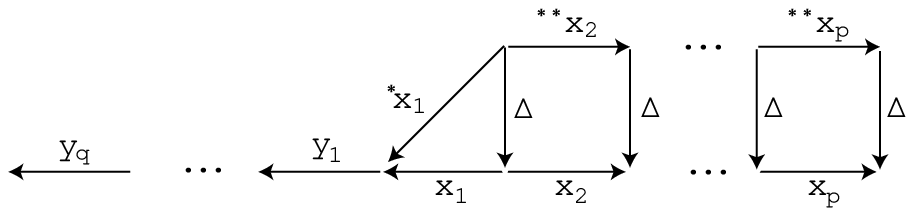}
\CloseGraphicBox$$ 
\vglue-3mm\noindent
qui donne $y_q\ii \ppp y_1\ii ({}^*x_1)\ii{}^{**}x_2
\ppp {}^{**}x_p = z  \D\ii$. Il s'agit encore de vérifier les
conditions de pgcd entre termes voisins: pour tout~$i$, on a 
$({}^{**}x_i)^* \gL {}^{**}x_{i+1} = {}^{**}(x_i^* \gL
x_{i+1}) = {}^{**}1 = 1$, et, enfin, $({}^*x_1)^* \gL y_1
= x_1 \gL y_1 = 1$.\FinDem

\Prop
{\sl Tout groupe de Garside est bi-automatique.}

\Dem Soit $M$ un monoïde de Garside, et $G$ son
groupe de fractions. Soit $S$ l'ensemble des éléments simples
de~$M$. Considérons le langage~$L$ formé par les formes
normales à gauche, considérées comme mots sur
l'alphabet~$S \cup S\ii$. D'abord $L$ est un langage régulier,
car l'appartenance d'un mot à~$L$ est définie par des
conditions locales consistant en une liste (finie) de lettres
permises après chaque lettre: après une lettre
négative~$y \ii$, les lettres autorisées sont les lettres
positives~$x$ vérifiant $x  \gL y  = 1$, et les lettres
négatives~$y '{}\ii$ vérifiant $y '{}^* \gL y  = 1$; après
une lettre positive~$x$, les lettres autorisées sont les lettres
positives~$x '$ vérifiant $x^* \gL x ' = 1$.

Ensuite la formule explicite du lemme~«\Normale»(i) montre
que la distance entre la forme normale d'un élément~$z $ et
celle de~$z  s \ii$ est uniformément bornée par~$2$ (en
termes de l'alphabet~$S$), ce qui établit la propriété du
$2$-compagnon de route --- ou \<<$2$-fellow traveller
property\>>. On déduit que le langage des formes normales à
gauche est associé à une structure automatique.

Enfin, la symétrie de la forme normale entraîne une propriété
de compagnon de route similaire pour la multiplication à
gauche, et donc le langage~$L$ est associé à une structure
bi-automatique. \FinDem

\bigskip A côté du résultat d'existence précédent, qui ne
décrit pas explicitement les automates mis en jeu, nous
allons maintenant construire des automates très simples calculant les
formes normales. Nous considérerons ici seulement le cas des
monoïdes. Le cas des groupes n'est pas directement couvert par cette
approche, mais les algorithmes décrits à la section~4
montrent comment se ramener à des fractions
irréductibles, et, de là, du groupe au monoïde.

Le cas de la forme normale à droite est le plus facile, et la
solution est une extension naturelle de celle décrite dans 
[\\\Eps, chapitre~9] pour les groupes de tresses. Nous allons
décrire un automate fini d'alphabet~$A$ et d'ensemble
d'états~$S$ calculant le dernier facteur de la forme normale à
droite, au sens où l'état final obtenu après lecture d'un mot~$u$
par l'automate est $\cl u \gR \D$, où $\cl u$ désigne la classe du
mot~$u$ dans le monoïde. On rappelle qu'un tel automate est la
donnée d'une fonction~$T: S \times A \ra S$ et d'un état initial~$q$;
l'état final de l'automate après lecture du mot~$u$ est, par définition,
l'état~$T(q, u)$ défini inductivement, en notant $\e$ le mot vide, par
$T(q, \e) = q$ et $T(q, ux) = T(T(q, u), x)$. On part de
formules de calcul de pgcd.

\Lem\«\Pgcd»
{\sl Soit $M$ un monoïde de Garside. Pour
$s , t $ simples, on a
$$s t  \gL \D = {}^*(t^* \cL s^{**}), 
\hbox{\qquad et \qquad}
s t  \gR \D = ({}^{**}t  \cR {}^*s )^*.\Eq\«200»$$

}\Dem La forme normale à gauche de~$s $ est $(s )$, donc,
par le lemme~«\Normale»(ii), celle de~$s  \D$ est $(\D,
s ^{**})$, et celle de~$s  t $, qui est $s  \D (t ^*)\ii$, est 
$({}^*(t^* \cL s^{**}), s^{**} \cL t^*)$. Or, par définition, le
premier terme de cette suite est~$s t  \gL \D$.

Un argument symétrique montre que la forme normale à droite
de~$\D t $ est $({}^{**}t , \D)$, et que celle de~$s t $ est
$({}^*s  \cR {}^{**}t , ({}^{**}t  \cR {}^*s )^*)$. Par
définition, le dernier terme de cette dernière suite est~$s t 
\gR \D$.\FinDem

\Prop
{\sl Soit $M$ un monoïde de Garside, $A$ l'ensemble de
ses atomes, $S$ l'ensemble de ses simples, et $\D$ l'élément
simple maximal. Définissons
$T: S \times A \ra S$ par $T(s , a  ) = ({}^{**}a   \cR
{}^*s)^*$. Alors, pour tout mot~$u$ sur~$A$, on a $\cl u \gR \D
= T(1, u)$, c'est-à-dire que le résultat de la lecture de~$u$
par l'automate~$(T, 1)$ est le dernier terme de la forme
normale à droite de~$\cl u$.}

\Dem On a trivialement $\cl\e \gR \D = 1 = T(1, \e)$,
donc, pour montrer le résultat inductivement, il suffit de
démontrer l'égalité
$$x a   \gR \D = T(x  \gR \D, a  )$$
pour $x$ dans~$M$ et $a  $ dans~$A$. Par le
lemme~«\Prelim»(ii) (en fait, sa contrepartie pour le pgcd à
droite), on a $(x a  ) \gR \D = ((x  \gR \D)a  ) \gR
\D$. Comme $x  \gR \D$ et $a  $ sont simples, le résultat
découle alors immédiatement de la formule~«200».\FinDem

\Rem
Comme dans~[\\\Eps], nous pouvons considérer un
transducteur, défini comme un automate muni d'une
fonction de sortie~$O$ envoyant $S \times A$ dans l'ensemble des
mots sur~$A$. On définit alors la sortie produite par la lecture du
mot~$u$ comme le mot~$O(q, u)$ inductivement défini par $O(q,
\e) = \e$ et $O(q, ua  ) = O(q, u) O(T(q, u), a  )$. La
démonstration précédente montre que, si nous posons $O(s, a  ) =
{}^*s \cR {}^{**}a    $, alors, pour tout  mot~$u$ sur~$A$, on a $\cl u
= \cl{O(1, u)} \om T(1, u)$: c'est dire que le mot~$O(1, u)$ représente
le reste du mot~$u$ lorsque le dernier terme de la forme
normale à droite a été retiré. En faisant lire par le
transducteur le mot~$O(1, u)$, on obtiendra de même
l'avant-dernier terme de la forme normale à droite de~$\cl u$,
et ainsi de suite. Cette itération fournit un algorithme de
complexité quadratique déterminant la forme normale à droite.

\Ex
La figure~3.1 représente le transducteur calculant la forme normale
à droite pour le monoïde $\Mon{a  , b  \; a  b a   =
b^2}$. Les états (domaines cerclés) correspondent aux éléments
simples. Les flèches pleines représentent la lecture de~$a  $, les
tiretées la lecture de~$b $; la présence d'une étiquette~$u$ sur
une flèche~$x$ entre les états~$s $ et~$t $ signifie qu'on a
$T(s , x) = t $ et $O(s , x) = u$, c'est-à-dire que, partant de
l'état~$s $ et lisant la lettre~$x$, on passe dans l'état~$t $ en
produisant le mot~$u$.

\midinsert
$$\OpenGraphicBox width 55mm height 48mm depth0mm;
\leavevmode\special{illustration 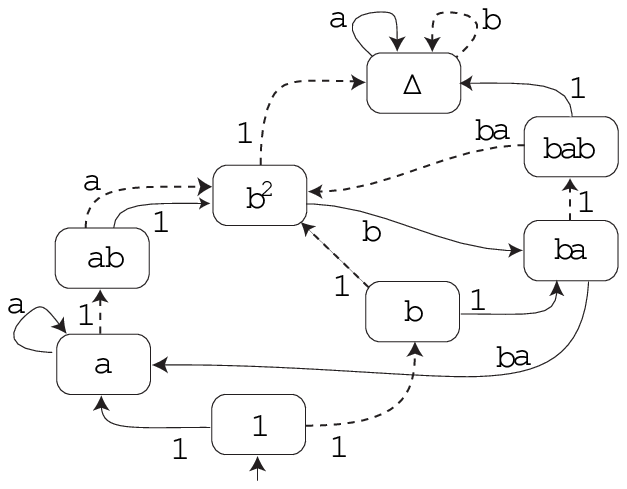}
\CloseGraphicBox$$ 
\centerline{{\bf Figure 3.2.} Transducteur pour $\Mon{a  , b  \;
a  b a   = b^2}$}
\endinsert

\bigskip Pour la forme normale à gauche, nous pouvons trivialement
obtenir des résultats analogues en considérant des
pseudo-automates lisant les mots de droite à gauche, mais la véritable
question est l'existence d'automates standard lisant de
gauche à droite. L'approche précédente échoue, car il est en général
faux que la valeur de~$(x a  ) \gL \D$ ne dépende que des valeurs
de~$x  \gL \D$ et de~$a  $: dans le groupe de
tresses~$B_3$, on a $\s_1\s_2 \gL \D = \s_1
\s_2^2 \gL \D = \s_1 \s_2$, mais $\s_1 \s_2 \s_1 \gL
\D = \D \not= \s_1 \s_2^2 \s_1 \gL \D = \s_1 \s_2$. Par
contre, la construction est encore possible en prenant comme
ensemble d'états l'ensemble~$P^{P}$, où $P$ est l'ensemble des
éléments primitifs à droite de~$M$, c'est-à-dire la clôture des atomes
par l'opération~$\cR$. 

\Prop
{\sl Soit $M$ un monoïde de Garside, $A$
l'ensemble de ses atomes, et $P$ l'ensemble de ses
éléments primitifs à droite. Définissons $T: P^{P} \times A
\ra P^{P}$ par 
$$T(f, a  ) = r_a   \comp f \Eq\«170»$$
où $r_a  $ est définie pour $a   \in A$ et $p \in P$ par
$r_a  (p) = a   \cR p$. Alors, pour tout mot~$u$ sur~$A$,
on a 
$$\cl u \gL \D = \bigvee \{p \in P \; T(\id, u)(p) = 1\}$$
--- et donc l'automate $(T, \id)$ détermine le premier
facteur de la forme normale à gauche de~$\cl u$.}

\Dem Pour~$x$ dans~$M$, notons $f_x$ l'application
de~$P$ dans~$P$ définie par $f_x (p) = x  \cR p$.
On a $f_1 = \id$ par définition, et, d'après les formules des
compléments, il vient, pour~$x  \in M$, $a   \in A$ et $p \in
P$,
$$f_{x a  }(p) = (x a  ) \cR p = a   \cR (x  \cR p) =
a   \cR f_x (p) = r_a  (f_x (p)),$$ 
soit $f_{x a  } = r_a   \comp f_x$, et, inductivement,
$f_{\overline u} = T(\id, u)$ pour tout mot~$u$
sur~$A$. Or, pour $p$ primitif, $p$ est un diviseur à
gauche de~$x$ si et seulement si on a $x  \cR p = 1$,
soit $f_x (p) = 1$. Le seul point à montrer est donc que
$x  \gL \D$ est égal au ppcm à droite~$x '$ des
éléments primitifs divisant~$x$ à gauche. Or, tout
élément primitif est simple, donc $x '$ est un diviseur à
gauche de~$x  \gL \D$. Inversement, $x  \gL \D$ est
simple, donc, d'après les résultats de la section~1, il est un
ppcm à droite d'élements primitifs, lesquels sont
nécessairement des diviseurs à gauche de~$x$.\FinDem

\Rem
L'approche précédente s'applique également à la
forme normale à droite. Introduisons l'ensemble~$\Til P$
des éléments {\it primitifs à gauche} comme la
clôture des atomes par l'opération~$\cL$. Pour~$x  \in M$,
définissons $\Til f_x  : \Til P \ra \Til P$ par $\Til f_x (p) = p \cL
x$. Alors la fonction~$\Til f_x$ est calculée inductivement par la
règle $\Til f_{x a  }(p) = \Til f_x (p \cL a  )$. On obtient donc
$\Til f_{\overline u} =
\Til T(\id, u)$, où
$\Til T = \Til P^{\Til P} \times A \ra \Til
P^{\Til P}$ est déterminée par
$$\Til T(f, a  ) = f \comp \Til r_a  , \Eq\«171»$$
avec $r_a   : p \mapsto p \cL a  $, formule qui est la
contre-partie exacte de~«170». On obtient alors, pour tout
mot~$u$ sur~$A$,
$$\cl u \gR \D = \Til{\bigvee} \{p \in \Til P \; \Til
T(\id, u)(p) = 1\}:$$
l'automate~$(\Til T, \id)$ détermine le dernier facteur de la
forme normale à droite de~$\cl u$.
Ces résultats montrent à nouveau le rôle central des
opérations~$\cR$ et~$\cL$: les automates précédents sont
essentiellement la table de ces opérations sur les clôtures des atomes. 

\Ex 
La figure~3.2 montre l'automate calculant le pgcd à gauche avec~$\D$
dans le cas du monoïde de tresses~$B_3^+$, de présentation
$\Mon{a  , b  \;
a  b a   = b a  b }$. Il y a deux atomes, à savoir $a  $ et~$b $,
et cinq éléments primitifs à droite (et à gauche), à savoir $1,
a  , b , a  b , b a  $. On obtient un automate à 20~états,
représenté sur la figure~3.2 (les flèches en plein
représentent~$a  $, les flèches en tireté représentent~$b $).

\midinsert
$$\OpenGraphicBox width 135mm height 45mm
depth0mm;
\leavevmode\special{illustration 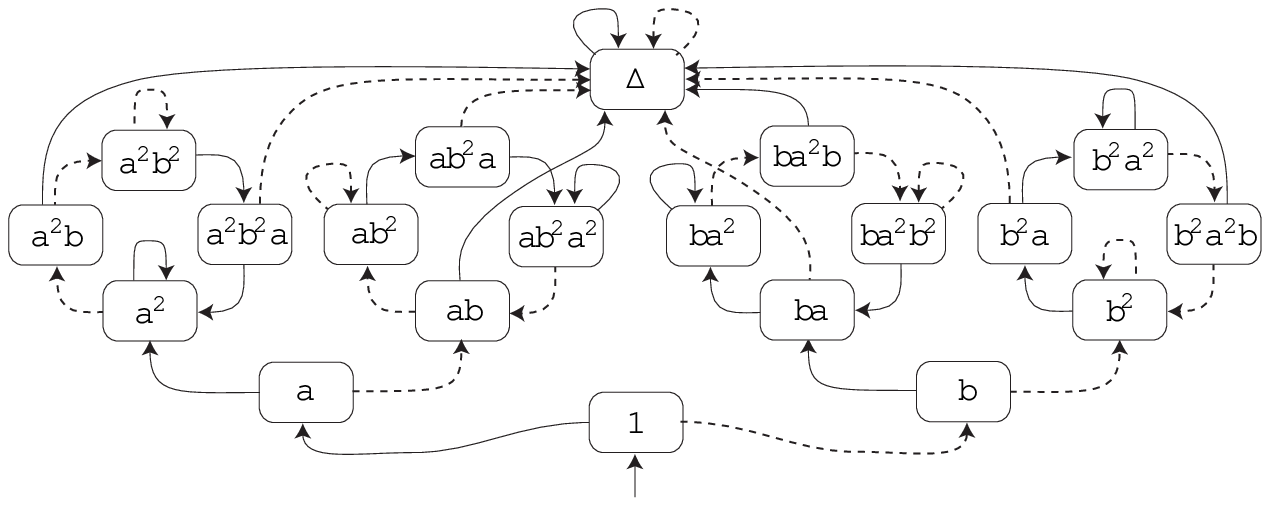}
\CloseGraphicBox$$ 
\centerline{{\bf Figure 3.2.} Automate pour le pgcd à gauche
dans~$B_3$}
\endinsert

\Sec Présentations

\bigskip On montre ici que tout groupe de Garside admet une
présentation d'un certain type, dit complémenté, et que toutes
les opérations du monoïde se calculent à partir d'une telle
présentation au moyen d'une opération combinatoire appelée
redressement de mots. Les résultats du début de la section
figurent, avec des notations différentes, dans~[\\\Dfx], 
auquel nous renvoyons pour d'autres détails et les
démonstrations manquantes. 

Pour tout ensemble~$\Sigma$, nous noterons~$\Mo(\Sigma)$
l'ensemble des mots sur~$\Sigma$, c'est-à-dire le monoïde libre
de base~$\Sigma$; le mot vide est noté~$\e$. Si $\Sigma$ est
une partie génératrice d'un monoïde~$M$, alors, par
définition, tout élément~$x$ de~$M$ admet des
décompositions comme produit d'éléments de~$\Sigma$:
pour $w$ dans~$\Mo(\Sigma)$, on notera comme
précédemment $\cl w$ l'image~$x$ de~$w$ dans~$M$; on
dira alors que $w$ est une expression de~$x$, et que $x$ est
l'évaluation de~$w$ dans~$M$.

\Def
Soit $M$ un monoïde vérifiant les
conditions~$(C_0)$, $(C_1)$ et $(C_2)$, et  $\Sigma$ est une
partie génératrice de~$M$. On appelle {\it sélecteur de ppcm}
sur~$\Sigma$ toute application (partielle)~$f$ de~$\Sigma \times
\Sigma$ dans~$\Mo(\Sigma)$ tel que, pour $a  , b  \in \Sigma$,
$f(a  , b )$ est une expression de~$a   \cR b $ lorsque ce
dernier existe.

\bigskip Par définition, on a $a   \cR a   = 1$
pour tout~$a  $ dans~$\Sigma$, donc, si $f$ est un
sélecteur de ppcm sur~$\Sigma$, et si $M$ satisfait à la
condition~$(C_0)$, c'est-à-dire n'admet pas
d'inversible autre que~$1$, alors, nécessairement,
$f(a  , a  )$ est le mot vide~$\e$. Par ailleurs, par
définition, $a   \cR b $ existe si et seulement si $b 
\cR a  $ existe, et, par conséquent, le domaine de tout
sélecteur de ppcm est un sous-ensemble symétrique
de~$\Sigma \times \Sigma$.

\Def 
Soit $\Sigma$ un ensemble non vide. On appelle
{\it fonction de redressement} sur~$\Sigma$ une application
partielle de~$\Sigma \times \Sigma$ dans~$\Mo(\Sigma)$ vérifiant
$f(a  , a  ) = \e$ pour tout~$a  $, et telle que le domaine
de~$f$ est symétrique. On note alors $\Rf$ la famille des
relations $a   f(a  , b ) = b  f(b , a  )$ pour $(a  , b )
\in \Dom(f)$, et $\eqfp$ la congruence sur~$\Mo(\Sigma)$
engendrée par~$\Rf$. On dit qu'une présentation de monoïde
est {\it complémentée} si elle est associée à une
(nécessairement unique) fonction de redressement.

\Prop [\\\Dfx]
{\sl Soit $M$ un monoïde gaussien --- ou,
plus généralement, un monoïde vérifiant les
conditions~$(C_0^+)$, $(C_1)$ et $(C_2)$. Soit
$\Sigma$ une partie génératrice quelconque de~$M$, et $f$ un
sélecteur de ppcm sur~$\Sigma$. Alors $f$ est une
fonction de redressement sur~$\Sigma$ et $\Mon{\Sigma
\; \Rf}$ est une présentation complémentée de~$M$.}

\Dem Notons $\sim$ la congruence sur~$\Mo(\Sigma)$ telle que
$M$ est $\Mo(\Sigma) / \! \sim$. Par définition d'un sélecteur
de ppcm, on a $a   f(a  , b ) \sim b  f(b , a  )$
pour $a  $, $b  \in \Dom(f)$, donc $u \eqfp v$
entraîne $u \sim v$ puisque les paires $\{a   f(a  ,
b ), b  f(b , a  )\}$ engendrent~$\eqfp$.

Réciproquement, nous montrons que $u \sim v$ entraîne $u
\eqfp v$ par récurrence sur~$\norm(\cl u)$. Pour $\norm(\cl
u) = 0$, on a  $\cl u = \cl v = 1$, donc, par~$(C_0)$, $u = v =
\e$, et $u \eqfp v$. Supposons
$u, v \not= \e$. On écrit $u = a   u_1$, $v = b  v_1$, avec
$a   , b  \in \Sigma$. L'hypothèse $\cl u = \cl v$ signifie que
$\cl u$ est un multiple commun à droite de~$a  $
et~$b $. Par~$(C_2)$, il existe un mot~$w$ vérifiant $\cl
u = \cl v = (a   \jR b ) \cl w$, d'où $u = a   u_1  \sim
a   f(a  , b ) w \sim b  f(b , a  ) w \sim b  v_1 = v$.
Par~$(C_1)$, on déduit $u_1 \sim f(a  , b ) w$ et $v_1 \sim
f(b , a  ) w$. Or, par~$(C_0^+)$, on a
$\norm(\cl{f(a  , b )w}) \le \norm(\cl u) - \norm(a  ) <
\norm(\cl u)$, et, de même, $\norm(\cl{f(b , a  )w}) <
\norm(\cl u)$. Appliquant l'hypothèse de
récurrence, on déduit $u_1 \eqfp f(a  , b ) w$, et $v_1 \eqfp
f(b , a  ) w$, d'où $u = a   u_1 \eqfp a   f(a  ,
b ) w \eqfp b  f(b , a  ) w \eqfp b  v_1 = v$, et donc $u
\eqfp v$.\FinDem

\Question
{\sl Le résultat précédent reste-t-il valable lorsque la
condition~$(C_0^+)$ (atomicité) est affaiblie en~$(C_0)$ (pas
d'inversible autre que~$1$)?}

\bigskip Si $M$ est un monoïde gaussien, $\Sigma$ une partie
génératrice de~$M$, et $f$ un sélecteur de ppcm sur~$\Sigma$,
alors toutes les opérations de~$M$ se calculent à partir de~$f$
de façon effective. Pour le démontrer, nous introduisons une
opération combinatoire sur les mots appelée {\it
redressement} (\<<reversing\>> en anglais). 

Par définition, si $f$ est un sélecteur de ppcm sur~$\Sigma$, le
mot $f(a  , b )$ est une expression de l'élément~$a   \cR
b $ de~$M$ pour tous~$a  $, $b $ dans~$\Sigma$. L'idée est
d'étendre l'application~$f$ en une opération binaire
partielle notée~$\cRf$ sur~$\Mo(\Sigma)$ de sorte que, pour tous
mots~$u$, $v$ sur~$\Sigma$, le mot~$u \cRf v$ soit une
expression du complément~$\cl u \cR \cl v$, quand ce dernier
existe.

Dans toute la suite, lorsque $\Sigma$ est un alphabet, nous
introduirons une copie disjointe notée~$a  \ii$ pour
chaque lettre~$a  $ de~$\Sigma$, et nous noterons~$\Sigma\ii$
l'ensemble des lettres~$a  \ii$. Pour tout mot~$w$
sur~$\Sigma \cup \Sigma\ii$, on note $w\ii$ le mot obtenu en
échangeant chaque lettre~$a  $ avec la lettre~$a  \ii$
correspondante, et en renversant l'ordre des lettres. Ainsi, si
$G$ est un groupe engendré par~$\Sigma$ et si le mot~$w$
représente l'élément~$x$ de~$G$, alors $w\ii$
représente~$x \ii$.

\Def [\\\Dff] [\\\Dfx]
Soit $f$ une fonction de redressement
sur~$\Sigma$. Pour~$w$, $w'$ mots sur~$\Sigma \cup
\Sigma\ii$, on dit que $w$ est {\it $f$-redressable} en~$w'$, noté
$w \vRf w'$, si on peut passer de~$w$ à~$w'$ en un nombre 
fini d'étapes consistant à remplacer un sous-mot de la
forme~$a  \ii b $ avec $a  , b  \in \Sigma$ par le mot $f(a  ,
b ) f(b , a  )\ii$ correspondant. Pour $u$, $v$ mots
sur~$\Sigma$, on définit $u \cRf v$ comme l'unique mot~$u'$
sur~$\Sigma$ tel qu'il existe un mot~$v'$ sur~$\Sigma$ vérifiant
$u\ii v \vRf v' {u'}\ii$, s'il existe, et $u \jRf v$ comme $u \om
(u \cRf v)$.

\bigskip On associe à chaque redressement de mot un graphe 
planaire (voisin d'un diagramme de Dehn) composé de flèches
verticales et horizontales étiquetées par des éléments
de~$\Sigma$, et tel que le redressement du mot~$a  \ii b $ en
$f(a  , b ) f(b , a  )\ii$ se traduit par la fermeture du motif
$$\OpenGraphicBox width 88mm height 19mm depth0mm;
\leavevmode\special{illustration 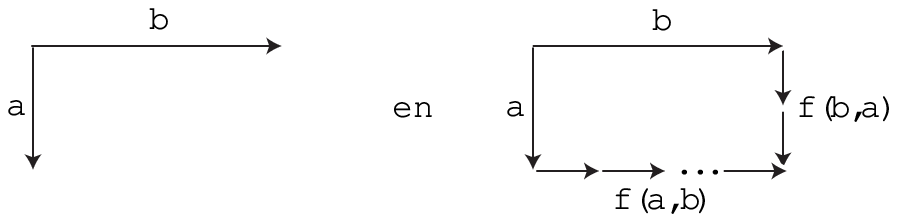}
\CloseGraphicBox$$ 
Nous renvoyons à~[\\\Dgd, chapitre~II], et
considérons ici simplement un exemple.

\Ex
Soit $\Sigma = \{a  , b \}$, et soit $f$ l'application de~$\Sigma^2$
dans~$\Mo(\Sigma)$ définie par $f(a  , b ) = b a  $, $f(b ,
a  ) =
b $, $f(a  , a  ) =  f(b , b ) = \e$. Le
monoïde~$\Mon{\Sigma \; \Rf}$ est ici
$\Mon{a  , b  \; a  b a   = b^2}$. La figure ci-dessous
illustre le redressement du mot $a  \ii
b  a  ^2$ en le mot $b  a   b^2 a  
b \ii a  \ii b \ii$, qui donne donc $a   \cRf b a  ^2 =
b a  b^2a  $ et $b a  ^2 \cRf a   = b a  b $.
$$\OpenGraphicBox width 50mm height 33mm
depth0mm;
\leavevmode\special{illustration 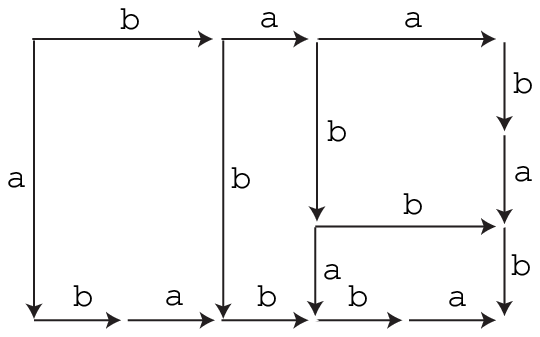}
\CloseGraphicBox$$ 

\bigskip L'opération~$\cRf$
est en général une opération partielle sur~$\Mo(\Sigma)$, non
nécessairement définie pour tous les mots: comme $\cRf$
coïncide avec~$f$ sur $\Sigma \times \Sigma$, c'est le cas si
$f$ n'est pas partout définie sur~$\Sigma
\times \Sigma$, mais, même si $f$ est partout définie, il se peut
que le redressement d'un mot ne se termine pas en un
nombre fini d'étapes. Comme dans la section~1, nous
noterons $u \cRf v = \bot$ quand $u \cRf v$ n'existe pas, et
nous étendons $\cRf$ en une opération partout définie
sur~$\Mo(\Sigma) \cup \{\bot\}$ en posant $u
\cRf \bot = \bot
\cRf u = \bot \cRf \bot = \bot$ pour tout~$u$, et de même
pour le produit des mots. De même, nous
étendons~$\eqfp$ en une congruence sur~$\Mo(\Sigma)
\cup \{\bot\}$ en déclarant $\bot \eqfp \bot$ vrai,
et $u \eqfp \bot$ faux pour tout~$u$ dans~$\Mo(\Sigma)$. 

\Lem\«\Equivalence» [\\\Dff] 
{\sl Soit $f$ une fonction de redressement
sur~$\Sigma$. Alors, pour tous mots~$u, v$ sur~$\Sigma$, on a $u(u
\cRf v) \eqfp v (v \cRf u)$.}

\bigskip D'après nos conventions, le résultat précédent
signifie que soit $u \cRf v$ et $v \cRf u$ sont définis et 
on a l'équivalence annoncée, soit ni l'un ni l'autre n'est
défini. La démonstration est une induction facile sur le
nombre d'étapes élémentaires de redressement.

\Prop [\\\Dfx] \«\Expression»
{\sl Soit $M$ un monoïde gaussien --- ou,
plus généralement, un monoïde vérifiant les
conditions~$(C_0^+)$, $(C_1)$ et $(C_2)$ ---, $\Sigma$ une
partie génératrice de~$M$, et $f$ un sélecteur de
ppcm sur~$\Sigma$. Alors, pour tous mots~$u$, $v$
dans~$\Mo(\Sigma)$, le mot~$u \cRf v$ existe
dans~$\Mo(\Sigma)$ si et seulement si l'élément $\cl u \cR
\cl v$ existe dans~$M$, et, dans ce cas, $u \cRf v$ est une
expression de $\cl u \cR \cl v$, et $u \jRf v$ est une
expression de~$\cl u \jR \cl v$.}

\Dem Le lemme~«\Equivalence» montre que, si $u \cRf v$
est défini, alors la classe de~$u \jRf v$ est un multiple à
droite commun des classes de~$u$ et~$v$, donc,
par~$(C_2)$, le ppcm à droite de ces classes existe.

Inversement, montrons par induction sur~$n$ que, si $u$ et
$v$ sont deux mots tels que $\cl u \jR \cl v$ existe et qu'on
ait $\norm(\cl u \jR \cl v) = n$, alors $u \cRf v$ existe
et représente~$\cl u \cR \cl v$. Pour $n = 0$, on a  $u = v =
\e$, et les résultats ont triviaux. Supposons $n \ge 1$. Pour $u
= \e$ ou $v = \e$, les résultats sont à nouveau triviaux.
Supposons
$u =
a   u_1$ et $v = b  v_1$ avec $a  , b  \in \Sigma$, et soit
$z  = \cl u \jR \cl v$. Par hypothèse, $z $ est
multiple à droite commun à~$a  $ et~$b $, donc $f(a  ,
b )$ et $f(b , a  )$ sont définis. Ensuite $z$ est multiple
à droite commum de $\cl{\mathstrut a   u_1}$ et $\cl{a  
f(a  , b )}$, donc $\cl{\mathstrut u_1}$ et $\cl{f(a  ,
b )}$ ont un multiple à droite commun, donc un ppcm à
droite. Or, par construction, on a $\norm(\cl{\mathstrut
u_1} \jR \cl{f(a  , b )}) \le \norm(z) - \norm(a  ) = n-1$,
donc, par hypothèse de récurrence, $u_1 \cRf f(a  , b )$
existe et il représente $\cl{\mathstrut u_1} \cR \cl{f(a  ,
b )}$. De même, $v_1 \cRf f(b , a  )$
existe et  représente $\cl{\mathstrut v_1} \cR
\cl{f(b , a  )}$; finalement, en posant $u_2 = f(a  , b )
\cRf u_1$ et $v_2 = f(xy, a  ) \cRf v_1$, on a de
même que $u_2 \cRf v_2$ existe et  représente
$\cl{u_2} \cR \cl{v_2}$. Donc  $u \cRf v$ est
défini, et il représente
$$(\cl{\mathstrut v_1} \cR
\cl{f(b , a  )}) \om ((\cl{f(b , a  )} \cR \cl{\mathstrut
v_1}) \cR (\cl{f(a  , b )} \cR \cl{\mathstrut u_1})),$$
lequel est, d'après les formules de complément, $\cl u \cR \cl
v$.\FinDem

\Def
Soit $f$ une fonction de redressement
sur~$\Sigma$. Pour $u, v$ mots sur~$\Sigma$, on dit que $u \eqfpp
v$ est vérifié si on a $u \cRf v = v \cRf u = \e$, c'est-à-dire
$u\ii v \vRf \e$. On étend~$\eqfpp$ de sorte que $\bot
\eqfpp \bot$ soit vrai, et $u \eqfpp \bot$ et $\bot \eqfpp u$
faux pour tout mot~$u$.

\bigskip Par le lemme~«\Equivalence», $u \eqfpp v$
entraîne~$u \eqfp v$. Dans le cas d'un monoïde
gaussien, la proposition~«\Expression» affirme que cette
implication est une équivalence. Nous allons en déduire une
solution pour le problème de mot du monoïde par
redressement, et, pour celui du groupe associé, par double
redressement. On notera
$\eqf$ la congruence sur~$\WX$ engendrée
par~$\eqfp$ et les paires $\{a   a  \ii, \e\}$ et $\{a  \ii a  ,
\e\}$ pour~$a  $ dans~$\Sigma$, de sorte que le groupe
$\Gr{\Sigma \; \Rf}$ est $\WX / \!\! \eqf$. Remarquons
que, par construction, $w \vRf w'$ entraîne
$w \eqf w'$ pour tous~$w, w'$ mots sur~$\Sigma \cup \Sigma\ii$.

\Prop\«\WordProblem»
{\sl Soit $M$ un monoïde gaussien --- ou,
plus généralement, un monoïde vérifiant les
conditions~$(C_0^+)$, $(C_1)$ et $(C_2)$ ---, $\Sigma$ une
partie génératrice de~$M$, et $f$ un sélecteur de
ppcm sur~$\Sigma$. 

(i) Deux mots~$u, v$ sur~$\Sigma$ représentent le
même élément de~$M$ si et seulement on a $u \eqfpp v$,
c'est-à-dire $u\ii v \vRf \e$.

(ii) Supposons que $M$ vérifie de surcroît la
condition~$(C_2^+)$, et soit $G$ le groupe de fractions
de~$M$. Alors un mot~$w$ sur~$\Sigma \cup \Sigma\ii$
représente~$1$ dans~$G$ si et seulement si il existe deux
mots~$u, v$ sur~$\Sigma$ vérifiant $w \vRf uv\ii$ et $u\ii v \vRf
\e$.}

\Dem (i) Par le lemme~«\Equivalence», la condition est
suffisante. Inversement, $\cl u = \cl v$ entraîne
$\cl u \cR \cl v = \cl v \cR \cl u = 1$, d'où $u \cRf v = v
\cRf u = \e$ par la proposion~«\Expression».

(ii) La condition est suffisante, car $w \vRf w'$ entraîne
$w \eqf w'$. Inversement, supposons $w \vRf u v\ii$. Alors $w \eqf \e$
entraîne $u \eqf v$, donc $u \eqfp v$ puisque $M$ se plonge
dans~$G$, d'où $u\ii v \vRf \e$ par~(i). \FinDem

\bigskip Si $M$ est un monoïde gaussien, les
opérations~$\cR$ et~$\jR$ de~$M$ se déterminent par
redressement à partir de tout sélecteur de ppcm. Le pgcd, ainsi
que les opérations symétriques de ppcm à gauche et de pgcd à
droite, peuvent également être calculées à l'aide de
redressements, à condition d'introduire, à côté d'un
sélecteur à droite associé au ppcm à droite comme
ci-dessus, une notion symétrique de sélecteur à gauche
et de redressement à gauche: supposant $\Sigma$ 
partie génératrice de~$M$, on appellera {\it
sélecteur de ppcm à gauche} sur~$\Sigma$ une
application~$\Til f$ de $\Sigma \times \Sigma$ dans~$\Mo(\Sigma)$
telle que, pour tous~$a  , b $ dans~$\Sigma$, $\Til
f(a  , b )$ représente l'élément	$a   \cL b $, s'il
existe. Si $\Til f$ est une fonction de redressement
sur~$\Sigma$, et si $w, w'$ sont deux mots sur~$\Sigma
\cup \Sigma\ii$, on dira que~$w$ est $\Til f$-redressable à
gauche en~$w'$, noté $w \vLf w'$, si on peut passer de~$w$
à~$w'$ en un nombre fini d'étapes consistant à remplacer un
factuer du type~$a   b \ii$ avec $a  , b  \in \Sigma$ par le
facteur $\Til f(b , a  )\ii \Til f(a  , b )$ correspondant. Les
résultats sont alors symétriques: introduisant $u \cLf v$
comme l'unique mot~$u'$ tel qu'il existe $v'$ vérifiant $u v\ii
\vLf v'{}\ii u'$, et $u \jLf v$ comme $(u \cLf v) \om v$, on
obtient que $u \cLf v$ représente l'élément~$\cl u \cL \cl v$,
et $u \jLf v$ l'élément~$\cl u \jL \cl v$.  On déduit alors du
lemme~«\Gcd»:

\Prop
{\sl Soit $M$ un monoïde gaussien, $\Sigma$
une partie génératrice de~$M$, et $f$ et $\Til f$
respectivement un sélecteur de ppcm à droite et à
gauche sur~$\Sigma$. Alors, pour tous mots~$u$, $v$
sur~$\Sigma$, le mot~$(u \jRf v) \cLf ((u
\cRf v) \jLf (v \cRf u))$ représente $\cl u \gL \cl v$.}

\bigskip On obtient donc le pgcd à gauche de~$\cl u$
et~$\cl v$ par un triple redressement: on redresse à droite
$u\ii v$ en $v' u'{}\ii$, qu'on redresse à gauche  en~$u''{}\ii v''$, et
un représentant du pgcd à gauche de~$\cl u$ et~$\cl v$ est obtenu en
redressant à gauche le mot~$u u''{}\ii$. Les opérations du
groupe de fractions se calculent de façon semblable. En
particulier, dénominateurs et numérateurs se déterminent par
double redressement.

\Prop
{\sl Soit $G$ le groupe de fractions d'un monoïde
gaussien~$M$, $\Sigma$ une partie génératrice de~$M$, et
$f$ et $\Til f$ respectivement un sélecteur de ppcm à droite et
à gauche sur~$\Sigma$. Soit $z $ un élément quelconque
de~$G$, et $w$ un mot sur~$\Sigma \cup \Sigma\ii$
représentant~$z $. Soient
$u, v, u', v'$ les mots sur~$\Sigma$ vérifiant $w \vRf v u\ii \vLf
u'{}\ii v'$. Alors $u'$ représente~$D (z )$, et $v'$
représente~$N (z )$.}

\Dem Par construction, on a $z  = \cl{u'}{}\ii \cl{v'}$, donc,
d'après le lemme~«\Fraction», il suffit de montrer $\cl{u'} \gL
\cl{v'} = 1$. Or, par construction, on a $u' = u \cLf v$ et $v' = v
\cLf u$: si $\cl{u'}$ et $\cl{v'}$ avaient un diviseur à gauche
non trivial commun, $u' v$ et $v' u$ ne pourraient  représenter
un ppcm à gauche de~$\cl u$ et~$\cl v$.
\FinDem

\bigskip Comme $z  = 1$ équivaut à $D (z ) = N (z ) =
1$, nous en déduisons une solution alternative du problème de
mot. 

\Cor
{\sl Sous les mêmes hypothèses, un mot~$w$ sur~$\Sigma \cup \Sigma\ii$
représente~$1$ dans~$G$ si, et seulement si, il existe deux mots~$u,
v$ sur~$\Sigma$ vérifiant $w \vRf u v\ii \vLf \e$.}

\bigskip Enfin, le lemme~«\Normale» montre comment calculer
par redressement la forme normale (à gauche) de tout élément d'un
groupe de Garside: supposant que $(\cl{v_q}\,\ii, \pp,
\cl{v_1}\,\ii, \cl{u_1}, \pp, \cl{u_p})$ soit la forme normale
de~$z $, on obtient la forme normale de~$z  \D$ en
déterminant les éléments~${\cl u_i}^{**}$, donc par double
redressement face à un mot représentant~$\D$. Pour~$w$
représentant un élément simple, on obtient le numérateur de
la forme normale de~$z  \cl w{}\ii$ en redressant à gauche
le mot $u_1 \ppp u_p w\ii$; le dénominateur s'obtient alors
par une succession de calculs de pgcd, donc également par
redressement. Noter que le lemme~«\Pgcd» permet de
déterminer successivement chaque facteur du dénominateur à
l'aide d'un unique redressement.

\Sec Reconnaître les monoïdes de Garside

\bigskip Nous avons vu dans la section précédente que tout
monoïde de Garside admet une présentation
complémentée $\Mon{\Sigma; \Rf}$, où $f$ est une
fonction de redressement sur~$\Sigma$, et qu'alors toutes les
opérations du monoïde et celles de son groupe de fractions se
déterminent par $f$-redressement. Nous abordons ici la
question réciproque de reconnaître quand une présentation
complémentée définit un monoïde de Garside.

Le fait d'admettre une présentation complémentée est
une hypothèse faible, et des conditions sur la fonction de
redressement considérée doivent être ajoutées pour que le
monoïde associé ait de bonnes propriétés. Nous allons établir
une liste de conditions nécessaires vérifiées par tout sélecteur
de ppcm dans un monoïde de Garside, puis montrer
que ces conditions sont suffisantes. Nous partirons
(une fois encore) des propriétés de l'opération~$\cR$.

\Lem\«\PreCoherence»
{\sl Soit $M$ un monoïde gaussien --- ou,
plus généralement, un monoïde vérifiant les
conditions~$(C_0)$, $(C_1)$ et $(C_2)$. Alors, pour
tous~$x , y , z $ dans~$M$, on a
$$(x  \cR y ) \cR (x  \cR z ) = (y  \cR x ) \cR (y  \cR
z ). \Eq\«200»$$

}\Dem D'après le lemme~«\Rules», le membre de gauche
dans~«200» est égal à $(x  \jR y ) \cR z $, et celui de droite
à $(y  \jR x ) \cR z $: comme $\jR$ est une opération
commutative, les deux sont égaux. Par ailleurs, si l'une des
expressions n'est pas définie, il en est de même de
l'autre.\FinDem

\Lem\«\Coherent»
{\sl Soit $M$ un monoïde gaussien --- ou,
plus généralement, un monoïde vérifiant les
conditions~$(C_0^+)$, $(C_1)$ et $(C_2)$ ---, 
$\Sigma$ une partie géné\-ratrice quelconque de~$M$, et $f$ un
sélecteur de ppcm sur~$\Sigma$. Alors, pour tous
mots~$u$, $v$, $w$ dans~$\Mo(\Sigma)$, on a
$$(u \cRf v) \cRf (u \cRf w) \eqfpp (v \cRf u) \cRf (v \cRf
w). \Eq\«\COH»$$

}\Dem Compte tenu de la proposition~«\Expression», la
formule~«\COH» est la traduction directe de la
formule~«200».\FinDem

\Def
Soit $f$ une fonction de redressement sur~$\Sigma$. 

(i) Pour $X \ince \Mo(\Sigma) \cup \{\bot\}$, on dit que $f$ 
a la {\it propriété du cube} (\resp la {\it propriété du cube
faible}) sur~$X$ si la relation~«\COH»\ (\resp la relation
$$(u \cRf v) \cRf (u \cRf w) \eqfp (v \cRf u) \cRf (v \cRf
w)~) \Eq\«\wCOH»$$
est vérifiée pour tous~$u$, $v$, $w$ dans~$X$.

(ii) On dit que $f$ vérifie la condition~$(\CCI)$ s'il existe une
partie finie~$X$ de~$\Mo(\Sigma)$ qui est close par~$\cRf$,
c'est-à-dire telle que, pour tous~$u, v$ dans~$X$, $u \cRf v$ existe
et appartient à~$X$, et si, de plus, la clôture~$X \vv$ de~$X$
par~$\jRf$ est finie et il existe~$\O$ dans~$X\vv$ tel que,
pour tout~$u$ dans~$X\vv$, on  ait $\O \cRf  u = \e$ et il
existe~$v$ dans~$X \vv$ vérifiant $v \cRf \O \eqfpp u$.

\bigskip Affirmer que $f$ a la propriété du cube sur
un triplet~$\{u, v, w\}$ signifie, lorsque tous les mots mis en
jeu sont définis, que le cube représenté sur la figure~5.1 se
ferme au sens où, partant des trois arêtes étiquetées~$u, v, w$
et utilisant le redressement de mots pour construire les six
faces, le redressement des trois petites faces triangulaires
restantes se termine avec des mots vides.

\midinsert
$$\OpenGraphicBox width 50mm height 48mm depth0mm;
\leavevmode\special{illustration 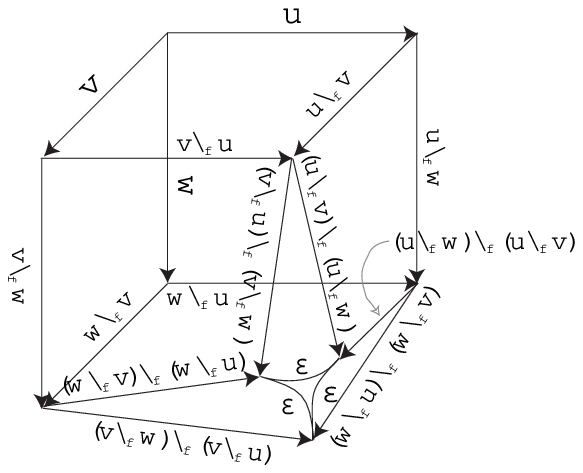}
\CloseGraphicBox$$
\centerline{{\bf Figure 5.1.} Propriété du cube}
\endinsert

\Prop\«\Necess»
{\sl Soit $M$ un monoïde de Garside, $\Sigma$ une partie
génératrice finie de~$M$, et $f$ un sélecteur de ppcm
sur~$\Sigma$. Alors la fonction~$f$ a la propriété du cube
sur~$\Mo(\Sigma)$ et elle vérifie la condition~$(\CCI)$.}

\Dem Le lemme~«\Coherent» exprime que $f$ a la propriété
 du cube sur~$\Mo(\Sigma)$. Notons~$X$ l'ensemble des
expressions des éléments primitifs à droite de~$M$. Comme
$\Sigma$ est fini, et que la longueur des expressions
dans~$\Mo(\Sigma)$ de chaque élément~$x$ de~$M$ est
bornée supérieurement par~$\norm(x )$, l'ensemble~$X$ est
fini. Supposons $u , v \in X$. Par la
proposition~«\Expression», le mot $u \cRf v$ est une
expression de~$\cl u \cR \cl v$: comme $\cl u$ et $\cl v$ sont
primitifs à droite, il en est de même de~$\cl u \cR
\cl v$, et $u \cRf v$ appartient à~$X$. Donc $X$
est clos par~$\cRf$.

Par le même raisonnement, l'ensemble~$Y$ des expressions
des éléments simples de~$M$ est un ensemble fini de mots
qui est clos par les opérations~$\cRf$ et~$\jRf$, et qui
inclut~$X$. Par conséquent, la clôture~$X \vv$ de~$X$
par~$\jRf$ est finie, et elle contient au moins une expression
pour chaque élément simple de~$M$. Soit $\O$ une
expression quelconque de l'élément~$\D$ dans~$X \vv$. Pour
tout mot~$u$ dans~$X \vv$, le mot~$\O
\cRf u$ représente $\D \cR \cl u$, qui est~$1$, et on a donc
nécessairement $\O \cRf u = \e$. Inversement, si $\O$ vérifie $\O
\cRf u = \e$ pour tout~$u$ dans~$X \vv$, $\O$ représente
nécesairement~$\D$, puisqu'il représente un élément
qui est un multiple de tous les éléments simples.

D'après la proposition~«\Char», l'application~$s  \mapsto s  \cR
\D$ est surjective sur les simples de~$M$. Donc, pour tout mot~$u$
dans~$X \vv$, il doit exister un élément simple~$s $ 
vérifiant $\cl u = s  \cR \D$. Cet élément~$s $ a au moins
une expression~$v$ dans~$X \vv$, et la condition $\cl u =
\cl v \cR \D$ entraîne $u \eqfpp v \cRf \O$. Par conséquent, la
condition~$(\CCI)$ est vérifiée.\FinDem

\bigskip Nous allons maintenant montrer que les conditions
nécessaires de la proposition~«\Necess» sont aussi
suffisantes, et, pour cela, utiliser la caractérisation de la
proposition~«\Char» en termes des conditions~$(C_i)$. La
condition~$(C_0)$ est gratuite.

\Lem\«\Conical»
{\sl Soit $f$ une fonction de redressement sur~$\Sigma$. Alors le
monoïde $\Mon{\Sigma \; \Rf}$ satisfait à la condition~$(C_0)$,
c'est-à-dire n'a d'autre inversible que~$1$.}

\Dem Par construction, $u \eqfp \e$ n'est possible que
pour $u = \e$.\FinDem

\bigskip Le point crucial est le résultat suivant, qui montre
l'importance de la propriété du cube.

\Prop\«\CohComp»
{\sl Soit $f$ une fonction de redressement sur~$\Sigma$. Alors
les conditions suivantes sont équivalentes:

(i) La fonction~$f$ a la propriété du cube sur~$\Mo(\Sigma)$;

(ii) Les relations~$\eqfp$ et~$\eqfpp$ coïncident;

(iii) La relation~$\eqfp$ est compatible avec
l'opération~$\cRf$, au sens où la conjonction de $u'
\eqfp u$ et $v' \eqfp v$ entraîne $u' \cRf v' \eqfp u
\cRf v$.}

\Dem Supposons~(i), et montrons~(ii). Par le
lemme~«\Equivalence», la relation~$\eqfpp$ est toujours
incluse dans la relation~$\eqfp$. Pour prouver l'inclusion
réciproque, puisque, par définition, $\eqfp$ est la relation
d'équivalence engendrée par les paires~$(ua   f(a  , b )v,
u b  f(b , a  )v)$ avec $a  , b  \in \Sigma$ et
$u, v \in \Mo(\Sigma)$, il suffit de montrer qu'on a $u a   f(a  ,
b )v
\eqfpp  u b  f(b , a  )v$, soit
$$v\ii f(b , a  )\ii b \ii u\ii u a   f(a  , b ) v \vRf \e,$$
ce qui est trivial, et que $\eqfpp$ est une relation
d'équivalence. Réflexivité et symétrie sont claires, et
seule la transitivité fait problème. Supposons $u
\eqfpp v \eqfpp w$. Par hypothèse, les mots $v \cRf u$, $u
\cRf v$, $v \cRf w$ et $w \cRf v$ existent et sont vides.
Donc $(v \cRf u) \cRf (v \cRf w)$ existe, et il est vide.
Comme $f$ a la propriété du cube sur~$\{v, u, w\}$, ceci
entraîne $(u
\cRf v) \cRf (u \cRf w) \eqfp \e$, donc $(u \cRf v) \cRf (u
\cRf w) = \e$. L'existence du mot $(u \cRf v) \cRf (u \cRf w)$
entraîne en particulier celle de $u \cRf w$. De plus, $u \cRf
v = \e$ entraîne $(u \cRf v) \cRf (u \cRf w) = u \cRf w$.
Nous avons donc $u \cRf w = \e$. Un argument symétrique
donne $w \cRf u = \e$, donc $u \eqfpp w$. Donc $\eqfpp$ est
transitive, et (i) entraîne~(ii).

Supposons maintenant (ii), et montrons~(iii). Puisque la
relation~$\eqfp$ est symétrique et transitive, il suffit de montrer
que, si $u \cRf v$ existe et qu'on a $v'
\eqfp v$, alors $u \cRf v'$ existe aussi et on a $u \cRf v' \eqfp u \cRf
v$ et $v' \cRf u \eqfp v \cRf u$. Or, sous ces hypothèse,  le
lemme~«\Equivalence» entraîne $u(u \cRf v) \eqfp v(v \cRf
u)$, donc 
$u(u \cRf v) \eqfp v'(v \cRf u)$,  soit, par~(ii), $u(u
\cRf v) \eqfpp v'(v \cRf u)$, et donc
$(v' (v \cRf u)) \cRf (u (u \cRf v)) = \e$. Ceci entraîne
en particulier $(v' (v \cRf u)) \cRf u = \e$, soit, par~«103», 
$$(v \cRf u) \cRf (v' \cRf u) = \e. \Eq\«147»$$ 
Nous déduisons que $v' \cRf u$, et, par conséquent, $u \cRf
v'$, existent. Comme $v$ et $v'$ jouent maintenant des
rôles symétriques, le calcul précédent donne $(v' \cRf u)
\cRf (v \cRf u) = \e$, ce qui, conjugué à~«147»,
démontre $v' \cRf u \eqfpp v \cRf u$, donc  $v'
\cRf u \eqfp v \cRf u$. Un calcul semblable donne
$(u \cRf v) \cRf (u \cRf v') = \e$ et
$(u \cRf v') \cRf (u \cRf v) = \e$, d'où  $u \cRf v
\eqfpp u' \cRf v$, et, finalement $u \cRf v \eqfp u'
\cRf v$. Donc (ii) implique~(iii).

On sait que $u' \eqfpp u$ entraîne $u' \eqfp u$. Inversement,
supposons $u' \eqfp u$. Par construction, on a 
$u\ii u \vRf \e$, c'est-à-dire $u \cRf u = \e$. Si (iii) est vraie, on
déduit que $u' \eqfp u$ entraîne l'existence de $u \cRf u'$ et
$u' \cRf u$, et qu'on a $u' \cRf u \eqfp u
\cRf u = \e$, et $u \cRf u' \eqfp u \cRf u = \e$, soit $u'
\eqfpp u$. Ceci montre que (iii) entraîne~(ii).

Finalement, supposons~(ii) et~(iii), et soient $u, v, w$ des
mots quelconques sur~$\Sigma$. Si $u \cRf v$ n'est pas défini, $v
\cRf u$ ne l'est pas non plus, et $(u \cRf v) \cRf
(u \cRf w)$ et $(v \cRf u) \cRf (v \cRf w)$ sont tous
deux~$\bot$, donc «\COH» est vraie. Supposons que $u
\cRf v$ et $v \cRf u$ soient définis. Par~«103», on a
$(u \cRf v) \cRf (u \cRf w) = (u(u \cRf v)) \cRf w$ et $(v
\cRf u) \cRf (v \cRf w) = (v(v \cRf u)) \cRf w$. Supposons
$(u(u \cRf v)) \cRf w$ défini. Par le
lemme~«\Equivalence», on a $u(u \cRf v) \eqfp v(v
\cRf u)$, donc $u \,(u \cRf v) \eqfp v\, (v \cRf u)$. Alors
(iii) entraîne que $(v(v \cRf u)) \cRf w$ est également
défini, et qu'on a $(v(v \cRf u)) \cRf w \eqfp (u(u \cRf v))
\cRf w$. Par~(ii), on déduit $(v(v \cRf u)) \cRf
w \eqfpp (u(u \cRf v)) \cRf w$, donc en particulier $((u(u
\cRf v)) \cRf w) \cRf ((v(v \cRf u)) \cRf w) = \e$. Un calcul
symétrique donne $((v(v \cRf u)) \cRf w) \cRf ((u(u \cRf
vu)) \cRf w) = \e$, d'où $(u(u \cRf v)) \cRf w \eqfpp (v(v
\cRf u)) \cRf w$, qui est la condition~«\COH» pour~$u, v,
w$. Donc (ii) et~(iii) entraînent~(i). \FinDem

\bigskip Il est alors facile d'établir les conditions~$(C_1)$
et~$(C_2)$.

\Lem\«\Cancel»
{\sl Soit $f$ une fonction de redressement sur~$\Sigma$
ayant la propriété du cube sur~$\Mo(\Sigma)$. Alors le
monoïde~$\Mon{\Sigma \; \Rf}$ satisfait à la
condition~$(C_1)$, c'est-à-dire admet la simplification à
gauche.}

\Dem Soit $M$ le monoïde $\Mon{\Sigma \; \Rf}$. Supposons 
$\cl u \, \cl v = \cl{\mathstrut u} \,
\cl{\mathstrut v'}$ dans~$M$, soit
$uv \eqfp uv'$. Par la proposition~«\CohComp», on a $uv
\eqfpp uv'$, c'est-à-dire $(uv)\ii(uv') \vRf \e$. Par
construction, on a $(uv)\ii(uv') \vRf v\ii v'$, donc, par
unicité du résultat du redressement, $v\ii v' \vRf \e$, soit $v
\eqfpp v'$, qui entraîne $v \eqfp v'$, et $\cl v =
\cl{v'}$. \FinDem

\Lem\«\Lcm»
{\sl Soit $f$ une fonction de redressement sur~$\Sigma$
ayant la propriété du cube sur~$\Mo(\Sigma)$. Alors, quels
que soient les mots~$u$, $v$ sur~$\Sigma$, il y a équivalence
entre 

(i) les mots $u \cRf v$ et $v \cRf u$ sont définis;

(ii) les éléments~$\cl u$ et $\cl v$ admettent un multiple à
droite commun dans le monoïde~$\Mon{\Sigma; \Rf}$.

Dans ce cas, $\cl u$ et $\cl v$ admettent un ppcm à
droite, à savoir $\cl{u \jRf v}$, et $u \cRf v$
représente~$\cl u \cR \cl v$.}

\Dem Supposons~(i). Par le lemme~«\Equivalence», on a
$u(u \cRf v) \eqfp v(v \cRf u)$, et la classe commune de ces
deux mots est un multiple à droite commun pour les
classes de~$u$ et~$v$. Donc (ii) est vérifiée.  

Inversement, supposons (ii). Il existe des mots
$u'$ et $v'$ vérifiant $u v' \eqfp v u'$. Puisque
$f$ a la propriété du cube, ceci implique $u v' \eqfpp v u'$
par la proposition~«\CohComp», soit ${v'}\ii u\ii v u' \vRf \e$.
Ainsi, le redressement du mot ${v'}\ii u\ii v u'$ converge,
et, {\it  a fortiori}, celui du sous-mot~$u\ii v$ converge
aussi, ce qui est dire que les mots $u \cRf v$ et $v \cRf u$
existent, et on a~(i). De plus, l'hypothèse ${v'}\ii u\ii v u'
\vRf \e$ implique l'existence de mots~$u''$, $v''$ satisfaisant
$${v'}\ii (u \cRf v) \vRf {v''}\ii, \quad
(v \cRf u)\ii u' \vRf u'', \quad
{v''}\ii u'' \vRf \e.$$
Ainsi, $\cl u'$ est un multiple à droite de la classe
de~$v \cRf u$, et, donc, $\cl{vu'}$ est un multiple à droite
de la classe de~$v \jRf u$: cette dernière est donc
ppcm à droite de~$\cl u$ et~$\cl v$.\FinDem

\Lem\«\CII»
{\sl Soit $f$ une fonction de redressement sur~$\Sigma$
ayant la propriété du cube sur~$\Mo(\Sigma)$. Alors le
monoïde~$\Mon{\Sigma \; \Rf}$ satisfait à la
condition~$(C_2)$: toute paire d'éléments ayant un multiple à
droite commun admet un ppcm à droite.}

\Dem Supposons que $\cl u$ et $\cl v$ admettent un
multiple à droite commun dans~$M$. Par le lemme~«\Lcm»,
les mots $u \cRf v$ et $v \cRf u$ sont définis, et $u \jRf v$
et $v \jRf u$ représentent un ppcm à droite de~$\cl u$
et~$\cl v$.\FinDem

\bigskip Nous passons à la condition~$(C_3)$.

\Lem\«\Convergent»
{\sl Soit $f$ une fonction de redressement sur~$\Sigma$, telle
qu'il existe une partie~$X$ de~$\Mo(\Sigma)$ close
par~$\cRf$. Alors le $f$-redressement converge toujours, et
donc l'opération~$\cRf$ est partout définie sur~$\Mo(\Sigma)$. 
 De plus, si le redressement de tout mot~$u\ii v$ avec $u, v
\in X$ requiert au plus $k$~étapes et se termine avec un
mot de longueur au plus~$\ell$, alors le redressement
d'un mot quelconque sur~$\Sigma \cup \Sigma\ii$ requiert au plus ${1
\over 4}k \om \lg(w)^2$ étapes, et il se termine avec un
mot de longueur au plus~$\ell \om \lg(w)$.}

\Dem Soit $w$ un mot sur $\Sigma \cup \Sigma\ii$. Ecrivons $w =
u_1^{e_1} \ppp u_r^{e_r}$ avec $u_i \in X$,
$e_i = \pm 1$ pour chaque~$i$, et $r \le \lg(w)$. Soit $p$
le nombre d'exposants~$e_i$ positifs. Un argument inductif
illustré sur la figure~5.2 montre qu'il existe des mots
 $v_1$, \pp, $v_r$ dans~$X$ tels que $w$
soit $f$-redressable en $v_1 \ppp v_p v_{p+1}\ii \ppp
v_r\ii$, et que le redressement se décompose en $p(r -p)$
redressements de mots de la forme $w_1\ii w_2$ avec
$w_1$, $w_2$ dans~$X$. Les bornes en résultent, puisqu'on
a toujours $p(r -p) \le r^2/4$.\FinDem

\midinsert
$$\OpenGraphicBox width 50mm height37mm depth0mm;
\leavevmode\special{illustration Convergence4.eps};
\CloseGraphicBox$$
\centerline{{\bf Figure 5.2.} Convergence du
redressement}
\endinsert

\Lem\«\CIII»
{\sl Soit $f$ une fonction de redressement sur~$\Sigma$,
ayant la propriété du cube sur~$\Mo(\Sigma)$, et telle que
$f$ satisfasse à la condition~$(\CCI)$. Alors le
monoïde~$\Mon{\Sigma \; \Rf}$ satisfait à la
condition~$(C_3)$.}

\Dem Soit $X$ la clôture de~$\Sigma$ par l'opération~$\cRf$.
D'après le lemme~«\Convergent», l'opération~$\cRf$ est
partout définie sur~$\Mon{\Sigma \; \Rf}$. Il résulte du
lemme~«\Lcm» que l'ensemble des éléments du monoïde
représentés par un mot de~$X$ est clos par l'opération~$\cR$,
et qu'il engendre le monoïde puisqu'il inclut~$\Sigma$. \FinDem

\bigskip Finalement, nous déduisons la réciproque de la
proposition~«\Necess»:

\Prop\«\Suffis»
{\sl Soit $f$ une fonction de redressement  sur~$\Sigma$,
ayant la propriété du cube sur~$\Mo(\Sigma)$, et satisfaisant
à la condition~$(\CCI)$. Alors le monoïde~$\Mon{\Sigma
\; \Rf}$ est un monoïde de Garside, et, donc, le
groupe~$\Gr{\Sigma \; \Rf}$ est un groupe de Garside.}

\Dem Par la proposition~«\Char», il suffit d'établir que le
monoïde~$\Mon{\Sigma \; \Rf}$ satisfait aux
conditions~$(C_0)$, $(C_1)$, $(C_2)$, $(C_3)$ et $(\Til C_1)$.
Les quatre premières conditions résultent des
lemmes~«\Conical», «\Cancel», «\CII» et «\CIII». Enfin,
pour~$(\Til C_1)$, d'après le lemme~«\Bijection», il suffit
d'établir la surjectivité de l'application $x  \mapsto x^*$
sur les éléments simples; compte tenu de la correspondance 
entre opération~$\cR$ sur~$M$ et opération~$\cRf$
sur~$\Mo(\Sigma)$, et entre égalité dans~$M$ et
$\eqfpp$-equivalence dans~$\Mo(\Sigma)$ due à la
satisfaction par~$f$ de la prioriété du cube, la
condition~$(\CCI)$ donne le résultat.
\FinDem

\Sec Critères pour la propriété du cube

\bigskip La propriété du cube pour la fonction~$f$ est
primordiale pour l'étude du monoïde $\Mon{\Sigma \; \Rf}$:
si elle n'est pas satisfaite, on ne sait essentiellement rien dire,
alors que, si elle l'est, on peut contrôler la divisibilité et
l'égalité à l'aide de redressements de mots et, en particulier,
établir assez simplement l'éventuel caractère petit gaussien du
monoïde. Il est donc crucial de savoir reconnaître si une
fonction de redressement donnée satisfait ou non la
propriété du cube. 

Si $f$ est une fonction de redressement sur un ensemble (fini)~$\Sigma$,
et $u, v, w$ des mots sur~$\Sigma$, établir la propriété du
cube en~$u, v, w$ se fait de façon effective à l'aide de
redressements. Plus précisément, il existe un processus
consistant en une suite finie de redressements tel que, si la
condition est vérifiée, alors le processus se termine en un
nombre fini d'étapes et il donne une preuve de la condition
cherchée. Cette situation, plus faible que la décidabilité
puisqu'il se peut, si la condition n'est pas vérifiée, qu'on
n'obtienne pas de réponse en un temps fini, est typique d'une
condition récursivement énumérable --- ou semi-décidable, ou
$\S^0_1$ [\\\BaM].

Considérons maintenant la propriété du cube sur
l'ensemble~$\Mo(\Sigma)$ entier: comme il existe une infinité de
triplets de mots dans~$\Mo(\Sigma)$, énumérer ceux-ci
systématiquement et vérifier la condition pour chacun d'eux
ne permet pas d'établir, même de façon théorique, la
propriété en un nombre fini d'étapes. On ne peut
obtenir un critère effectif que si on sait par avance que la
propriété du cube sur un certain sous-ensemble fini
de~$\Mo(\Sigma)$ entraîne la propriété du cube partout. 

Nous rappellons d'abord deux résultats partiels  dans cette
direction obtenus antérieurement. Le premier est implicite
dans la preuve de Garside que le monoïde des
tresses~$B_n^+$ satisfait aux conditions~$(C_1)$ et~$(C_2)$.

\Prop\«\CritI» [\\\Dff] [\\\Dgc]
{\sl Soit $f$ une fonction de redressement 
sur~$\Sigma$ telle que le monoïde~$\Mon{\Sigma \; \Rf}$
satisfasse~$(C_0^+)$. Alors $f$ a la propriété du cube
sur~$\Mo(\Sigma)$ si et seulement si elle a la propriété du
cube sur~$\Sigma$, si et seulement si elle a la propriété du
cube faible sur~$\Sigma$.}

\bigskip Le critère précédent s'applique aux présentations
standards des groupes d'Artin, la condition~$(C_0^+)$ étant
vérifiée puisque les relations préservent la longueur: les
résultats de [\\\BrS] montrent alors la propriété du cube sur
les générateurs. Il en est de même pour les présentations
de groupe considérés dans~[\\\BuM]: elles sont associées à
des fonctions de redressement~$f$ telles que $f(a, b)$ n'est 
défini que pour $a$ dans~$A$ et $b$ dans~$B$, $\{A, B\}$
étant une partition de l'ensemble des générateurs: la
propriété du cube sur les lettres est triviale car tout triplet
contient nécessairement deux lettres de~$A$, ou deux lettres
de~$B$, pour lesquelles aucun redressement n'est possible; les
groupes ainsi présentés ne sont pas gaussiens, puisque, dans le
monoïde associé, deux éléments de~$A$, ou de~$B$, n'ont pas
de multiple commun.

Mis à part de tels cas, il est en général difficile d'établir 
la condition~$(C_0^+)$, et aucune méthode uniforme
n'est connue --- voir [\\\Pic] pour un exemple mettant en jeu
l'algorithme de Knuth-Bendix.

Un autre critère pour la propriété du cube est établi
dans~[\\\Dgc], à partir d'une définition alternative de celle-ci:

\Lem\«\Alternative»
{\sl Soit $f$ une fonction de redressement sur~$\Sigma$, et
$X$ une partie de~$\Mo(\Sigma)$. Alors $f$ a la propriété du
cube sur~$X$ si et seulement si on a
$$w \cRf (u(u \cRf v)) \eqfp w \cRf (v(v \cRf u)), \quad
\hbox{et}
\quad (u(u \cRf v)) \cRf w \eqfp (v(v \cRf u)) \cRf
w\Eq\«55»$$
pour tous~$u$, $v$, $w$ dans~$X$.} 

\Dem La relation~«103» donne $(u \cRf v) \cRf (u \cRf
w) = (u(u \cRf v)) \cRf w$ et $(v \cRf u) \cRf (v \cRf w) =
(v(v \cRf u)) \cRf w$, donc «\COH» est équivalente à la
seconde des relations~«55». Pour la première, en
utilisant~«103» et le lemme~«\Equivalence», et en supposant 
«\COH» satisfaite pour $(u, w, v)$ et pour $(w, v, u)$, on
obtient
$$\Eqalign{
w \cRf (u(u \cRf v)) 
&= (w \cRf u) \cdot ((u \cRf w) \cRf (u \cRf v)) \cr
&\eqfp (w \cRf u) \cdot ((w \cRf u) \cRf (w \cRf v)) \cr
&\eqfp (w \cRf v) \cdot ((w \cRf v) \cRf (w \cRf u)) \cr
&\eqfp (w \cRf v) \cdot ((v \cRf w) \cRf (v \cRf u)) 
=w \cRf (v(v \cRf u)). &\Block\cr}$$

\Def
Soit $f$ une fonction de redressement sur~$\Sigma$, et 
$X$ une partie de~$\Mo(\Sigma)$. Pour $u \eqfp v$,
définissons $\df(u, v)$ comme le nombre minimal de
$\Rf$-relations nécessaires pour transformer~$u$ en~$v$. On
dit alors que $f$ est {\it $k$-cohérente} sur~$X$ si, pour
tous~$u, v, w$ dans~$X$, on a
$$\df(w \cRf u(u \cRf v), w \cRf v(v \cRf u)) + \df((u(u \cRf
v)) \cRf w, (v(v \cRf u)) \cRf w) \le k \Eq\«57»$$

\bigskip Par définition, si $X$ est fini, et si $f$ a la propriété
du cube (faible) sur~$X$, alors elle est $k$-cohérente pour un
certain~$k$. Un critère suffisant mais non nécessaire est:

\Prop\«\CritII» [\\\Dgc]
{\sl Si $f$ est une fonction de redressement sur~$\Sigma$ 
qui est $1$-cohérente sur~$\Sigma$, alors $f$ a la propriété du
cube sur~$\Mo(\Sigma)$.}

\bigskip Le critère de la proposition~«\CritII» est
automatiquement vérifié dans le cas de deux générateurs,
mais il l'est rarement dans le cas général. 

Nous allons maintenant établir un nouveau critère, qui
s'applique dans tous les cas et ne nécessite aucune vérification
préalable.

\Prop\«\CritIII»
{\sl Soit $f$ une fonction de redressement sur~$\Sigma$, et
$X$ une partie de~$\Mo(\Sigma) \cup \{\bot\}$ qui
inclut~$\Sigma$ et est close par~$\cRf$. Alors $f$ a la propriété du
cube sur~$\Mo(\Sigma)$ si, et seulement si, elle a la propriété du
cube sur~$X$.}

\bigskip Il est évident que la condition est nécessaire, et tout
le travail consiste à montrer qu'elle est suffisante. La
démonstration est décomposée en plusieurs étapes.
Jusqu'à la fin de la section, nous supposons que $f$ est une
fonction de redressement fixé sur~$\Sigma$, que $X$ est
une partie de~$\Mo(\Sigma) \cup \{\bot\}$ qui inclut~$\Sigma$ et
est close par~$\cRf$, et que $f$ a la propriété du cube
sur~$X$. D'après la proposition~«\CohComp», il suffit que
nous montrions la compatibilité de l'opération~$\cRf$ et de la
relation~$\eqfp$.

\Lem\«\CompatI» 
{\sl Pour $u, v, v' \in X$, la relation $v' \eqfpp v$
entraîne $v' \cRf u \eqfpp v \cRf u$ et $u \cRf v' \eqfpp u
\cRf v$.}

\Dem Suppsons que $v \cRf u$ existe. Par hypothèse, nous
avons $v' \cRf v = v \cRf v' = \e$. Comme $f$ a la propriété du
cube sur~$\{v, v', u\}$, nous déduisons
$$v' \cRf u = \e \cRf (v' \cRf u)  = (v' \cRf v) \cRf (v' \cRf u)
\eqfpp (v \cRf v') \cRf (v \cRf u) = \e \cRf (v \cRf u) =  v
\cRf u.$$ 
Ceci entraîne en particulier que $v' \cRf u$ et donc $u \cRf
v'$ existent. Utilisant la propriété du cube en~$\{u, v, v'\}$, nous
trouvons
$$(u \cRf v) \cRf (u \cRf v') \eqfpp (v \cRf u) \cRf (v \cRf v') = (v
\cRf u) \cRf \e = \e,$$
qui entraîne $(u \cRf v) \cRf (u \cRf v')  = \e$. On obtient
$(u \cRf v') \cRf (u \cRf v)  = \e$ symétriquement, d'où
on déduit $u \cRf v' \eqfpp u \cRf v$. \FinDem

\Lem\«\Transit» 
{\sl (i) La relation~$\eqfpp$ est transitive sur~$X$.

(ii) Pour $u, v, u', v'$ dans~$X$, la conjonction de $u'
\eqfpp u$ et $v' \eqfpp v$ entraîne $v' \cRf u' \eqfpp v \cRf
u$.}

\Dem (i) Supposons $u, v, w \in X$ et $u \eqfpp v \eqfpp
w$.  Par hypothèse, on a $u \cRf v = v \cRf u = \e$.
Appliquant le lemme~«\CompatI»\ à $v \eqfpp w$, on
obtient $u
\cRf w \eqfpp \e$ et $w \cRf u \eqfpp \e$, donc $u \cRf w = w
\cRf u = \e$, soit $u \eqfpp w$. 

(ii) Supposons que $v \cRf u$ existe. En appliquant
le lemma~«\CompatI», on obtient
$$v' \cRf u' \eqfpp v' \cRf u \eqfpp v \cRf u,$$
qui entraîne $v' \cRf u' \eqfpp v \cRf u$ par~(i), puisque
$v' \cRf u'$, $v' \cRf u$ et $v \cRf u$ appartiennent à~$X$, ce
dernier étant supposé clos par~$\cRf$.\FinDem

\bigskip Nous introduisons maintenant des raffinements de la
relation~$\eqfpp$.

\Def Pour $u$, $u'$ dans~$\Mo(\Sigma) \cup \{\bot\}$, on dira
que $u \spe(0) u'$ est vérifié si on a  $u = u' = \e$, et, pour 
$p \ge 1$, que $u \spe(p) u'$ est vérifié s'il existe deux
décompositions $u = u_1 \ppp u_p$, $u'= u'_1 \ppp u'_p$
telles qu'on ait $u_1$, \pp, $u'_p \in X$ et $u'_j \eqfpp u_j$
pour tout~$j$.

\bigskip Ainsi, $u' \spe(1) u$ est la conjonction de~$u' \eqfpp
u$ et de $u, u' \in X$. 

\Lem\«\FineI» 
{\sl (i) Toute relation $u \spe(p) u'$ entraîne $u \eqfpp u'$.

(ii) La conjonction des relations $u' \spe(p) u$ et $v' \spe(q) v$
entraîne $u' \cRf v' \spe(q) u \cRf v$ et $v'
\cRf u' \spe(p) v \cRf u$.}

\Dem (i) Le résultat est clair par récurrence sur~$p$
(noter que, même pour $u, u' \in X$, il n'y a en général
aucune raison pour que, réciproquement, $u \eqfpp u'$
entraîne $u \spe(p) u'$).

(ii) Observons d'abord que le résultat est vrai pour $p =
0$ et pour $q = 0$. En effet, pour $p = 0$, on a $u' = u =
\e$, donc $u' \cRf v' = \e = u \cRf v$, et, pour $q = 0$, on a
$v' = v = \e$, donc $v' \cRf u' = u'$ et $v \cRf u = u$.

Ensuite, nous utilisons une récurrence sur~$p + q$.
Supposons $u' \spe(p) u$ et $v' \spe(q) v$. D'après ce qui
précède, nous pouvons supposer $p \ge 1$ et
$q \ge 1$. Alors, on peut écrire $u = u_1 u_2$, $u' = u'_1
u'_2$ avec $u'_1 \eqfpp u_1$ et $u'_2 \spe(p-1) u_2$,
et, de même, $v = v_1 v_2$, $v' = v'_1 v'_2$ avec $v'_1
\eqfpp v_1$ et $v'_2 \spe(q-1) v_2$. Posons $u_{0, j} = u_j$
et $v_{i, 0} = v_i$, puis $u_{i, j} = v_{i-1, j} \cRf u_{i, j-1}$ et
$v_{i, j} = u_{i, j-1} \cRf v_{i-1, j}$ (Figure~6.1). L'hypothèse
que $X$ est close par~$\cRf$ entraîne inductivement que
tous les mots~$u_{i, 1}$, $u'_{i, 1}$, $v_i$, et $v'_i$ sont
dans~$X  \cup \{\bot\}$. Le lemme~«\Transit»(ii) donne
$u'_{1, 1} \eqfpp u_{1, 1}$ et $v'_{1, 1} \eqfpp v_{1, 1}$.
Ensuite, l'hypothèse de récurrence donne $u'_{1, 2}
\spe(p-1) u_{1, 2}$ et $v'_{1, 2} \eqfpp v_{1, 2}$. Elle donne de
même $u'_{2, 1} \eqfpp  u_{2, 1}$ et $v'_{2, 1}
\spe(q-1) v_{2, 1}$, et, finalement, $u'_{2, 2} \spe(p-1)
u_{2, 2}$ et $v'_{2, 2} \spe(q-1) v_{2, 2}$. Réunissant les
relations, on déduit $u' \cRf v' = u'_{2, 1} u'_{2, 2}
\spe(p) u_{2, 1} u_{2, 2} = u \cRf v$, et $v' \cRf u' = v'_{1, 2}
v'_{2, 2} \spe(q) v_{1, 2} v_{2, 2} = v \cRf u$.
\FinDem

\midinsert
$$\OpenGraphicBox width 45mm height24mm depth0mm;
\leavevmode\special{illustration 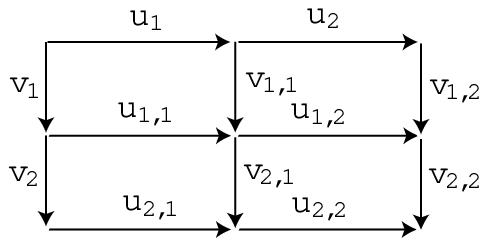}
\CloseGraphicBox$$ 
\centerline{{\bf Figure 6.1.} Démonstration du
lemme~«\FineI»}
\endinsert

\Def\«\DefSPE»
Pour $v$, $v'$ dans~$\Mo(\Sigma) \cup \{\bot\}$, on dira
que $v \SPE v'$ est vérifié s'il existe deux décompositions
$v = v_0 v_1 v_2 v_3$, $v'= v'_0 v'_1 v'_2 v'_3$ et
deux entiers~$q, r$ vérifiant

(i) $v'_0 \spe(q) v_0$, et $v'_3 \spe(r) v_3$, 

(ii) $v_1, v_2, v'_1, v'_2 \in X$, et $v'_2 \eqfpp
v_1 \cRf v'_1$, et $v_2 \eqfpp v'_1 \cRf v_1$.

\Lem 
{\sl La relation $v \SPE v'$ entraîne $v \eqfpp v'$.}

\Dem Avec les notations de la définition, on a d'abord $v'_0
\eqfpp v_0$ et $v'_3 \eqfpp v_3$ par le lemme~«\FineI»(i),
soit $v_0\ii v'_0 \vRf \e$ et $v_3\ii v'_3 \vRf \e$. De même,
par hypothèse, on a $v_2 \eqfpp v'_1 \cRf v_1$, soit
$v_2\ii (v'_1 \cRf v_1) \vRf \e$, et, symétriquement, 
${v'_2}\ii(v'_1 \cRf v_1) \vRf \e$. On trouve alors
$$\eqalign{
v\ii v' 
&= v_3\ii v_2\ii v_1\ii v_0\ii v'_0 v'_1 v'_2 v'_3 \cr
& \vRf v_3\ii v_2\ii v_1\ii v'_1 v'_2 v'_3
\vRf v_3\ii v_2\ii (v_1 \cRf v'_1) (v'_1 \cRf v_1)\ii v'_2 v'_3
\vRf v_3\ii v'_3 \vRf \e, \cr}$$ 
lorsque tous les mots~$v_i$, $v'_i$ sont dans~$\Mo(\Sigma)$, et on
conclut que $v \eqfpp v'$ est vrai. D'un autre côté, si au
moins un des mots~$v_i$ est~$\bot$, le mot~$v'_i$
correspondant est également~$\bot$, et on a
$v = v' = \bot$, donc $v \eqfpp v'$.\FinDem

\Lem\«\FineII»
{\sl La conjonction de $u' \spe(p) u$ et $v' \SPE v$ entraîne
$u' \cRf v' \SPE u \cRf v$ et $v' \cRf u' \spe(p) v \cRf u$.}

\Dem Nous raisonnons par récurrence sur~$p$. Supposons
$u' \spe(p) u$ et $v' \SPE v$. Pour $p = 0$, la seule
possibilité est $u' = u = \e$, donc $u' \cRf v' = v' \SPE v = u
\cRf u$, et le résultat est vrai. Supposons $p = 1$. Fixons
des décompositions $v = v_0 v_1 v_2 v_3$, $v' = v'_0 v'_1
v'_2 v'_3$ comme dans la définition~«\DefSPE». On pose
$u_0 = u$ et $u'_0 = u'$, et on définit de proche en proche
$u_j = v_{j-1} \cRf u_{j-1}$ et $w_{j} = u_j \cRf v_j$, et, de
même, $u'_j = v'_{j-1} \cRf u'_{j-1}$ et $w'_{j} = u'_j \cRf
v'_j$. 
$$\OpenGraphicBox width 85mm height 15mm depth0mm;
\leavevmode\special{illustration 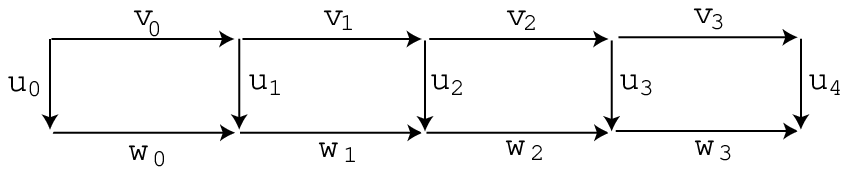}
\CloseGraphicBox$$ 
Par hypothèse, nous avons $u'
\spe(1) u$ et $v'_0 \spe(q) v_0$, donc le lemme~«\FineI»
donne
$$\Displaylines{
u'_1 \spe(1) u_1, \Eq\«69» \cr
w'_{0} \spe(q) w_{0}.\Eq\«70» \cr}$$
Par hypothèse, $v_1$ et~$v'_1$ sont dans~$X$, donc il en
est de même de $v_1 \cRf v'_1$ et $v'_1 \cRf v_1$.
Toujours par hypothèse, nous avons $v_2 \eqfpp v_1
\cRf v'_1$ et $v'_2 \eqfpp v'_1 \cRf v_1$, donc, par le
lemme~«\CompatI», nous obtenons
$$\Displaylines{
u_3 = v_2 \cRf u_2 \eqfpp (v_1 \cRf v'_1) \cRf u_2
= (v_1 \cRf v'_1) \cRf (v_1 \cRf u_1), \Eq\«71»\cr
w_{2} = u_2 \cRf v_2 \eqfpp u_2 \cRf (v_1 \cRf v'_1)
= (v_1 \cRf u_1) \cRf (v_1 \cRf v'_1), \Eq\«72»\cr
}$$
et, {\it mutatis mutandis},
$$\Displaylines{
u'_3 \eqfpp = (v'_1 \cRf v_1) \cRf (v'_1 \cRf u'_1), \Eq\«73»\cr
w'_{2} \eqfpp (v'_1 \cRf u'_1) \cRf (v'_1 \cRf v_1), \Eq\«74»\cr
}$$
Comme $f$ a la propriété du cube sur~$X$, le mot de droite
dans~«71» est $\eqfpp$-équivalent à $(v'_1 \cRf v_1) \cRf
(v'_1 \cRf u_1)$, et, comme tous les mots concernés sont
dans~$X$ et que $\eqfpp$ est transitive sur~$X$ d'après le
lemme~«\Transit»(ii), nous déduisons
$$u_3 \eqfpp (v'_1 \cRf v_1) \cRf (v'_1 \cRf u_1). \Eq\«75»$$
Nous avons vu que $u'_1 \eqfpp u_1$ est vrai, et $u_1$,
$u'_1$, $v'_1$ et $v'_1 \cRf v_1$ sont dans~$X$,
donc, en appliquant le lemme~«\CompatI» deux fois
en partant de~«73» et de~«75», nous obtenons
$u'_3 \eqfpp u_3$, et même $u'_3 \spe(1) u_3$
puisque $u_3$ et $u'_3$ sont dans~$X$ par construction.

En appliquant de même la propriété du cube et la transitivité
de~$\eqfpp$ sur~$X$, nous obtenons à partir de~«72» et
«74» les relations
$$\Displaylines{
w_{2} \eqfpp (u_1 \cRf v_1) \cRf (u_1 \cRf v'_1) 
= w_{1} \cRf (u_1 \cRf v'_1),\Eq\«76»\cr
w'_{2} \eqfpp (u'_1 \cRf v'_1) \cRf (u'_1 \cRf v_1) 
= w'_{1} \cRf (u'_1 \cRf v_1).\Eq\«77»\cr}$$ 
Comme on a $u'_1 \eqfpp u_1$, le lemme~«\CompatI»
donne d'abord $u_1 \cRf v'_1 \eqfpp u'_1 \cRf v'_1 = w'_{1}$, d'où, 
en utilisant~«76», le lemme~«\CompatI»\ à nouveau, et la transitivité
de~$\eqfpp$ sur~$X$, 
$$w_{2} \eqfpp w_{1} \cRf w'_{1},\Eq\«78»$$
et, par un argument symétrique en partant de~«77»,
$$w'_{2} \eqfpp w'_{1} \cRf w_{1}.\Eq\«79»$$
Maintenant, nous avons vu que $u'_3 \spe(1) u_3$ est
vérifiée, et, par hypothèse, nous avons $v'_3
\spe(r) v_3$. En appliquant le lemme~«\FineI», nous
déduisons
$$\Displaylines{
u'_4 \spe(1) u_4, \Eq\«80» \cr
w'_{3} \spe(r) w_{3}. \Eq\«81» \cr}$$
Par construction, on a $u \cRf v = w_{0} w_1 w_2 w_{3}$
et $v \cRf u = u_4$, et, de même, $u' \cRf v' = w'_{0} w'_1 w'_2
w'_{3}$ et $v' \cRf u' = u'_4$. Alors, la conjonction
de~«69», «78», «79» et «81» donne $u' \cRf v' \SPE u
\cRf v$, et «80» donne $v' \cRf u' \spe(1) v \cRf u$, ce qui
est le résultat escompté.

Supposons finalement $p \ge 2$. On écrit $u = u_1 u_2$, $u'
= u'_1 u'_2$ avec $u'_1 \spe(1) u_1$ et $u'_2 \spe(p-1)
u_2$. En appliquant l'hypothèse de récurrence à~$u_1$,
$u'_1$, $v$ et~$v'$, nous obtenons
$$\Displaylines{
u'_1 \cRf v' \SPE u_1 \cRf v \Eq\«80»\cr
v' \cRf u'_1 \spe(1) v \cRf u_1. \Eq\«81»\cr}$$
Ensuite, en utilisant~«80» et en appliquant
l'hypothèse de récurrence à~$u_2$, $u'_2$, $u_1 \cRf v$
et $u'_1 \cRf v'$, nous obtenons
$$\Displaylines{
u'_2 \cRf (u'_1 \cRf v') \SPE u_2 \cRf (u_1 \cRf v) \Eq\«82»\cr
(u'_1 \cRf v') \cRf u'_2 \spe(p) (u_1 \cRf v) \cRf u_2.
\Eq\«83»\cr}$$  
La relation~«82» est $u' \cRf v' \SPE u \cRf v$, tandis que
la concaténation de~«81» et~«83» donne $v' \cRf
u' \spe(p) v \cRf u$. La démonstration est donc
complète.\FinDem

\bigskip Nous pouvons maintenant compléter la
démonstration de la Proposition~«\CritIII». 

\Dem Soit $f$ une fonction de redressement sur~$\Sigma$,
$X$ une partie de~$\Mo(\Sigma) \cup \{\bot\}$ close
par~$\cRf$, et telle que $f$ ait la propriété du cube sur~$X$.
D'après la proposition~«\CohComp», il suffit que nous
montrions que la conjonction de $u' \eqfp u$ et $v' \eqfp v$
entraîne $u' \cRf v' \eqfp u \cRf v$ et $v' \cRf u' \eqfp v
\cRf u$. Comme $\eqfp$ est une relation transitive, il suffit
de prouver l'implication dans le cas où on a utilisé exactement
une relation de~$\Rf$ pour transformer~$uv$ en~$u'v'$. 
Supposons donc, sans perte de généralité puisque la
conclusion cherchée est symétrique, $u' = u$ et que $v'$ soit
obtenu à partir de~$v$ en utilisant une relation de la
présentation, c'est-à-dire en remplaçant un facteur~$a   f(a  ,
b )$ par le facteur~$b  f(b , a  )$ correspondant. On a donc
des décompositions $v = v_0 a   f(a  , b ) v_3$ et $v' = v_0
b  f(b , a  ) v_3$. Soient $p$, $q$, $r$ les longueurs des
mots~$u$, $v_0$ et $v_3$. Par hypothèse, $\Sigma$ est inclus
dans~$X$, donc $u' = u$ entraîne $u' \spe(p) u$, et, de
même, on a $v_0 \spe(q) v_0$ et $v_3 \spe(r) v_3$, et, par
conséquent, on a $v' \SPE v$ par définition. Par le
lemme~«\FineII», on déduit $u \cRf v' \SPE u \cRf v'$ et
$v' \cRf u \spe(p) v \cRf u$, et, de là, en particulier, $u \cRf v'
\eqfp u \cRf v'$ et $v' \cRf u \eqfp v \cRf u$, ce qui
est exactement ce que nous voulions démontrer.\FinDem 

\Rem
Au lieu d'utiliser la relation~$\eqfpp$, et ses
raffinements~$\spe(p)$ et~$\SPE$, nous pourrions espérer
n'utiliser partout que la relation~$\eqfp$, ce qui
simplifierait la démonstration, notamment
parce que $\eqfp$ est partout transitive. Il est douteux
qu'une telle approche naïve puisse aboutir, car la
contrepartie du lemme~«\CompatI» n'a aucune raison
d'être vraie: {\it a priori}, l'hypothèse $u' \eqfp u$
n'implique aucune conséquence directe pour les mots~$u
\cRf u'$ et $u' \cRf u$. 

\Question
{\sl Soit $f$ une fonction de redressement 
sur~$\Sigma$, et $X$ une partie de~$\Mo(\Sigma) \cup
\{\bot\}$ close par~$\cRf$. Que peut-on déduire de
l'hypothèse que $f$ a la propriété du cube faible sur~$X$ ---
en particulier dans le cas $X  = \Mo(\Sigma)$?}

\bigskip Cette question, qui est ouverte, est de peu d'intérêt
algorithmique car la seule méthode systématique connue
pour établir la propriété du cube faible~«\wCOH» est d'établir
la propriété du cube forte~«\COH».

En rapprochant les propositions~«\Necess», «\Suffis», et
«\CritIII», nous obtenons une version explicite du
théorème~B' de l'introduction:

\Prop\«\CritFinal»
{\sl Pour~$f$ fonction de redressement sur~$\Sigma$,
notons $(\CCII)$ la conjonction des conditions
suivantes:

(i) la clôture~$\Til\Sigma$ de~$\Sigma$ pour~$\cRf$ existe et
elle est finie;

(ii) la fonction~$f$ a la propriété du cube sur~$\Til\Sigma$;

(iii) notant~$\Til\Sigma\vv$ la clôture de~$\Til\Sigma$
par~$\jRf$ --- qui nécessairement existe et est finie si
(i) et (ii) sont vérifiées --- il existe un mot~$\O$
dans~$\Til\Sigma\vv$ tel que, pour tout~$u$
dans~$\Til\Sigma\vv$, on ait $\O \cRf u = \e$ et il
existe~$v$ dans~$\Til\Sigma\vv$ vérifiant
$v \cRf \O \eqfpp u$.

Si $M$ est un monoïde de Garside, et si
$f$ est un sélecteur de ppcm sur une partie
génératrice~$\Sigma$ de~$M$,  alors $f$ est une fonction de 
redressement sur~$\Sigma$, elle satisfait aux
conditions~$(\CCII)$ et $M$ admet la
présentation~$\Mon{\Sigma; \Rf}$. Inversement, si $f$ est une
fonction de redressement sur~$\Sigma$ satisfaisant aux
conditions~$(\CCII)$, alors le monoïde~$\Mon{\Sigma; \Rf}$ est
un monoïde de Garside, et $f$ est un sélecteur de ppcm
sur~$\Sigma$.}

\bigskip Les conditions mises en jeu dans la proposition
précédente sont de type~$\S^0_1$: si le monoïde~$\Mon{\Sigma \;
\Rf}$ est un monoïde de Garside, alors les conditions~(i), (ii)
et~(iii) sont vérifiées, donc, partant de l'ensemble fini~$\Sigma$, on
déterminera en un nombre fini d'étapes la clôture de~$\Sigma$
pour~$\cRf$, puis on établira la propriété du cube pour~$f$
sur cette clôture, et, enfin, la condition~(iii), toujours par un
nombre fini de redressements qui, par hypothèse, convergent
tous en un nombre fini d'étapes.

\Ex
Considérons la présentation $\Mon{a  , b  \; a  b a   = b^2}$,
associée au complément~$f$ défini par $f(a  , b ) = b a  $ et
$f(b , a  ) = b $. La clôture de~$\{a  , b \}$ par~$\cRf$
est l'ensemble $\{a  , b , b a  , b a  b , \e, a  b \}$. Que $f$ ait
la propriété du cube sur cet ensemble se vérifie
directement --- en fait, on peut aussi appliquer ici le critère de
la proposition~«\CritII» car l'alphabet est réduit à deux lettres.
Enfin, la clôture de~$\{a  , b , b a  , b a  b , \e, a  b \}$
par~$\jRf$ est l'ensemble
$\{a  , b , b a  , b a  b , \e, a  b , a  b a  ,
b b , a  b a  b , b a  b a  \}$, et il est immédiat de
vérifier que les choix$\O = a  b a  b $ et $\O =
b a  b a  $ conviennent pour la condition~(iii). (Noter que
le groupe $\Gr{a  , b  \; a  b a   =
b^2}$ est le groupe de tresses~$B_3$ rapporté aux
générateurs $a = \s_1$ et $b =
\s_2\s_1$.)

\bigskip Il est alors facile de déduire le théorème~B de l'introduction.
Partant d'une présentation supposée finie  d'un
groupe~$G$ (\resp d'un monoïde~$M$), nous
pouvons énumérer systématiquement toutes les
présentations de~$G$ (\resp de~$M$) en appliquant
les transformations de Tietze, et, pour chacune d'elles,
tester si elle est complémentée, c'est-à-dire associée à une
fonction de redressement, et si cette dernière satisfait aux
conditions de la proposition~«\CritFinal». Alors $G$ est un
groupe de Garside (\resp $M$ est un monoïde de Garside) si et seulement si la réponse est positive
pour au moins une des présentationss, ce qui se
trouvera établi en un nombre fini d'étapes pour autant
que les tests des diverses présentations soient menés
en parallèle, c'est-à-dire sans attendre l'hypothétique
fin d'un test pour passer au suivant.

\Rem Un groupe de Garside peut
s'exprimer comme groupe de fractions de plusieurs
monoïdes: par définition, au moins un de ceux-ci est
un monoïde de Garside, mais ce n'est pas
nécessairement le cas des autres. Par exemple, le
groupe de tresses~$B_3$ est à la fois le groupe de
fractions des monoïdes $\Mon{a  , b  \;
a  b a   =
b  a   b }$ et $\Mon{a  , b  \; a  b a   = b^2}$, qui sont
petits gaussiens, et du monoïde $\Mon{a  , b  \; a  b a   =
b  a  a   b }$, qui n'est pas atomique. De la même
façon, un monoïde de Garside peut admettre
plusieurs présentations associées à des compléments
différents, et les conditions de la
proposition~«\CritFinal» ne valent que pour ceux de ces
compléments qui sont des sélecteurs de ppcm: ainsi,
le monoïde de Garside $\Mon{a, b, c; ab = bc =
ca}$ admet aussi la présentation $\Mon{a, b, c; ab =
bc = ca, aba = caa}$, pour laquelle le redressement n'est
pas toujours convergent. Il serait donc incorrect
d'énoncer la proposition~«\CritFinal»  sous la forme
\<<Le monoïde~$\Mon{\Sigma; \Rf}$ est un monoïde de Garside si et seulement si  $f$ satisfait aux
conditions~$(\CCII)$\>>.

\penalty-2000\bigskip \centerline{{\sc Références}}\bigskip
\parindent=20pt

\Ref \\\Adj; S.I. Adjan;  On the embeddability of
semigroups ; Soviet. Math. Dokl., 1-4; 1960;
819--820.

\Reff \\\Adk; ---; Fragments of the word Delta in a
braid group; Mat. Zam. Acad. Sci. SSSR; 36-1; 1984;
25--34; traduction: Math. Notes of the Acad. Sci. USSR;
36-1 (1984) 505--510.

\Ref \\\BaM; G. Baumslag \& C.F. Miller III (eds); {\rm Algorithms
and Classification in Combinatorial Group Theory}; MSRI
Publications 23, Springer Verlag (1992).

\Ref\\\BDM; D. Bessis, F. Digne, \& J. Michel;
Springer theory in braid groups and the
Birman-Ko-Lee monoid; prépublication (2000).

\Reff \\\BrS; E. Brieskorn \& K. Saito; Artin-Gruppen und
Coxeter-Gruppen; Invent. Math.; 17; 1972; 245--271.

\Reff \\\BMR; M. Brou\'e, G. Malle \& R. Rouquier; Complex
reflection groups, braid groups, Hecke algebras; J. Reine
Angew. Math.; 500; 1998; 127--190.

\Reff\\\BuM; M. Burger \& S. Mozes; Finitely presented
simple groups and product of trees; C. R. Acad. Sci. Paris;
324--1; 1997; 747--752.

\Reff \\\Cha; R. Charney; Artin groups of finite type are
biautomatic; Math. Ann.; 292-4; 1992; 671--683.

\Reff \\\Chb; ---; Geodesic automation and growth
functions for Artin groups of finite type; Math. Ann.;
301-2; 1995; 307--324.

\Ref \\\ClP; A.H. Clifford \& G.B. Preston; {\rm The algebraic
theory of semigroups, vol.~1}; AMS Surveys {\bf 7}, (1961).
  
\Reff \\\Dfa; P. Dehornoy; Deux propri\'et\'es des
groupes de tresses; C. R. Acad. Sci. Paris; 315; 1992;
633--638.

\Reff \\\Dfb; ---; Braid groups and left
distributive operations; Trans. Amer. Math. Soc.; 345-1;
1994; 115--151. 

\Reff \\\Dff; ---; Groups with a complemented
presentation; J. Pure Appl. Algebra; 116; 1997; 115--137.

\Reff \\\Dfz; ---; Gaussian groups are torsion free;
J. of Algebra; 210; 1998; 291--297.

\Reff \\\Dgc; ---; On completeness of word
reversing; Discrete Math.; 225; 2000; 93--119.

\Ref \\\Dgd; ---; {\rm Braids and Self-Distributivity}; Progress 
in Math. vol.~192, Birkh\"auser (2000).

\Reff \\\Dfx; P. Dehornoy \& L. Paris; Garside groups, a
generalization of Artin groups; Proc. London Math. Soc.;
79-3; 1999; 569--604.

\Reff \\\Dlg; P. Deligne; Les immeubles des groupes de
tresses g\'en\'eralis\'es; Invent. Math.; 17; 1972;
273--302.

\Reff \\\ElM; E. A. Elrifai \& H. R. Morton; Algorithms for
positive braids; Quart. J. Math. Oxford; 45-2; 1994;
479--497.

\Ref \\\Eps; D. Epstein \& {\it al.}; Word Processing in
Groups; Jones \& Barlett Publ. (1992).

\Reff \\\Gar; F. A. Garside; The braid group and
other groups; Quart. J. Math. Oxford; 20 {\rm No.78}; 1969;
235--254.

\Reff \\\Mic; J. Michel; A note on words in braid monoids; J.
of Algebra; 215; 1999; 366--377.

\Ref\\\Pic; M. Picantin; The conjugacy problem in small
Gaussian groups; Comm. in Algebra, à paraître.

\Ref\\\Pid; ---; The center of small Gaussian groups; J. of
Algebra, à paraître.

\Reff \\\Rmm;  J.H. Remmers; On the geometry of semigroup
presentations; Advances in Math.; 36; 1980; 283--296.

\Reff \\\Tat; K.~Tatsuoka; An isoperimetric inequality for Artin
groups of finite type; Trans. Amer. Math. Soc.; 339--2;
1993; 537--551. 

\Ref \\\Thu; W.~Thurston; Finite state algorithms for the
braid group; notes en circulation (1988).

\bigskip\hfill SDAD FRE 2271 CNRS, Math\'ematiques

\hfill Universit\'e de Caen, 14~032 Caen, France

\hfill dehornoy@math.unicaen.fr

\hfill http://www.math.unicaen.fr/$\sim$dehornoy/

\bye